\documentclass[10pt]{article}

\usepackage{amsmath,amsfonts,amssymb,cases,mathdots,color,colordvi, url}

\definecolor{lred}{rgb}{1,0.8,0.8}
\definecolor{lblue}{rgb}{0.8,0.8,1}
\definecolor{dred}{rgb}{0.6,0,0}
\definecolor{dblue}{rgb}{0,0,0.5}
\definecolor{red}{rgb}{0.9,0,0}
\definecolor{blue}{rgb}{0,0,0.9}

\newtheorem{proposition}{Proposition}[section]
\newtheorem{theorem}{Theorem}[section]
\newtheorem{lemma}{Lemma}[section]

\newtheorem{definition}{Definition}[section]
\newtheorem{remark}{Remark}

\def\V{{\mathbb{V}}}
\def\R{{\mathbb{R}}}
\def\C{{\mathbb{C}}}
\def\S{{\mathbb{S}}}
\def\T{{\mathbb{T}}}
\renewcommand{\Re}{{\mathbb{R}}}

\title{Spectral Operators of Matrices}

\author{Chao Ding\footnote{National Center for Mathematics and Interdisciplinary Sciences,  Chinese Academy of Sciences, Beijing,  China. This work was initiated while C. Ding was with Department of Mathematics,   National University of Singapore during 2007 to 2012.
  Email: dingchao@amss.ac.cn.}\,,  \ \ Defeng Sun\footnote{Department of Mathematics and Risk Management
Institute,   National University of Singapore,  
Singapore.    Email:
matsundf@nus.edu.sg.   }\,,\ \ Jie Sun\footnote{Department of Mathematics and Statistics, Curtin University, Australia. Email: sun.curtin@gmail.com.}\ \ and \ Kim-Chuan Toh\footnote{Department of
Mathematics,   National University of Singapore,   
Singapore.   Email: mattohkc@nus.edu.sg.}}

\date{January 10, 2014}

\begin{document}
\maketitle

\begin{abstract}

The {class of matrix optimization problems (MOPs) has been recognized in 
recent years to be a powerful tool by} researchers far beyond the optimization community to model many important applications involving  structured low rank matrices. This trend
can be credited to some extent to the exciting developments in the emerging field of compressed sensing. The L\"owner operator, which generates a matrix valued function by applying a single-variable function to each of the singular values of a matrix,
has played an important role for a long time in solving matrix optimization problems.
{However, the classical theory developed for L\"owner operators
 has become inadequate in these recent applications.}  
The main objective of this paper is to provide some necessary theoretical foundations for designing numerical methods for solving the MOP.  This goal is achieved by introducing and conducting a thorough study on a new class of matrix valued functions, coined as spectral operators of matrices.
Several fundamental properties of spectral operators, including the well-definedness, continuity,  directional differentiability,   Fr\'{e}chet-differentiability,  locally Lipschitzian continuity,   $\rho$-order B(ouligand)-differentiability ($0<\rho\leq 1$),  $\rho$-order G-semismooth ($0<\rho\leq 1$) and the characterization of Clarke's generalized Jacobian,  are systematically studied.

\end{abstract}

\vskip 5 true pt

{\bf AMS subject classifications}: 90C25, 90C06, 65K05, 49J50, 49J52

{\bf OR/MS subject classifications}: Primary: Mathematics/matrices; Secondary: Mathematics/functions

\textbf{Key Words}:  directional differentiability; Fr\'{e}chet differentiability;  matrix valued functions; proximal mappings; semismoothness; spectral operators

\section{Introduction}
\label{section:introduction}

 Let $\R^{m\times n}$ and $\C^{m\times n}$ be the vector spaces of  $m\times n$ real and complex matrices over the scalar field of real numbers $\R$, respectively.
 For any   $X\in\C^{m\times n}$,  we denote the conjugate  transpose of $X$ by
$X^\T $.  If $X\in \R^{m\times n}$, then  $X^\T $ is just the transpose of $X$.
 We use  $\V^{m\times n}$ to represent either  the real Euclidean vector space    $\R^{m\times n}$ or $\C^{m\times n}$  with the trace  inner product $\langle X,Y\rangle:={\rm Re}({\rm trace}(X^\T Y))$ for $X,Y\in\V^{m\times n}$ and its induced norm $\|\cdot\|$, where ``${\rm Re}$" means the real part of a complex number.
Without loss of generality, we assume that $m\leq n$ throughout this paper.
For convenience, we also call $\V^{m\times n}$   a matrix space.

 Let $\S^{m}\subseteq \V^{m\times m} $ be the real vector subspace of
$m\times m$ real symmetric matrices or complex Hermitian matrices. For any given $Y \in \S^m$,  we use $\lambda_{1}(Y)\ge \lambda_2(Y) \ge \ldots \ge \lambda_{m}(Y)$ to denote the eigenvalues of $Y$ (all real and counting multiplicity) and use $\lambda(Y)$ to denote the vector of eigenvalues of $Y$. For any given $Z\in\V^{m\times n}$, we use $\sigma_{1}(Z)\geq\sigma_{2}(Z)\geq\ldots\geq\sigma_{m}(Z)$ to denote the singular values of $Z$ (always nonnegative and counting multiplicity)  and use $\sigma(Z)$ to denote the vector of the singular values of $Z$. We use ${\mathbb O}^{p}$ ($p=m,n$) to denote the set of  $p\times p$ orthogonal matrices  in  $ \R^{p\times p}$ if
 $\V^{m\times n} =\R^{m\times n}$ and the set of  $p\times p$ unitary matrices in  $\C^{p\times p}$ if $\V^{m\times n}= \C^{m\times n}$. For  $X\in\V^{m\times m}$, ${\rm diag}(X)$ denotes the column vector   consisting of  all the diagonal entries of $X$ being arranged from the first to the last and  for  $x\in\Re^m$, ${\rm Diag}(x)$  denotes  the $m$ by $m$ diagonal matrix whose $i$-th diagonal entry is $x_i$, $i=1,\ldots,m$.

In this paper, we  shall introduce and  study a class of matrix valued functions, to be called  spectral operators of matrices.
This class of matrix valued functions frequently arise in various applications.  Our first motivating application comes from   matrix optimization problems (MOPs). Let ${\cal X}$ be the vector space $\V^{m\times n}$ or $\S^n$.
Suppose that $f:{\cal X}\to (-\infty,\infty]$ is a closed proper convex function. One simple class of
 MOPs just mentioned take the form of
\begin{equation}\label{eq:def-MOP-primal}
\begin{array}{cl} \min &  \displaystyle f_0(X) + f(X)
\\[1.2mm]
{\rm s.t.} & {\cal A} X = b,\quad X\in{\cal X},
\end{array}
\end{equation}
where $f_0: {\cal X}\to \Re$ is a smooth function whose gradient is Lipschitz continuous, e.g., a linear function $f_0(\cdot)= \langle C, \cdot\rangle$ for some $C\in{\cal X}$,
 ${\cal A}:{\cal X}\to \Re^{p}$ is a linear operator, and $b\in\Re^{p}$ is a given vector.
The above MOPs cover many problems as special cases.  For example,  by considering the particular case that $f\equiv \delta_{\S^{m}_{+}}$, the indicator function of the positive semidefinite matrix cone $\S^{m}_{+}$,
we  can see that the extensively  studied  semidefinite programming (SDP) \cite{Todd01}  is in the form of (\ref{eq:def-MOP-primal}) with a linear function $f_0$.
MOPs also arise frequently from  other applications such as the  matrix norm approximation, matrix completion, rank minimization, graph theory, machine learning, etc  \cite{GT94,Toh97,TT98,RFP07,CanRec08,CanTao09,CLMW09,CSPW09,WMGR09,CFP03,GSun10,Lovasz79,Dobrynin04,KLVempala97}. See \cite{Ding12} for more details.

The Karush-Kuhn-Tucker (KKT) condition of  \eqref{eq:def-MOP-primal} can be written in the following form \cite[Corollary 28.3.1]{Rockafellar70}:
\begin{equation}\label{eq:KKT-MOP-1}
\left\{\begin{array}{l}
\nabla f_0(X)-{\cal A}^*y+\Gamma=0\,,\\
{\cal A}X-b=0\,,\\
\Gamma \in \partial f(X).
\end{array}\right.
\end{equation}
 Let $\psi_{f}:{\cal X}\to \Re$ be the Moreau-Yosida regularization of the closed proper convex function $f$, i.e.,
\begin{equation}\label{eq:def-MY}
\psi_{f}(X):=\min_{Y\in{\cal X}}\Big\{ f(Y)+\frac{1}{2}\|Y-X\|^{2}\Big\},\quad X\in{\cal X}\,,
\end{equation} and $P_{f}(X)$ be the proximal mapping of $f$ at $X$, the unique optimal solution to (\ref{eq:def-MY}).
 It is well-known (see e.g., \cite[Proposition 12.19]{RWets98}) that the mapping $P_{f}:{\cal X}\to{\cal X}$ is globally Lipschitz continuous on ${\cal X}$ with modulus $1$ and $\psi_{f}$ is continuously differentiable on ${\cal X}$ with
$\nabla\psi_{f}(X)=X-P_{f}(X).$
From \cite{Moreau65} (see also \cite[Theorem 31.5]{Rockafellar70})  we know that the KKT condition \eqref{eq:KKT-MOP-1} is equivalent to the following system of Lipschitzian equations
\[
\left[\begin{array}{c}
\nabla f_0(X)-{\cal A}^*y+\Gamma\\
{\cal A}X-b\\
X-P_{f}(X+\Gamma)
\end{array}\right]=0\,.
\]
Thus, the study of MOPs depends crucially on the study of various differential properties of $P_{f}$.
In \cite{ZSToh10,CLSToh12,LSToh12}, Newton-CG based proximal-point algorithms
have been designed to solve large scale SDPs, matrix spectral norm approximation,
and nuclear norm minimization problems, respectively. Those algorithms and
their convergence analyzes all depend crucially on understanding the various
differential properties of the associated proximal mappings $P_f$.


For any given $Z\in\V^{m\times n}$, let  ${\mathbb O}^{m,n}(Z)$   denote the set of  matrix pairs $(U,V)\in{\mathbb O}^{m}\times{\mathbb O}^{n}$ satisfying the singular value decomposition $Z=U\left[ \Sigma(Z)\quad  0 \right]V^\T $, where  $\Sigma(Z)$ is an $m\times m$  diagonal matrix whose $i$-th diagonal entry is $\sigma_i(Z)\ge 0$. For any given $Y\in\S^m$, we use ${\mathbb O}^{m}(Y)$ to denote the set of matrices $P\in{\mathbb O}^{m}$ satisfying the eigenvalue decomposition $Y=P \Lambda(Y) P^\T $, where $\Lambda(Y)$ is an $m\times m$  diagonal matrix whose $i$-th diagonal entry is $\lambda_i(Y)$, a real number.
Assume that the  closed proper convex function $f:{\cal X}\to(-\infty,\infty]$ is unitarily invariant, i.e., for any $X\in{\cal X}\equiv\V^{m\times n}$, $U\in{\mathbb O}^{m}$ and $V\in{\mathbb O}^{n}$, $f(X)=f(U^\T XV)$, or for any $X\in{\cal X}\equiv\S^m$, $P\in{\mathbb O}^{m}$, $f(X)=f(P^\T XP)$. For example, for a given  $k\in\{1,\ldots,m\}$,  Ky Fan's $k$-norm matrix function
$\|X\|_{(k)}={\sum_{i=1}^k}\sigma_i(X)$, $X\in{\mathbb{V}}^{m\times n}$ and the indicator function $\delta_{\S^{m}_+}$ are unitarily invariant. Recall that a function $\phi:\Re^{m}\to(-\infty,\infty]$ is said to be {\it symmetric}  if $\phi(x)=\phi(Qx)$ $\forall$  $x\in\Re^{m}$ and any { permutation matrix} $Q$, and is said to be {\it absolutely symmetric}  if $\phi(x)=\phi(Qx)$ $\forall$  $x\in\Re^{m}$ and any {\it signed permutation matrix} $Q$, which has exactly one nonzero entry in each row and each column, and that entry being $\pm 1$. For the unitarily invariant function $f:{\cal X}\to(-\infty,\infty]$, we know from Lewis \cite[Proposition 2.1]{Lewis95} and Davis \cite{Davis57} that there exists an absolutely symmetric function  $\theta:\Re^m\to (-\infty, +\infty]$ such that $f(\cdot) \equiv   \theta (\sigma(\cdot)) $ if ${\cal X}\equiv\V^{m\times n}$    and  a symmetric function $\theta:\Re^m\to (-\infty, +\infty]$ such that $f(\cdot) \equiv   \theta(\lambda(\cdot)) $ if ${\cal X}\equiv\S^m$, respectively.
 {}Furthermore, from
   \cite{Lewis96,Lewis95,LSendov05a}, we know that the proximal mapping
$P_{f}:{\cal X}\to{\cal X}$ can be written as
\[
P_{f}(X)=\left\{  \begin{array}{ll}
U\left[{\rm Diag}\big(P_{\theta}(\sigma(X))\big)\quad 0\right]V^\T  & \mbox{if $X\in{\cal X}\equiv\V^{m\times n}$}, \\[3pt]
P{\rm Diag}\big(P_{\theta}(\lambda(X))\big)P^\T  & \mbox{if $X\in{\cal X}\equiv\S^m$}\,,
\end{array} \right.
\]
where $(U,V)\in{\mathbb O}^{m,n}(X)$ if $X\in\V^{m\times n}$ and $P\in{\mathbb O}^m(X)$ if $X\in\S^m$.  The proximal mapping $P_{f}$ is a spectral operator with respect to the mixed symmetric mapping $P_{\theta}: \R^m \rightarrow \R^m$ (see Definition \ref{def:mixed_symmetric} in Section \ref{subsec:nonsymmetric}).

Proximal mappings of unitarily invariant proper closed convex functions
belong to a class of matrix functions studied previously in two seminal papers by Lewis \cite{Lewis96}, and Lewis and Sendov \cite{LSendov03}.
In  \cite{Lewis96}, Lewis defined a  Hermitian matrix valued function  by using the gradient mapping $g(\cdot)=\nabla \phi(\cdot):\Re^{m}\to\Re^{m}$ of a symmetric function $\phi:\Re^{m}\to(-\infty,\infty]$.
The corresponding Hermitian matrix valued function $G:\S^{m}\to\S^{m}$  is defined by $G(Y)={\sum_{i=1}^m} g_{i}(\lambda(Y)) p_ip_i^\T $,
where $P\in{\mathbb O}^m(Y)$ and $p_i$ is the $i$-th column of $P$.
Lewis \cite{Lewis96} proved that such a   function $G$ is well-defined, by using the ``block-refineness'' property of $g$.  It was further shown by Lewis and Sendov  in \cite{LSendov03} that $G$ is (continuously) differentiable at $X$ if and only if $g$ is (continuously) differentiable at $\lambda(X)$. Qi and Yang \cite{QYang03} proved that the locally Lipschitz continuous function $G$ is (strongly) semismooth at $X$ if and only if $g$ is (strongly) semismooth at $\lambda(X)$. Note that if  the function $g$ has the form $g(y)=(h(y_1),\ldots,h(y_m))$  $\forall\;y\in\Re^m$ for some given real valued functional $h:\Re\to\Re$, then the corresponding Hermitian matrix valued function $G$ is called   L\"owner's (Hermitian) operator \cite{Lowner34}, which has been well-studied in the literature.  See e.g., \cite{CQTseng03,SSun08} for more details.
For the non-Hermitian case, by considering the gradient mapping $g(\cdot)=\nabla \phi(\cdot):\Re^{m}\to\Re^{m} $ of an  absolutely symmetric function $\phi:\Re^{m}\to(-\infty,\infty]$,  Lewis \cite{Lewis95} studied the corresponding matrix valued function by $G(Z)={\sum_{i=1}^m} g_{i}(\sigma(Z)) u_iv_i^\T $ for $Z\in\V^{m\times n}$,
where $(U,V)\in{\mathbb O}^{m,n}(Z)$ and $u_i$ and $v_i$ are the $i$-th column of $U$ and $V$, respectively. See also Lewis and Sendov \cite{LSendov05a} for more details. If the function $g$ has the form $g(z)=(h(z_1),\ldots,h(z_m))^T\; \forall\; z\in\Re^m$  for some given real valued functional $h:\Re\to\Re$ satisfying $h(0)=0$, then the corresponding non-Hermitian matrix valued function $G$ is called L\"owner's (non-Hermitian) operator in Yang's thesis \cite{Yang09}. Some important properties of L\"owner's (non-Hermitian) operators have been studied by Yang in \cite{Yang09}, including the well-definedness,  (continuous) differentiability and   (strongly) semismoothness.

Besides MOPs, the proximal mapping $P_f$  has also played a   crucial role in some recent applications. For instance, by employing its differentiability and the corresponding derivative formulas, one can derive the divergence expression of the proximal mapping $P_f$, which can be used to obtain the Stein unbiased risk estimator (SURE) \cite{Stein81} of the regularization solution of the matrix recovery problem involving the nuclear norm regularization. Cand\'{e}s et al. \cite{CSTrzasko12} provided a parameter selection method based on the SURE for the singular value thresholding (SVT) operator. See also \cite{DVPFDossal12} for more details.  Although some partial work has been done on different cases,  many fundamental problems of the proximal mapping are unsolved. For example, even in the Hermitian case,  there still remain two important issues  to be explored (i) an explicit formula for the directional derivative of $G$ and (ii) the characterization of Clarke's generalized Jacobian of the general Hermitian matrix
valued function $G$.

The spectral operators of matrices to be considered in this paper  go much beyond proximal mappings. As a matter of fact,
 the spectral operators arising from applications may  {not even}  be the gradient mapping of any scalar valued (absolutely) symmetric 
function.
Therefore, the theoretical results on the spectral operators obtained in this paper are    not covered by the previous works just mentioned \cite{Lewis96,LSendov03,QYang03,Yang09}. For instance, such
spectral operators have already been used  in
low-rank matrix completion problems with fixed basis
coefficients \cite{MSPan12}. The problem of low-rank matrix completion
aims to recover an unknown low-rank matrix from some under-sampled observations with or without noises. A basic approach to solve a low-rank matrix completion problem is to minimize the rank of a matrix subject to certain constraints consistent with the sampled
observations. Since  minimizing  a  rank function with constraints  is generally NP-hard, a widely-used convex relaxation approach is to replace the rank function with the nuclear norm.
For various theoretical breakthroughs along this direction, we refer the readers to  \cite{CanRec08,CanTao09,Gro09,KMO09,Rec09,RFP07} and references therein.
However,  since for many situations, such as the correlation matrix completion in statistics and the density matrix completion in the quantum state tomography where the nuclear norm is a constant, the efficacy of the nuclear norm minimization approach for recovery is inadequate,   Miao et al. \cite{MSPan12} proposed a rank-corrected procedure to generate an estimator of high accuracy and low rank, in which
non-traditional
spectral operators play a pivotal role. A rank-correction term of the form  $-\langle G(\widetilde{X}), X \rangle$ was  added to the nuclear norm penalized least squares model,
where $\widetilde{X}\in\V^{m\times n}$ is a given initial estimator and $G:\V^{m\times n}\to\V^{m\times n}$ is a matrix-valued function defined by
\[
G(X)=U\left[{\rm Diag}\big(g(\sigma(X))\big)\quad 0\right]V^\T ,\quad X\in\V^{m\times n}
\]
with $(U,V)\in{\mathbb O}^{m,n}(X)$, and $g:\Re^{m}\to\Re^m$ is given by
\begin{equation}\label{eq:example-Miao}
g_i(x)=
h\left({\frac{x_i}{\|x\|_{\infty}}}\right) \; \mbox{if $x\in\Re^m\setminus\{0\}$}, \quad
g_i(0) = 0
\end{equation}
for some   scalar valued function $h:\R\rightarrow \R$.
For example, for given $\varepsilon, \tau>0$, the following $h$ was  considered in \cite{MSPan12}:
\begin{equation}\label{eq:h-example-Miao}
h(t)={\rm sgn}(t)(1+\varepsilon^{\tau})\frac{|t|^{\tau}}{|t|^{\tau}+\varepsilon^{\tau}},\quad t\in\Re\, .
\end{equation}
 It can be checked  that  $G$ is the spectral operator with respect to the absolutely symmetric mapping $g$  (Definition \ref{def:mixed_symmetric}). Note that  for such a spectral operator $G$, there does not exist a function $\psi:\V^{m\times n}\to\Re$ such that $G$ is derived through the gradient mapping of $\psi$ because the Jacobian of $G$ at $X$,
when it exists, is not self-adjoint.
By using the rank-correction term, Miao et al. \cite{MSPan12} established a non-asymptotic recovery error result and provided necessary
and sufficient conditions for rank consistency.
Various properties of spectral operators such as well-definedness and continuity play  an important role in their study. More discussions on the rank-correction function can be found in Miao's PhD thesis \cite{Miao13}.

Spectral operators of matrices can  also be used in some other related areas such as in statistical shape analysis, {which involves low rank matrices}. For instance, in order to establish  necessary and sufficient conditions on the existence of the extrinsic mean shape for the reflection shape space (see e.g., \cite{BPatrangenaru03} for the definition) and to provide the corresponding explicit formula (which has important applications in biology, medicine, image analysis, archeology, etc (cf. \cite{DMardia98})), very recently Ding and Qi \cite{DQi13} used  the following matrix valued function {$G:\S^{m}\rightarrow \S^m$ defined by}
\[
G(X) = P{\rm Diag}(g(\lambda(X)))P^{\T}, \quad X\in\S^{m}
\]
with $P\in{\mathbb O}^{m}(X)$, and $g:\Re^{m}\to\Re^{m}$ being given by $g(x) = Q^\T p(x)$,
where for $x\in \Re^m$, $Q$ is an $m$ by $m$ permutation matrix such that $Qx = x^{\downarrow}$,  the vector of entries of $x$ being arranged in the non-increasing order $x^{\downarrow}_{1}\ge\ldots\ge x^{\downarrow}_{m}$ and $p(x)$ is the unique optimal solution to the following convex optimization problem
\[
\min \left\{\frac{1}{2}\|y-x^{\downarrow}\|^{2}\mid \sum_{i=1}^{m}y_{i}=1,\ y_{1}\ge\ldots\ge y_{k}\ge 0, \  y_{k+1}=\ldots=y_{m}=0 \right\}\, ,
\]
where $1\le k \le m$ is a given integer to indicate the  rank of a desired matrix. For a certain nonempty open set ${\cal N}\in\S^{m}$, e.g., ${\cal N}=\left\{X\in{\mathbb{S}}^m\mid \lambda_{k-1}(X)>\lambda_k(X)>\lambda_{k+1}(X) \right\}$, one can easily check that   $g:\Re^m\to \Re^m$ is symmetric (see Definition \ref{def:mixed_symmetric}) on $\lambda_{\cal N}:=\left\{\lambda(X)\mid X\in{\cal N}\right\}$ and the  defined matrix function $G$ is a spectral operator on ${\cal N}$.


The remaining parts of this paper are organized as follows. In
Section \ref{subsec:nonsymmetric}, we give the definition of spectral operators of   matrices and study their well-definedness. Some preliminary results on the differential
properties of singular values and vectors  of matrices are also given in this section.
We study
the continuity,  directional and Fr\'{e}chet-differentiability of spectral operators defined on the single matrix space $\V^{m\times n}$ in Section \ref{section:spectraloperator}. More sophisticated differential properties such as Bouligand-differentiability and $G$-semismoothness of spectral operators are presented in Section \ref{section:semismoothness}.
In Section \ref{section:extension}, we study the spectral operators defined on the Cartesian product of several matrix spaces, and list the main results
corresponding to those derived  in Sections \ref{section:spectraloperator} and
\ref{section:semismoothness}. 
We conclude our paper in the final section.

\section{Spectral operators of matrices}
\label{subsec:nonsymmetric}

In this section, we will first define spectral operators on the Cartesian product
of several  real or complex matrix spaces. The study of spectral operators under this general setting is not only useful but also necessary. In fact, spectral operators defined on the Cartesian product of several matrix spaces  appear naturally in the study of the differentiability of spectral operators, {even if they are} only defined on a single matrix space (see Section \ref{subsection:Directional differentiability}).   Moreover, the spectral operators used in many applications are defined on the Cartesian product of several matrix spaces. See e.g., \cite{DSToh10,WDSToh11} for more details.

Let $s$ be a positive integer and $0\leq s_{0} \leq s$ be a nonnegative integer. For given positive integers $m_{1},\ldots,m_{s}$ and $n_{s_{0}+1},\ldots,n_{s}$, define the finite dimensional real vector space ${\cal X}$ by
\begin{equation*}\label{eq:def-Xspace}
{\cal X}:=\S^{m_{1}} \times\ldots\times\S^{m_{s_{0}}}\times\V^{m_{s_{0}+1}\times n_{s_{0}+1}} \times\ldots\times\V^{m_{s}\times n_{s}}\,.
\end{equation*}
Without loss of generality, we assume that $m_{k}\leq n_{k}$, $k=s_{0}+1,\ldots,s$. For any $X=(X_1,\ldots,X_s)\in{\cal X}$, we have for $1\le k\le s_0$, $X_k\in\S^{m_k}$  and $s_0+1\le k\le s$, $X_k\in\V^{m_k\times n_k}$. Denote $m_{0}:=\sum_{k=1}^{s_{0}}m_{k}$ and $m=\sum_{k=s_{0}+1}^{s}m_{k}$. For any $X\in{\cal X}$, define $\kappa(X)\in\Re^{m_0+m}$ by
\[
\kappa(X):=\left(\lambda(X_1),\ldots,\lambda(X_{s_0}),\sigma(X_{s_0+1}),\ldots,\sigma(X_s)\right)\,.
\]

Recall that a matrix $Q\in\Re^{p\times p}$ is said to be a {\it signed permutation matrix} if $Q$ has exactly one nonzero entry in each row and each column and that entry being $\pm 1$. Let $\mathbb{P}^{p}$ and $\pm\mathbb{P}^{p}$ be the sets of all $p\times p$ permutation matrices and signed permutation matrices, respectively. For ${\cal X}$, define the set ${\cal P}$ by
\begin{equation*}\label{eq:def-permutation-Q}
{\cal P}:=\left\{\left(Q_{1}, \ldots,Q_{s}\right)\,|\, Q_{k}\in\mathbb{P}^{m_{k}},\ 1\leq k\leq s_{0}\ {\rm and}\ Q_{k}\in\pm\mathbb{P}^{m_{k}},\ s_{0}+1\leq k\leq s\right\}\,.
\end{equation*}
Let ${\bf g}:\Re^{m_0+m}\to\Re^{m_0+m}$ be a given mapping. For any    ${\bf x}=({\bf x}_1,\dots,{\bf x}_s)\in  \Re^{m_{0}+m}$ with ${\bf x}_k\in\Re^{m_{k}}$,   rewrite ${\bf g}({\bf x})\in\Re^{m_{0}+m}$ in the form ${\bf g}({\bf x})=\left({\bf g}_{1}({\bf x}), \ldots,{\bf g}_{s}({\bf x})  \right)$ with  ${\bf g}_{k}({\bf x})\in\Re^{m_{k}}$ for $1\le k\le s$.
\begin{definition}\label{def:mixed_symmetric}
The given mapping ${\bf g}:\Re^{m_0+m}\to\Re^{m_0+m}$  is said to be {\it mixed symmetric}, with respect to ${\cal P}$, at     ${\bf x}=({\bf x}_1,\dots,{\bf x}_s)\in  \Re^{m_{0}+m}$ with ${\bf x}_k\in\Re^{m_{k}}$, if
\begin{equation}\label{eq:def-symmetric}
{\bf g}(Q_1{\bf x}_1,\ldots,Q_s{\bf x}_s)=\left(Q_1{\bf g}_1({\bf x}),\ldots,Q_{s}{\bf g}_{s}({\bf x})\right)\quad \forall\,\left(Q_{1}, \ldots,Q_{s}\right)\in{\cal P}\,.
\end{equation}
The mapping ${\bf g}$ is said to be mixed symmetric,  with respect to ${\cal P}$,  over a set ${\cal D} \subseteq \Re^{m_0+m}$ if \eqref{eq:def-symmetric} holds for every ${\bf x}\in{\cal D}$. We call ${\bf g}$ a {\it mixed symmetric} mapping,  with respect to ${\cal P}$,  if (\ref{eq:def-symmetric}) holds for every ${\bf x}\in \Re^{m_{0}+m}$.
\end{definition}

Note that for each $k\in\{1,\ldots,s\}$, the function value ${\bf g}_k({\bf x})\in\Re^{m_0+m}$ is dependent on all $\bf x_1,\ldots,{\bf x}_s$.  With causing no confusion, in later discussions  we often   drop ``with respect to ${\cal P}$" from    Definition \ref{def:mixed_symmetric}.
The following result on ${\bf g}$ can be checked directly from the definition.

\begin{proposition}\label{prop:smallgprop} Suppose that the mapping ${\bf g}:\Re^{m_0+m}\to\Re^{m_0+m}$  is mixed  symmetric at  ${\bf x}=({\bf x}_1,\dots,{\bf x}_s)\in  \Re^{m_{0}+m}$ with ${\bf x}_k\in\Re^{m_{k}}$.     Then,
 for any $i,j\in\{1,\ldots,m_k\}$,
 \[({\bf g}_k({\bf x}))_i=({\bf g}_k({\bf x}))_j\quad  {\rm  if} \quad  ({\bf x}_k)_{i}=({\bf x}_k)_{j}, \quad  \forall \,  1\le k\le s
\]
and
\[
({\bf g}_k({\bf x}))_i =0 \quad {\rm if} \quad ({\bf x}_{k})_i=0, \quad \forall\,  s_0+1\le k\le s.
\]
\end{proposition}

 Let ${\cal N}$ be a given nonempty set in ${\cal X}$. Define
 \[
 \kappa_{\cal N}: =\left\{ \kappa(X)\mid   X \in {\cal N}\right\}\, .
 \]
\begin{definition}\label{def:def-spectral-op}
Suppose that  ${\bf g}:\Re^{m_0+m}\to\Re^{m_0+m}$  is   mixed  symmetric on $\kappa_{\cal N}$.
The spectral operator $G:{\cal N}\to{\cal X}$ with respect to   ${\bf g}$ is defined by
\[
G(X):=\left(G_{1}(X),\ldots,G_{s}(X)\right),\quad X=(X_1,\ldots,X_s)\in{\cal N}
\] with
\begin{equation*}\label{eq:def-Gk-two-parts}
G_{k}(X):=\left\{  \begin{array}{ll}
P_{k}{\rm Diag}\big({\bf g}_{k}(\kappa(X))\big)P_{k}^\T  & \mbox{if $1\le k\le s_{0}$,}
\\[3pt]
U_{k}\left[{\rm Diag}\big({\bf g}_{k}(\kappa(X))\big)\quad 0\right]V_{k}^\T  & \mbox{if $s_{0}+1\le k\le s$,}
\end{array} \right.
\end{equation*} where $P_{k}\in{\mathbb O}^{m_{k}}(X_{k})$, $1\leq k\leq s_{0}$, $(U_{k},V_{k})\in{\mathbb O}^{m_{k},n_{k}}(X_{k})$, $s_{0}+1\leq k\leq s$.
\end{definition}

Before studying the well-definedness of spectral operators, it is worth   mentioning that for the case that ${\cal X}\equiv\S^m$ (or $\V^{m\times n}$) if ${\bf g}$ has the form ${\bf g}(y)=(h(y_1),\ldots,h(y_m))\in\Re^m$ with $y_i\in\Re$ for some given scalar valued  functional $h:\Re\to\Re$, then the corresponding spectral operator $G$ is    called  L\"{o}wner operator by Sun and Sun \cite{SSun08} in recognitions of L\"{o}wner's original contribution on this topic in  \cite{Lowner34} (or the L\"{o}wner non-Hermitian operator by Yang in her thesis \cite{Yang09} if $h(0) =0$).

\subsection{The well-definedness}

In order to show the well-definedness of spectral operators, we need the following two simple propositions.

Let $\overline{Y}\in\S^m$ be given. Denote $\overline{\mu}_{1} > \overline{\mu}_{2}>\ldots> \overline{\mu}_{r}$ the distinct eigenvalues of $\overline{Y}$. Define the index sets
 \begin{equation*}\label{eq:ak-symmetric}
 \alpha_{l} :=\{i\,|\, \lambda_{i}(\overline{Y})= \overline{\mu}_{l},\ 1\leq i\leq m\},\quad l=1,\ldots,r\,.
 \end{equation*}
  Let  $\Lambda(\overline{Y})$ be the $m\times m$  diagonal matrix whose $i$-th diagonal entry is $\lambda_i(\overline{Y})$. Then, the following elementary property on the eigenvalue decomposition of $\overline{Y}$ can be checked directly.
\begin{proposition}\label{prop:QLambdaQ}
 The matrix $Q\in {\mathbb O}^m$ satisfies $Q\Lambda(\overline{Y})=\Lambda(\overline{Y})Q$ if and only if there exist $Q_l\in{\mathbb Q}^{|\alpha_l|}$, $l=1,\ldots,r$ such that $Q$ is a block diagonal matrix whose $l$-th diagonal block is $Q_l$, i.e.,
\[
Q={\rm Diag}(Q_{1},Q_{2},\ldots,Q_{r})\,.
\]
\end{proposition}

Let $\overline{Z}\in\V^{m\times n}$  be  given. We  use $\overline{\nu}_{1}>\overline{\nu}_{2}>\ldots>\overline{\nu}_{r}>0$ to denote the nonzero distinct singular values of $\overline{Z}$.
Define
\begin{equation}\label{eq:ak-nonsymmetric}
a_l:=\{i\,|\,\sigma_{i}(\overline{Z})=\overline{\nu}_{l}, \ 1\le i\le m\}, \quad l=1,\ldots,r\, \quad {\rm and} \quad  b:= \{i\,|\,\sigma_{i}(\overline{Z})=0, \ 1\le i\le m\}.
\end{equation}
The following observation can be derived easily. For the real case, the proof can be found in \cite[Theorem 3.7]{LSendov05b}, and the corresponding result for the complex case can be obtained similarly.

\begin{proposition}\label{prop:PSigma=SigmaW}
Let  $\overline{\Sigma}:=\Sigma(\overline{Z})$. Then,
$P\in{\mathbb O}^{m}$ and $W\in{\mathbb O}^{n}$ satisfy
\begin{equation*}\label{eq:twodecomp}
P\left[ \overline{\Sigma} \quad 0 \right]=\left[  \overline{\Sigma}  \quad  0 \right]W
\end{equation*} if and only if there exist $Q\in{\mathbb O}^{|a|}$, $Q'\in{\mathbb O}^{|b|}$ and $Q''\in{\mathbb O}^{n-|a|}$ such that
\[
P=\left[\begin{array}{cc}Q & 0 \\0 & Q'\end{array}\right] \quad {\rm and} \quad W=\left[\begin{array}{cc}Q & 0 \\0 & Q''\end{array}\right]\,,
\]
where $|a|= |a_1| +\ldots +|a_r|$ and  $Q={\rm Diag}(Q_{1},Q_{2},\ldots,Q_{r})$ is a   block diagonal  matrix whose $l$-th  diagonal block is $Q_{l}\in{\mathbb O}^{|a_{l}|}$.
\end{proposition}

By combining Propositions \ref{prop:QLambdaQ} and \ref{prop:PSigma=SigmaW} with the mixed symmetric property of ${\bf g}$, we are able to obtain the following result on the well-definedness of spectral operators.
\begin{theorem}\label{thm:well-def-spectral-op}  Let ${\bf g}:\Re^{m_0+m}\to\Re^{m_0+m}$  be mixed  symmetric on $\kappa_{\cal N}$.
Then the spectral operator $G:{\cal N}\to{\cal X}$ defined in Definition \ref{def:def-spectral-op} with respect to   ${\bf g}$ is well-defined.
\end{theorem}
\noindent {\bf  Proof.} Let  $X=(X_1,\ldots,X_s)\in{\cal N}$ be arbitrarily chosen with   $X_k\in\S^{m_k}$  for $1\le k\le s_0$ and $X_k\in\V^{m_k\times n_k}$ for $s_0+1\le k\le s$.
Let ${\bf x}=({\bf x}_{1},\ldots,{\bf x}_{s}) :=\kappa(X)$ with   ${\bf x}_k\in\Re^{m_{k}}$. Then we know from Proposition \ref{prop:smallgprop}  that for any $i,j\in\{1,\ldots,m_k\}$,
 \[({\bf g}_k({\bf x}))_i=({\bf g}_k({\bf x}))_j\quad  {\rm  if} \quad  ({\bf x}_k)_{i}=({\bf x}_k)_{j}, \quad  \forall \,  1\le k\le s
\]
  and
  \[
  ({\bf g}_k({\bf x}))_i =0 \quad {\rm if} \quad ({\bf x}_{k})_i=0, \quad \forall\,  s_0+1\le k\le s,
   \]
   which, together with Propositions \ref{prop:QLambdaQ} and   \ref{prop:PSigma=SigmaW}, imply that
   the  matrix $G(X)$ is independent of the choices of  $P_{k}\in{\mathbb O}^{m_{k}}(X_{k})$, $1\leq k\leq s_{0}$, $(U_{k},V_{k})\in{\mathbb O}^{m_{k},n_{k}}(X_{k})$, $s_{0}+1\leq k\leq s$. That is,     $G$ is well defined at $X$.  Since $X$ is arbitrarily chosen from ${\cal N}$,  the spectral operator $G$ is well-defined on ${\cal N}$. $\hfill \Box$

\medskip

\subsection{Differential properties of singular values and vectors}

In this subsection, we collect  some  useful preliminary results on the
singular value decomposition (SVD) of matrices.
 Let $\overline{Z}\in\V^{m\times n}$  be  given. Consider the following SVD of $\overline{Z}$:
\begin{equation}\label{eq:SVD}
\overline{Z}=\overline{U}\left[ \Sigma(\overline{Z})\quad  0 \right]\overline{V}^\T =\overline{U}\left[ \Sigma(\overline{Z})\quad  0\right]\left[ \overline{V}_{1} \quad  \overline{V}_{2} \right]^\T =\overline{U}\Sigma(\overline{Z}) \overline{V}_{1}^\T \,,
\end{equation} where
 $\overline{U}\in{\mathbb O}^{m}$ and $\overline{V}=\left[ \overline{V}_{1} \quad  \overline{V}_{2} \right] \in{\mathbb O}^{n}$ with
  $\overline{V}_{1}\in\V^{n\times m}$ and
  $\overline{V}_{2}\in\V^{n\times (n-m)}$.
  Define the index sets $a$ and $c$ by
\begin{equation}\label{eq:def-a-b-c}
a:=\{i \,|\,\sigma_{i}(\overline{Z})>0,\ 1\le i\le m\}  \   \  {\rm and} \ c:=\{m+1,\ldots,n\}\,.
\end{equation}
Let the  index sets $a_l$, $l=1,\ldots,r$ and $b$ be  defined by \eqref{eq:ak-nonsymmetric}. For each $  i\in \{1,\ldots,m\}$, we also define $l_{i}(\overline{Z})$ to be the number of singular values which are equal to $\sigma_ i(\overline{Z})$ but are ranked  before $i$ (including $i$), and $\tilde{l}_{i}(\overline{Z})$ to be the number of singular values which are equal to $\sigma_i(\overline{Z})$ but are ranked  after $i$ (excluding $i$),  i.e.,
 define $l_{i}(\overline{Z})$ and  $\tilde{l}_{i}(\overline{Z})$  such that
\begin{eqnarray}
&&\sigma_{1}(\overline{Z})\geq\ldots\geq\sigma_{i-l_{i}(\overline{Z})}(\overline{Z})>\sigma_{i-l_{i}(\overline{Z})+1}(\overline{Z})=\ldots=\sigma_{i}(\overline{Z})=\ldots=\sigma_{i+\tilde{l}_{i}(\overline{Z})}(\overline{Z})\nonumber\\
&&>\sigma_{i+\tilde{l}_{i}(\overline{Z})+1}(\overline{Z})\geq\ldots\geq\sigma_{m}(\overline{Z})\,.\label{eq:onsymmetric-l_i}
\end{eqnarray}
In later discussions, when the dependence of $l_{i}$ and $\tilde{l}_{i}$ on $\overline{Z}$
are clear  from the context, we often drop $\overline{Z}$ from these notations.

For any $Y\in\V^{m\times n}$, let $Y_{ij}$ be the $(i,j)$-th entry of $Y$. For any $Y\in\V^{m\times n}$ and the given index sets ${\cal I}\subseteq \{1,\ldots, m\}$ and  ${\cal J}\subseteq \{1,\ldots, n\}$,  we use $ Y_{{\cal J}}$ to denote the sub-matrix of $Y$ obtained by removing all the columns of $Y$ not in ${\cal J}$ and use $Y_{{\cal I}{\cal J}}$ to
  denote the $|{\cal I}|\times|{\cal J}|$ sub-matrix of $Y$ obtained  by removing all the rows of $Y$ not in ${\cal I}$ and all the columns of $Y$ not in  ${\cal J}$. For notational convenience, we define two linear  matrix operators
$S :\V^{p\times p}\to \S^{p}$, $T :\V^{p\times p}\to\V^{p\times p}$  by
\begin{equation}\label{eq:maps-ST}
S(Y):=\frac{1}{2}(Y+Y^\T ) \quad
T(Y):=\frac{1}{2}(Y-Y^\T ), \quad Y\in\V^{p\times p}\,.
\end{equation}
 The following proposition can be derived directly from the directional differentiability (e.g., see \cite[Theorem 7]{Lancaster64} and \cite[Proposition 1.4]{Torki01}) of the eigenvalues of a Hermitian matrix. For more details, see \cite[Section 5.1]{LSendov05b}.

\begin{proposition}\label{prop:easy}
Suppose that $\overline{Z}\in\V^{m\times n}$ has the SVD (\ref{eq:SVD}).  For any $\V^{m\times n}\ni H\to 0$, we have
\begin{equation}\label{eq:directi-diff}
\sigma_{i}(\overline{Z}+H)-\sigma_{i}(\overline{Z})-\sigma'_{i}(\overline{Z};H)=O(\|H\|^{2})\,,  \quad i=1,\ldots,m\,,
\end{equation} where
\begin{equation}\label{eq:directi-diff-sigma}
\sigma'_{i}(\overline{Z};H)=\left\{ \begin{array}{lcl}
\lambda_{l_{i}}\left( S(\overline{U}_{a_{l}}^\T H \overline{V}_{a_{l}})\right) & {\rm if} & i\in a_{l},\ l=1,\ldots,r\,,\\[3pt]
\sigma_{l_{i}}\Big(\left[\overline{U}_{b}^\T H \overline{V}_{b}\quad \overline{U}_{b}^\T H \overline{V}_{2} \right]\Big) & {\rm if} & i\in b\,,
\end{array}\right.
\end{equation}
where for each $i\in\{1,\ldots,m\}$, $l_{i}$ is defined in (\ref{eq:onsymmetric-l_i}).
\end{proposition}
\vskip 5 true pt

The following results are also needed for subsequent discussions. For the real case, the detailed proof can be found in \cite[Proposition 7]{DSToh10}.   The results for the complex case can be derived in a similar manner.

\begin{proposition}\label{prop:diag-U-V}
 For any $\V^{m\times n}\ni H\to 0$, let $Z:=\left[ \Sigma(\overline{Z})\quad  0 \right]+H$. Suppose that $U\in{\mathbb O}^{m}$ and $V=[V_{1}\quad V_{2}]\in{\mathbb O}^{n}$ with $V_{1}\in\V^{n\times m}$ and $V_{2}\in\V^{n\times(n-m)}$ satisfy
\[
Z=\left[ \Sigma(\overline{Z})\quad  0 \right]+H=U\left[ \Sigma(Z)\quad  0 \right]V^\T =U\left[ \Sigma(Z)\quad  0 \right][V_{1}\quad V_{2}]^\T \,.
\] Then, there exist $Q\in{\mathbb O}^{|a|}$, $Q'\in{\mathbb O}^{|b|}$ and $Q'' \in{\mathbb O}^{n-|a|}$  such that
\begin{equation}\label{eq:diag-U-V}
U=\left[\begin{array}{cc}Q & 0 \\0 & Q'\end{array}\right]+O(\|H\|) \quad {\rm and} \quad V=\left[\begin{array}{cc}Q & 0 \\0 & Q''\end{array}\right]+O(\|H\|)\,,
\end{equation} where $Q={\rm Diag}(Q_{1},Q_{2},\ldots,Q_{r})$,  $Q_{l}\in{\mathbb O}^{|a_{l}|}$. Furthermore, we have
\begin{eqnarray}\label{eq:diag-H-decomp-ak}
\Sigma(Z)_{a_{l}a_{l}}-\Sigma(\overline{Z})_{a_{l}a_{l}} &=&Q_{l}^\T S(H_{a_{l}a_{l}})Q_{l}+O(\|H\|^{2}), \quad l=1,\ldots,r,
\\[3pt]
\label{eq:diag-H-decomp-b}
\left[\Sigma(Z)_{bb}-\Sigma(\overline{Z})_{bb}\quad 0\right] &=&
Q'^\T \left[H_{bb}\quad H_{bc}\right]Q''+O(\|H\|^{2})\,.
\end{eqnarray}
\end{proposition}

Given the index set  $a_l$ for  $l\in\{1,\ldots,r\}$ at $\overline{Z}\in\V^{m\times n}$, we  define ${\cal U}_{l} :\V^{m\times n}\to\V^{m\times n}$ by
 \begin{equation}\label{eq:Pk-nonsymmetric}
 {\cal U}_{l}(Z)=\sum_{i\in a_{l}}u_{i}v_{i}^\T ,\quad  Z \in \V^{m\times n}\,,
 \end{equation}
where $u_{i}$ and $v_{i}$ are the $i$-th column of $U$ and $V$, respectively, and $(U,V)\in {\mathbb O}^{m,n}(Z)$.
Let ${\cal B}\subseteq\V^{m\times n}$ be an open neighborhood of $\overline{Z}$. By shrinking ${\cal B}$ if necessary, we may assume that for any  $Z\in{\cal B}$, if $i\in a_{l}$, $1\leq l \leq r$, then $\sigma_{i}(Z)>0$, and if $i\in a_{l}$, $j\in a_{l}$ and $1\leq l\neq l'\leq r$, then $\sigma_{i}(Z)\neq\sigma_{j}(Z)$. Therefore, for any $Z\in{\cal B}$, we may define matrices
     $\Gamma_{l}(Z)$,  $\Xi_{l}(Z)\in\Re^{m\times m}$  and $\Upsilon_{l}(Z)\in\Re^{m\times(n-m)}$, $l=1,\ldots,r$ by
 \begin{eqnarray}
(\Gamma_{l}(Z))_{ij}&=&\left\{ \begin{array}{ll}
                \displaystyle{\frac{1}{\sigma_{i}(Z)-\sigma_{j}(Z)}} & {\rm if} \ i\in a_{l},\ j\in a_{l'},\ l\neq l',\ l'=1,\ldots,r+1 \,,\\
                 \displaystyle{\frac{-1}{\sigma_{i}(Z)-\sigma_{j}(Z)}} & {\rm if} \ i\in a_{l'},\ j\in a_{l},\ l\neq l',\ l'=1,\ldots,r+1 \,,\\
                 0 & {\rm otherwise} \,,
                \end{array}\right.
\label{eq:Gamma_k}
\\[3pt]
(\Xi_{l}(Z))_{ij}&=&\left\{ \begin{array}{ll}
                \displaystyle{\frac{1}{\sigma_{i}(Z)+\sigma_{j}(Z)}} & {\rm if}\ i\in a_{l},\ j\in a_{l'},\ l\neq l',\ l'=1,\ldots,r+1 \,,\\
                 \displaystyle{\frac{1}{\sigma_{i}(Z)+\sigma_{j}(Z)}} & {\rm if}\ i\in a_{l'},\ j\in a_{l},\ l\neq l',\ l'=1,\ldots,r+1 \,,\\
                 \displaystyle{\frac{2}{\sigma_{i}(Z)+\sigma_{j}(Z)}} & {\rm if}\ i,j\in a_{l},\\
                 0 & {\rm otherwise},
                \end{array}\right.
\label{eq:Xi_k}
\\[3pt]
(\Upsilon_{l}(Z))_{ij}&=&\left\{ \begin{array}{ll}
                 \displaystyle{\frac{1}{\sigma_{i}(Z)}} & {\rm if}\ i\in a_{l},\;
j=1,\ldots,n-m,
             \\
                 0 & {\rm otherwise}.
                \end{array}\right.
\label{eq:tildeXi_k}
\end{eqnarray}
We use $``\circ"$ to denote the usual Hadamard product between two matrices, i.e., for any two matrices $A$ and $B$ in $\V^{m\times n}$ the $(i,j)$-th entry of $  Z:= A\circ B \in \V^{m\times n}$ is
$Z_{ij}=A_{ij} B_{ij}$.
 We have the following differential properties of  ${\cal U}_{l} $, $l=1, \ldots, r$. For the real case, the results have been shown in \cite[Proposition 2.11]{DSToh10}. By using similar arguments to the real case,  one can derive the corresponding results for the complex case.

 \begin{proposition}\label{prop:P''_k}
   Let   ${\cal U}_{l} $, $l=1, \ldots, r$ be
defined by (\ref{eq:Pk-nonsymmetric}).
 Then, there exists an open neighborhood
${\cal B}$ of $\overline{Z}$ such that
${\cal U}_{l} $ is at least twice continuously differentiable in ${\cal
B}$, and for   any
$H\in\V^{m\times n}$, the first order derivative of ${\cal U}_{l} $ at $Z\in{\cal B}$ is given by
\begin{equation}\label{eq:U'_k}
{\cal U}^\prime _{l}(Z)H= U[\Gamma_{l}(Z)\circ S(U^\T HV_{1})+\Xi_{l}(Z)\circ
T(U^\T HV_{1})]V_{1}^\T +U(\Upsilon_{l}(Z)\circ U^\T HV_{2})V_{2}^\T \,,
\end{equation}
 where $(U,V)\in{\mathbb O}^{m,n}(Z)$ and the
linear operators $S $ and $T $ are defined by
(\ref{eq:maps-ST}).
 \end{proposition}


\section{Continuity, directional and Fr\'{e}chet differentiability}
\label{section:spectraloperator}

 In this and the next section, we will first  focus on the study of spectral operators for the case that ${\cal X}\equiv\V^{m\times n}$. The corresponding extensions for  the spectral operators defined on the general Cartesian product  of several matrix spaces will be presented  in Section \ref{section:extension}.

Let ${\cal N}$ be a given nonempty open set in $\V^{m\times n}$. Suppose that $g:\Re^m\to\Re^m$ is  mixed symmetric, with respect to ${\cal P}\equiv\pm{\mathbb P}^m$ (called absolutely symmetric in this case), on an open set $\hat{\sigma}_{ {\cal N}}$ in $\Re^{m}$ containing $\sigma _{\cal N}:=\left\{\sigma(X)\mid X\in{\cal N}\right\}$. The spectral operator $G:{\cal N}\to\V^{m\times n}$ with respect to $g$  defined in Definition \ref{def:def-spectral-op} then takes the form of
\[
G(X)=U\left[{\rm Diag}(g(\sigma(X)))\quad 0\right]V^\T ,\quad X\in{\cal N}\,,
\] where $(U,V)\in{\mathbb O}^{m,n}(X)$. Let  $\overline{X}\in{\cal N}$ be given. Consider the SVD  (\ref{eq:SVD}) for $\overline{X}$, i.e.,
\begin{equation}\label{eq:Y-eig-Z-SVD}
\overline{X}=\overline{U}\left[ \Sigma(\overline{X})\quad  0 \right]\overline{V}^\T \,,
\end{equation} where $\overline{V}=\left[ \overline{V}_{1} \quad  \overline{V}_{2} \right] \in{\mathbb O}^{n}$ with
  $\overline{V}_{1}\in\V^{n\times m}$ and
  $\overline{V}_{2}\in\V^{n\times (n-m)}$. Let $\overline{\sigma}:=\sigma(\overline{X})\in\Re^{m}$. Let $a$, $b$, $c$, $a_{l}$, $l=1,\ldots,r$ be the index sets defined by (\ref{eq:def-a-b-c})  and (\ref{eq:ak-nonsymmetric}) with $\overline{Z}$ being replaced by $\overline{X}$. Denote $\bar{a}:=\{1,\ldots,n\}\setminus a$. For any given vector $y\in\Re^m$, let $|y|^{\downarrow}$ be the vector of entries of $|y|=(|y_1|,\ldots,|y_m|)$ being arranged in the non-increasing order $|y|^{\downarrow}_1\ge\ldots\ge|y|^{\downarrow}_m$. The following result follows from the absolutely symmetric property of $g$ on $\hat{\sigma}_{\cal N}$.

\begin{proposition}\label{prop:G-g}
Let $U\in{\mathbb O}^m$ and $V=\left[ V_{1} \quad  V_{2} \right] \in{\mathbb O}^{n}$ with
  $V_{1}\in\V^{n\times m}$ and
  $V_{2}\in\V^{n\times (n-m)}$ be given. Let $y\in\hat{\sigma}_{\cal N}$.   Then, for $Y:=U\left[{\rm Diag}(y)\quad 0\right]V^\T$ it always holds that
\[
G(Y)=U\left[{\rm Diag}(g(y))\quad 0\right]V^\T= U{\rm Diag}(g(y))V_1^\T \, .
\]
\end{proposition}
\noindent {\bf  Proof.} Let   $P\in\pm{\mathbb P}^m$ be a signed permutation matrix such that $Py=|y|^{\downarrow}$. Then, we know that $\sigma(Y)=|y|^{\downarrow}$ and  $Y$ has the following SVD
\[
Y=U[P^\T{\rm Diag}(|y|^{\downarrow})W \quad 0]V^\T= UP^\T\left[{\rm Diag}(|y|^{\downarrow})\quad 0 \right][V_1W^\T \quad V_2]^\T\,,
\]
where $W:=|P|\in{\mathbb P}^m$ is the $m$ by $m$ permutation matrix whose $(i,j)$-th element is the absolute value of the $(i,j)$-th element of $P$. Then, we know from Definition \ref{def:def-spectral-op} that
\[
G(Y)=UP^\T\left[{\rm Diag}(g(|y|^{\downarrow})) \quad 0\right][V_1W^\T \quad V_2]^\T\,.
\]
Since $g$ is absolutely symmetric at $y$, one has
\[
{\rm Diag}(g(|y|^{\downarrow}))={\rm Diag}(g(Py))={\rm Diag}(Pg(y))=P{\rm Diag}(g(y))W^\T\,.
\]
Thus,
\[
G(Y)=UP^\T\left[P{\rm Diag}(g(y))W^\T \quad 0\right][V_1W^\T \quad V_2]^\T=U\left[{\rm Diag}(g(y))\quad 0\right]V^\T\,,
\]
 which, proves the conclusion. $\hfill \Box$
\vskip 10 true pt

By using Proposition \ref{prop:diag-U-V}, we have the following result on the continuity of the spectral operator $G$.

\begin{theorem}\label{thm:continuity}
Suppose that $\overline{X}\in{\cal N}$ has the SVD (\ref{eq:Y-eig-Z-SVD}). The spectral operator $G$ is continuous at $\overline{X}$ if and only if   $g$ is continuous at $\sigma(\overline{X})$.
\end{theorem}
\noindent {\bf  Proof.} $``\Longleftarrow"$  Let $X\in{\cal N}$. Denote $H=X-\overline{X}$ and $\sigma=\sigma(X)$. Let $U\in{\mathbb{O}^m}$ and $V\in{\mathbb{O}^n}$ be such that $X=\overline{X}+H=U\left[ \Sigma(X)\quad  0 \right]V^\T$.
Then, we know from \eqref{eq:Y-eig-Z-SVD} that
\[
\left[ \Sigma(\overline{X})\quad  0 \right]+\overline{U}^\T H\overline{V}=\overline{U}^\T U\left[ \Sigma(X)\quad  0 \right]V^\T\overline{V}\,.
\]
From \eqref{eq:diag-U-V} in Proposition \ref{prop:diag-U-V}, we know that for any $X$ sufficiently close to $\overline{X}$,  there exist $Q\in{\mathbb O}^{|a|}$, $Q'\in{\mathbb O}^{|b|}$ and $Q'' \in{\mathbb O}^{n-|a|}$  such that
\begin{equation}\label{eq:con-UQ}
\overline{U}^\T U=\left[\begin{array}{cc}Q & 0 \\0 & Q'\end{array}\right]+O(\|H\|) \quad {\rm and} \quad \overline{V}^\T V=\left[\begin{array}{cc}Q & 0 \\0 & Q''\end{array}\right]+O(\|H\|)\,,
\end{equation} where $Q={\rm Diag}(Q_{1},Q_{2},\ldots,Q_{r})$,  $Q_{l}\in{\mathbb O}^{|a_{l}|}$. On the other hand,   from the definition of the spectral operator $G$ one has
\[
U^\T \left(G(X)-G(\overline{X})\right)V=\left[{\rm Diag}(g(\sigma))\quad 0\right]- U^\T \overline{U} \left[{\rm Diag}(g(\overline{\sigma}))\quad 0\right]\overline{V}^\T V\,.
\]
Thus, we obtain from \eqref{eq:con-UQ} and Proposition \ref{prop:smallgprop} that for any  $X$ sufficiently close to $\overline{X}$, 
\[
U^\T \left(G(X)-G(\overline{X})\right)V= \left[{\rm Diag}(g(\sigma)-g(\overline{\sigma}))\quad 0\right]+O(\|H\|)\,.
\]
Thus, since $g$ is assumed to be continuous at $\overline{\sigma}$, we can conclude that the spectral operator $G$ is   continuous at $\overline{X}$.

 $``\Longrightarrow"$ Suppose that $G$ is continuous at $\overline{X}$. Let $(\overline{U},\overline{V})\in{\mathbb O}^{m\times n}(\overline{X})$ be fixed.  Choose any $\sigma\in\hat{\sigma}_{ {\cal N}}$ and denote $X:=\overline{U}[{\rm Diag}(\sigma)\quad 0]\overline{V}^{\T}$. We know from Proposition \ref{prop:G-g} that  $G(X)=\overline{U}{\rm Diag}(g(\sigma))\overline{V}_1^{\T}$ and
\[
 {\rm Diag}(g(\sigma)-g(\overline{\sigma})) = \overline{U}^\T\left(G(X)-G(\overline{X})\right)\overline{V}_1\,.
 \]
Hence, we know from the assumption that $g$ is continuous at $\overline{\sigma}$.
 $\hfill \Box$
 \vskip 10 true pt

Next, we introduce some notations which are frequently used in   later discussions. For any given $X\in{\cal N}$, let $\sigma=\sigma(X)$. For the mapping $g$, we define three matrices ${\cal E}_{1}({\sigma}),{\cal E}_{2}({\sigma})\in\Re^{m\times m}$ and ${\cal F}({\sigma})\in\Re^{m\times(n-m)}$ (depending on $X\in{\cal N}$) by
\begin{eqnarray}
({\cal E}^0_{1}({\sigma}))_{ij} &:=&
\left\{ \begin{array}{ll}  \displaystyle{\frac{g_i(\sigma)-g_j({\sigma})}{\sigma_{i}-\sigma_{j}}} & \mbox{if $\sigma_{i}\neq\sigma_{j}$}\,,\\[3pt]
0 & \mbox{otherwise}\,,
\end{array} \right.  \quad i,j\in\{1,\ldots,m\}\,,
\label{eq:def-matric-E1}
\\[3pt]
({\cal E}^0_{2}({\sigma})_{ij} &:=&
\left\{ \begin{array}{ll}  \displaystyle{\frac{g_i({\sigma})+g_j({\sigma})}{\sigma_{i}+\sigma_{j}}} & \mbox{if $\sigma_{i}+\sigma_{j}\neq 0$}\,,\\[3pt]
0 & \mbox{otherwise}\,,
\end{array} \right.  \quad i,j\in\{1,\ldots,m\}\,,
\label{eq:def-matric-E2}
\\[3pt]
({\cal F}^0({\sigma}))_{ij} &:=&
\left\{\begin{array}{ll}
\displaystyle{\frac{g_i({\sigma})}{\sigma_{i}}} & \mbox{if $\sigma_{i}\neq 0$}\,,\\
0 & \mbox{otherwise},
\end{array}\right. \quad  i\in\{1,\ldots,m\},\quad j\in\{1,\ldots, n-m\}\,.
\label{eq:def-matric-F}
\end{eqnarray}
Note that when the dependence of ${\cal E}^0_{1}({\sigma})$, ${\cal E}^0_{2}({\sigma})$ and ${\cal F}^0({\sigma})$  on $\sigma$ are clear from the context, we often drop ${\sigma}$ from these notations. In particular, let $\overline{\cal E}^0_{1}$, $\overline{\cal E}^0_{2}\in\V^{m\times m}$ and $\overline{\cal F}^0\in\V^{m\times(n-m)}$ be the matrices defined by (\ref{eq:def-matric-E1})-(\ref{eq:def-matric-F}) with respect to $\overline{\sigma}=\sigma(\overline{X})$. Since $g$ is absolutely symmetric at $\overline{\sigma}$, we know that for all $i\in a_l$, $1\leq l\leq r$, the function values $g_i(\overline{\sigma})$ are the same (denoted by $\bar{g}_l$). Therefore, for any $X\in{\cal N}$, define
\begin{equation}\label{eq:def-GS}
G_{S}({X}):=\sum_{l=1}^{r}\bar{g}_{l}{\cal U}_{l}(X) \quad {\rm and}\quad G_{R}(X):=G(X)-G_{S}(X)\,,
\end{equation}
where ${\cal U}_{l}(X)$ is given by (\ref{eq:Pk-nonsymmetric}).
The following lemma follows from  Proposition \ref{prop:P''_k} directly.
\begin{lemma}\label{lem:Gs-diff-formula}
Let ${G}_{S} :{\cal N}\to\V^{m\times n}$ be defined by (\ref{eq:def-GS}). Then, there exists an open neighborhood ${\cal B}$ of $\overline{X}$ in ${\cal N}$ such that $G_{S} $ is twice continuously differentiable on ${\cal B}$, and for any $\V^{m\times n}\ni H\to 0$,
\[
G_{S}(\overline{X}+H)-G_{S}(\overline{X})=G'_{S}(\overline{X})H+O(\|H\|^{2})
\] with
\begin{equation}\label{eq:Gs'-formula}
G'_{S}(\overline{X})H=\overline{U}\left[\overline{\cal E}^0_1\circ S(\overline{U}^\T H\overline{V}_{1})+\overline{\cal E}^0_2\circ T(\overline{U}^\T H\overline{V}_{1})
 \quad \overline{\cal F}^0\circ (\overline{U}^\T H\overline{V}_{2})  \right]\overline{V}^\T \,.
\end{equation}
\end{lemma}

Lemma \ref{lem:Gs-diff-formula} says that in an open  neighborhood of $\overline{X}$,   $G(\cdot)$ can be decomposed into a ``smooth part" $G_S(\cdot)$ plus a ``nonsmooth part" $G_R(\cdot)$. As we will see in the later developments, this decomposition   simplifies many of our proofs.

\subsection{Directional differentiability}\label{subsection:Directional differentiability}

Let ${\cal Z}$ and ${\cal Z}^\prime$ be two finite dimensional real Euclidean spaces and  ${\cal O}$ be an open set in ${\cal Z}$. A function $F: {\cal O}  \subseteq {\cal Z}  \to  {\cal Z}^\prime$ is said to be  {\it Hadamard directionally differentiable} at $z\in{\cal O}$ if the limit
\begin{equation}\label{eq:def-H-dir-diff}
\lim_{t\downarrow 0,\; h'\to h}\; \frac{F(z+th')-F(z)}{t} \quad \mbox{exists for any $h\in{\cal Z}$}\,.
\end{equation}
It is clear that if $F$ is Hadamard directionally differentiable at $z$, then $F$ is directionally differentiable at $z$, and the limit in (\ref{eq:def-H-dir-diff}) equals the directional derivative $F'(z;h)$ for any $h\in{\cal Z}$.

Assume that the $g$ is directionally differentiable at $\overline{\sigma}$. Then, from the definition of directional derivative and the absolutely symmetry of $g$ on the nonempty open set $\hat{\sigma}_{ {\cal N}}$, it is easy to see that the directional derivative $\phi :=g'(\overline{\sigma};\cdot):\Re^{m}\to\Re^{m}$ satisfies
\begin{equation}\label{eq:dir-diff-symmetric}
\phi(Qh)=Q\phi(h)\quad \forall\, Q\in\pm\mathbb{P}_{\overline{\sigma}}^m\quad {\rm and} \quad \forall\,h\in\Re^{m}\,,
\end{equation} where $\pm\mathbb{P}_{\overline{\sigma}}^m$ is the subset defined with respect to $\overline{\sigma}$ by
\begin{equation}\label{eq:def-Q-permu}
\pm\mathbb{P}_{\overline{\sigma}}^m:=\left\{{Q}\in\pm\mathbb{P}^m\,|\, \overline{\sigma}={Q}\overline{\sigma}\right\}\,.
\end{equation}
Note that $Q\in\pm\mathbb{P}_{\overline{\sigma}}^m$ if and only if
\begin{equation}\label{eq:block-structure-P_sigma}
Q= {\rm Diag}\big(Q_1,\dots,Q_{r},Q_{r+1} \big)
\quad {\rm with} \quad Q_{l}\in\mathbb{P}^{|a_{l}|}, \quad l=1,\ldots,r \quad {\rm and} \quad Q_{r+1}\in\pm\mathbb{P}^{|b|}\,.
\end{equation}
For any $h\in\Re^{m}$, we rewrite $\phi(h)$  in the following form $\phi(h)=\left(\phi_{1}(h), \ldots,\phi_{r}(h),\phi_{r+1}(h)  \right)$
with $\phi_{l}(h)\in\Re^{|a_{l}|}$, $l=1,\ldots,r$ and $\phi_{r+1}(h)\in\Re^{|b|}$. Therefore, we know from \eqref{eq:dir-diff-symmetric} and \eqref{eq:block-structure-P_sigma} that the function $\phi:\Re^{m}\to\Re^{m}$ is a mixed symmetric mapping, with respect to
$ \mathbb{P}^{|a_1|}\times\ldots\times\mathbb{P}^{|a_r|}\times\pm\mathbb{P}^{|b|}$,
 over $\Re^{|a_{1}|}\times \ldots \times  \Re^{|a_{r}|}\times\Re^{|b|}$.
Let ${\cal W}:=\S^{|a_{1}|}\times\ldots\times\S^{|a_{r}|}\times\V^{|b|\times(n-|a|)}$.  We can define the spectral operator $\Phi :{\cal W}\to{\cal W}$ with respect to the symmetric mapping $\phi$ as follows: for any $W=\left(W_{1}, \ldots,  W_{r}, W_{r+1}\right)\in{\cal W}$,
\begin{equation}\label{eq:dir-derivative-spectral-op}
\Phi(W):=\Big(\Phi_{1}(W),\ldots, \Phi_{r}(W), \Phi_{r+1}(W) \Big)
\end{equation} with
\[
\Phi_{l}(W):=\left\{ \begin{array}{ll}
\widetilde{P}_{l}{\rm Diag}(\phi_{l}({\kappa}({W})))\widetilde{P}_{l}^\T  & \mbox{if $1\leq l\leq r$,}\\[3pt]
\widetilde{M} {\rm Diag}(\phi_{l}({\kappa}({W}))) \widetilde{N}_{1}^\T  & \mbox{if $l=r+1$,}
\end{array}\right.
\]
 where ${\kappa}({W}):=\left(\lambda({W}_{1}),\ldots,\lambda({W}_{r}),\sigma({W}_{r+1})\right)\in\Re^{m}$; $\widetilde{P}_{l}\in{\mathbb O}^{|a_{l}|}(W_{l})$; and $(\widetilde{M},\widetilde{N})\in{\mathbb O}^{|b|,n-|a|}(W_{r+1})$, $\widetilde{N}:=\big[\widetilde{N}_{1}\quad \widetilde{N}_{2}\big]$ with $\widetilde{N}_{1}\in\V^{(n-|a|)\times|b|}$, $\widetilde{N}_{2}\in\V^{(n-|a|)\times(n-m)}$. From Theorem \ref{thm:well-def-spectral-op}, we know that $\Phi$ is well defined on ${\cal W}$.

In order to present   the directional differentiability results for the spectral operator $G$, we define the following first divided directional difference $g^{[1]}(\overline{X};H)\in\V^{m\times n}$ of $g$ at $\overline{X}$ along the direction $H\in\V^{m\times n}$ by
\begin{eqnarray}
g^{[1]}(\overline{X};H)&:=&\left[\overline{\cal E}^0_1\circ S(\overline{U}^\T H\overline{V}_{1})+\overline{\cal E}^0_2\circ T(\overline{U}^\T H\overline{V}_{1})
\quad \overline{\cal F}^0\circ \overline{U}^\T H\overline{V}_{2}  \right] + {\widehat \Phi}(D(H)),
\label{eq:divid-dir-diff-matrix-sym-op-2}
\end{eqnarray}
where $\overline{\cal E}_1, \overline{\cal E}_2, \overline{\cal F}$ are
defined as in (\ref{eq:def-matric-E1})--(\ref{eq:def-matric-F}) at $\overline \sigma=\sigma(\overline{X})$,
\begin{equation}\label{eq:def-D-mapping}
D(H):=\left(  S(\overline{U}_{a_1}^\T H\overline{V}_{a_1}),\ldots,S(\overline{U}_{a_r}^\T H\overline{V}_{a_r}),\overline{U}_{b}^\T H[\overline{V}_{b}\quad \overline{V}_{2}] \right)\in{\cal W}
\end{equation}
and for any $W=\left(W_{1}, \ldots,  W_{r}, W_{r+1}\right)\in{\cal W}$, ${\widehat \Phi}(W)\in \V^{m\times n} $ is defined by
\begin{equation}\label{eq:def-Diag-Phi}
{\widehat \Phi}(W):= \left[\begin{array}{cc} {\rm Diag}\left(\Phi_{1}(W), \dots, \Phi_{r}(W)\right)  & 0
\\[2mm]
0 &  \Phi_{r+1}(W)
\end{array}\right].
\end{equation}

For the directional differentiability of the spectral operator, we have the following result.

\begin{theorem}\label{prop:H-dir-diff-spectr-op}
Suppose that $\overline{X}\in{\cal N}$ has the SVD (\ref{eq:Y-eig-Z-SVD}). The spectral operator $G$ is Hadamard directionally differentiable at $\overline{X}$ if and only if $g$ is Hadamard directionally differentiable at $\overline{\sigma}=\sigma(\overline{X})$. In that case,   the directional derivative of $G$ at $\overline{X}$ along any direction $H\in\V^{m\times n}$ is given by
\begin{equation}\label{eq:dir-diff-spectr-op}
G'(\overline{X};H)=\overline{U}g^{[1]}(\overline{X};H) \overline{V}^\T \,.
\end{equation}
\end{theorem}
\noindent {\bf  Proof.} $``\Longleftarrow"$ Let $H\in\V^{m\times n}$ be any given direction. For any $\V^{m\times n}\ni H'\to H$ and $\tau>0$, denote $X:=\overline{X}+\tau H'$. Consider the SVD of $X$, i.e.,
\begin{equation}\label{eq:def-Y-Z-dir-diff}
X=U[\Sigma(X)\quad 0]V^\T \,.
\end{equation}
Denote $\sigma=\sigma(X)$. For $\tau$ and $H'$ sufficiently close to $0$ and $H$, let $G_{S}$ and $G_{R}$ be the mappings defined in (\ref{eq:def-GS}).
Then, by Lemma \ref{lem:Gs-diff-formula}, we know that
\begin{equation}\label{eq:smoothpart-dir-diff-Spectral-op}
\lim_{\tau\downarrow 0,\; H'\to H}\frac{1}{\tau}(G_{S}(X)-G_{S}(\overline{X}))=G'_{S}(\overline{X})H\,,
\end{equation}
where $G'_{S}(\overline{X})H$ is given by \eqref{eq:Gs'-formula}.
On the other hand, for $\tau$ and $H'$ sufficiently close to $0$ and $H$,
we have ${\cal U}_{l}(X)={\sum_{i\in a_{l}}}u_{i}v_{i}^\T $, $l=1,\ldots,r$ and
\begin{equation}\label{eq:GR-def-dir-diff}
G_{R}(X)=G(X)-G_{S}(X)= \sum_{l=1}^{r}\sum_{i\in a_{l}}[g_i(\sigma)-g_i(\overline{\sigma})]u_{i}v_{i}^\T +\sum_{i\in b}g_i(\sigma)u_{i}v_{i}^\T  \,.
\end{equation}
For $\tau$ and $H'$ sufficiently close to $0$ and $H$,  let
\[
\Delta_{l}(\tau,H')=\frac{1}{\tau}\sum_{i\in a_{l}}[g_i(\sigma)-g_i(\overline{\sigma})]u_{i}v_{i}^\T ,\quad l=1,\ldots,r \quad {\rm and} \quad \Delta_{r+1}(\tau,H')=\sum_{i\in b}g_i(\sigma)u_{i}v_{i}^\T \,.
\]

Firstly, consider the case that $\overline{X}=[\Sigma(\overline{X})\quad 0 ]$. Then, from (\ref{eq:directi-diff}) and (\ref{eq:directi-diff-sigma}), we know that for any $\tau$  and $H'\in\V^{m\times n}$ sufficiently close to $0$ and $H$,
\begin{equation}\label{eq:limit-eig-singular-dir-diff-1} \sigma(X)=\sigma(\overline{X})+\tau\sigma'(\overline{X};H')+O(\tau^{2}\|H'\|^{2})\,,
\end{equation} where $(\sigma'(\overline{X};H'))_{a_{l}}=\lambda(S(H'_{a_{l}a_{l}}))$, $l=1,\ldots,r$ and $(\sigma'(\overline{X};H'))_{b}=\sigma([H'_{bb}\quad H'_{bc}])$. Denote $h':= \sigma'(\overline{X};H')$  and $h:=\sigma'(\overline{X};H)$. By using the fact that the singular value functions of a general matrix are globally Lipschitz continuous, we know that
\begin{equation}\label{eq:limit-eig-singular-dir-diff-2}
\lim_{\tau\downarrow 0, \; H'\to H}\; (h'+O(\tau\|H'\|^{2}))=h\,.
\end{equation}   Since $g$ is assumed to be Hadamard directionally differentiable at $\overline{\sigma}$, we have
\[
\lim_{\tau\downarrow 0,\; H'\to H}\; \frac{g(\sigma)-g(\overline{\sigma})}{\tau}= \lim_{\tau\downarrow 0,\; H'\to H}\; \frac{1}{\tau} [g(\overline{\sigma}+\tau(h'+O(\tau\|H'\|^{2})))-g(\overline{\sigma})]= g'(\overline{{\sigma}};h)=\phi(h)\,,
\]
where $\phi\equiv g'(\overline{\sigma};\cdot):\Re^{m}\to\Re^{m}$ satisfies the
condition (\ref{eq:dir-diff-symmetric}). Since $u_{i}v_{i}^\T $, $i=1,\ldots,m$ are uniformly bounded, we know that for $\tau$ and $H'$ sufficiently close to $0$ and $H$,
\begin{eqnarray*}
\Delta_{l}(\tau,H') &=& U_{a_{l}}{\rm Diag}(\phi_{l}(h))V_{a_{l}}^\T +o(1) \quad l=1,\ldots,r\,,
\\
\Delta_{r+1}(\tau,H') &=&U_{b}{\rm Diag}(\phi_{r+1}(h))V_{b}^\T +o(1)\,.
\end{eqnarray*}
By (\ref{eq:diag-U-V}) in Proposition \ref{prop:diag-U-V}, we know that there exist $Q_{l}\in{\mathbb O}^{|a_{l}|}$, $l=1,\ldots,r$, $M\in{\mathbb O}^{|b|}$ and $N=[N_{1}\quad N_{2}]\in{\mathbb O}^{n-|a|}$ with $N_{1}\in\V^{(n-|a|)\times |b|}$ and $N_{2}\in\V^{(n-|a|)\times(n-m)}$ (depending on $\tau$ and $H'$) such that
\begin{eqnarray*}
U_{a_{l}}&=&
\left[\begin{array}{c} O(\tau\|H'\|)  \\ Q_{l}+O(\tau\|H'\|) \\ O(\tau\|H'\|) \end{array}\right], \quad  V_{a_{l}}=\left[\begin{array}{c} O(\tau\|H'\|)  \\ Q_{l}+O(\tau\|H'\|) \\ O(\tau\|H'\|) \end{array}\right]  \ l=1,\ldots,r\,,
\\
U_{b}&=&\left[\begin{array}{c} O(\tau\|H'\|) \\ M+O(\tau\|H'\|)\end{array}\right],\quad  [V_{b}\quad V_{c}]=\left[\begin{array}{c} O(\tau\|H'\|) \\
N+O(\tau\|H'\|)\end{array}\right]  \,.
\end{eqnarray*}
Thus,  we have
\begin{eqnarray}
\Delta_{l}(\tau,H') &=&\left[\begin{array}{ccc} 0 & 0  & 0  \\[3pt] 0 & Q_{l}{\rm Diag}(\phi_{l}(h))Q_{l}^\T   & 0   \\[3pt] 0 & 0  & 0 \end{array}\right]
+O(\tau\|H'\|)+o(1),\quad l=1,\ldots,r\,,
\label{eq:Delta-k-dir-diff-1}
\\
\Delta_{r+1}(\tau,H')&=&
\left[\begin{array}{cc}0 & 0  \\0 & M{\rm Diag}(\phi_{r+1}(h))N_{1}^\T  \end{array}\right]+O(\tau\|H'\|)+o(1)\,.
\label{eq:Delta-k-dir-diff-2}
\end{eqnarray}
We know from \eqref{eq:diag-H-decomp-ak} and \eqref{eq:diag-H-decomp-b} that
\begin{eqnarray}
S(H'_{a_{l}a_{l}}) &=&
S(H_{a_{l}a_{l}})+o(1)=\frac{1}{\tau}Q_{l}[\Sigma(X)_{a_{l}a_{l}}-\overline{\nu}_{l}I_{|a_{l}|}]Q_{l}^\T +O(\tau\|H'\|^{2}),\quad l=1,\ldots,r\,,\quad
\label{eq:H-H'-2-dir-diff}\\ [1pt]
[H'_{bb}\quad H'_{bc}] &=&
[H_{bb}\quad H_{bc}]+o(1)=\frac{1}{\tau}M[\Sigma(X)_{bb}-\overline{\nu}_{r+1}I_{|b|}]N_{1}^\T +O(\tau\|H'\|^{2})\,.
\label{eq:H-H'-3-dir-diff}
\end{eqnarray}
Since $Q_{l}$, $l=1,\ldots,r$, $M$ and $N$ are uniformly bounded, by taking  subsequences if necessary, we may assume that when $\tau\downarrow 0$ and $H'\to H$, $Q_{l}$, $M$ and $N$ converge to  $\widetilde{Q}_{l}$, $\widetilde{M}$ and $\widetilde{N}$, respectively.  Therefore, by taking limits in (\ref{eq:H-H'-2-dir-diff}) and (\ref{eq:H-H'-3-dir-diff}), we obtain from (\ref{eq:limit-eig-singular-dir-diff-1}) and (\ref{eq:limit-eig-singular-dir-diff-2}) that
\begin{eqnarray*}
 S(H_{a_{l}a_{l}}) &=& \widetilde{Q}_{l}\Lambda(S(H_{a_{l}a_{l}})) \widetilde{Q}_{l}^\T , \quad l=1,\ldots,r\,,
\\[3pt]
\left[ H_{bb}\quad H_{bc} \right] &=&
\widetilde{M} \left[ \Sigma(\left[ H_{bb}\quad H_{bc} \right])\quad 0\right] \widetilde{N}^\T = \widetilde{M} \Sigma(\left[ H_{bb}\quad H_{bc} \right]) \widetilde{N}_{1}^\T .
\end{eqnarray*}
Hence, by using the notation (\ref{eq:dir-derivative-spectral-op}), we know from (\ref{eq:GR-def-dir-diff}), (\ref{eq:Delta-k-dir-diff-1}), (\ref{eq:Delta-k-dir-diff-2}) and (\ref{eq:def-Diag-Phi})   that
\begin{eqnarray}
\lim_{\tau\downarrow 0, \; H'\to H}\frac{1}{\tau}G_{R}(X)
&=&\lim_{\tau\downarrow 0, \; H'\to H}\sum_{l=1}^{r+1}\Delta_{l}(\tau,H')
\;=\;
{\widehat \Phi}(D(H))\,,
\label{eq:nonsmoothpart-dir-diff-Spectral-op-diag}
\end{eqnarray}
where $D(H)=\left(S(H_{a_{1}a_{1}}),\ldots,S(H_{a_{r}a_{r}}), H_{b\bar{a}}\right)$.

To prove the conclusion for  the general case of $\overline X$, rewrite (\ref{eq:def-Y-Z-dir-diff}) as
\[
\left[\Sigma(\overline{X})\quad 0 \right]+\overline{U}^\T H'\overline{V}=\overline{U}^\T U[\Sigma(X)\quad 0]V^\T \overline{V}\,.
\]
Let  $\widetilde{U}:=\overline{U}^\T U$, $\widetilde{V}:=\overline{V}^\T V$
and $\widetilde{H} = \overline{U}^\T H \overline{V}$. Denote $\widetilde{X}:=[\Sigma(\overline{X})\quad 0 ]+\overline{U}^\T H'\overline{V}$. Then, we obtain that $G_{R}(X)=\overline{U} G_{R} (\widetilde{{X}}) \overline{V}^\T $.
Thus, we know from (\ref{eq:nonsmoothpart-dir-diff-Spectral-op-diag}) that
\begin{equation}\label{eq:nonsmoothpart-dir-diff-Spectral-op}
\lim_{\tau\downarrow 0, \; H'\to H}\frac{1}{\tau}G_{R}(X)
= \overline{U}  {\widehat \Phi}(D(\widetilde{H}))    \overline{V}^\T \,.
\end{equation}
Therefore, by combining (\ref{eq:smoothpart-dir-diff-Spectral-op}) and (\ref{eq:nonsmoothpart-dir-diff-Spectral-op}) and noting that $G(\overline{X})=G_{S}(\overline{X})$, we obtain that for any given $ H\in\V^{m\times n}$,
\[
\lim_{\tau\downarrow 0,\; H'\to H}\frac{G(X)-G(\overline{X})}{\tau}=\lim_{\tau\downarrow 0,\; H'\to H}\frac{G_{S}(X)-G_{S}(\overline{X})+G_{R}(X)}{\tau}=\overline{U}g^{[1]}(\overline{X};\widetilde{H}) \overline{V}^\T \,,
\] where $g^{[1]}(\overline{X};\widetilde{H})$ is given by (\ref{eq:divid-dir-diff-matrix-sym-op-2}). This implies that $G$ is Hadamard directionally differentiable at $\overline{X}$ and (\ref{eq:dir-diff-spectr-op}) holds.

$``\Longrightarrow"$ Suppose that $G$ is Hadamard directionally differentiable at $\overline{X}$. Let $(\overline{U},\overline{V})\in{\mathbb O}^{m\times n}(\overline{X})$ be fixed.  For any given direction $h\in\Re^{m}$, suppose that $\Re^{m}\ni h'\to h$. Denote $H':=\overline{U}[{\rm Diag}(h')\quad 0]\overline{V}^\T \in\V^{m\times n}$ and   $H:=\overline{U}[{\rm Diag}(h)\quad 0]\overline{V}^\T \in\V^{m\times n}$. Then, we have $H'\to H$ as $h'\to h$.  Since for all $\tau>0$ and $h'$ sufficiently close to $0$ and $h$, $\sigma:=\overline{\sigma}+\tau h'\in\hat{\sigma}_{\cal N}$, we know from Proposition \ref{prop:G-g} that for all $\tau>0$ and $h'$ sufficiently close to $0$ and $h$, $G(\overline{X}+\tau H')=\overline{U}{\rm Diag}(g(\overline{\sigma}+\tau h'))\overline{V}_1^\T$. This  implies that
 \[
{\rm Diag}\Big(\lim_{\tau\downarrow 0,\; h'\to h}\frac{g(\overline{\sigma}+\tau h')-g(\overline{\sigma})}{\tau}\Big)= \overline{U}^\T\left(\lim_{\tau\downarrow 0,\; H'\to H}\frac{G(\overline{X}+\tau H')-G(\overline{X})}{\tau}\right)\overline{V}_1\,.
 \]
 Thus, we know from the assumption that $\displaystyle{\lim_{\tau\downarrow 0,\; h'\to h}}\frac{g(\overline{\sigma}+\tau h')-g(\overline{\sigma})}{\tau}$ exists and  that $g$ is Hadamard directionally differentiable at $\overline{\sigma}$.
$\hfill \Box$
\vskip 10 true pt

\begin{remark}\label{remark:Lip-H-diff-spectral-op}
Note that for a general spectral operator $G$, we cannot obtain the directional differentiability at $\overline{X}$ if we only assume that $g$ is directionally differentiable at $\sigma({\overline{X}})$. In fact, a counterexample can be found in \cite{Lewis96}. However, since $\V^{m\times n}$ is a finite dimensional Euclidean space, it is well-known that for locally Lipschitz continuous functions,  the directional differentiability in the sense of Hadamard and G\^{a}teaux are equivalent (see e.g. \cite[Theorem 1.13]{Nashed71}, \cite[Lemma 3.2]{DRubinov83}, \cite[p.259]{Flett80}). Therefore, if   $G$ and $g$ are locally Lipschitz continuous near $\overline{X}$ and $\sigma({\overline{X}})$, respectively (e.g., the proximal mapping $P_{f}$ and its vector counterpart $P_{\theta}$), then $G$ is directionally differentiable at $\overline{X}$ if and only if $g$ is directionally differentiable at $\sigma(\overline{X})$.
\end{remark}

\subsection{Fr\'{e}chet differentiability}

For a given $X\in{\cal N}$, suppose that the given absolutely symmetric mapping $g$    is F(r\'{e}chet)-differentiable at $\sigma=\sigma(X)$. The following results on the Jacobian matrix $g'(\sigma)$ can be obtained directly from the assumed  absolute symmetry of $g$ on $\hat{\sigma}_{ {\cal N}}$ and the block structure \eqref{eq:block-structure-P_sigma} for any
$Q\in\pm\mathbb{P}_{{\sigma}}^m$.
\begin{lemma}\label{lem:structure-g'-gen}
For any  $X\in{\cal N}$, suppose that $g$ is F-differentiable at $\sigma=\sigma(X)$. Then, the Jacobian matrix $g'(\sigma)$ has the following property
\begin{equation*}\label{eq:struc-symm-g'-gen}
g'(\sigma)=Q^\T g'(\sigma)Q\quad \forall\, Q\in\pm\mathbb{P}_{{\sigma}}^m\,.
\end{equation*}
In particular,
\[
\left\{
\begin{array}{ll}
(g'(\sigma))_{ii}=(g'(\sigma))_{i'i'} & \mbox{if $\sigma_{i}=\sigma_{i'}$ and $i,i'\in\{1,\ldots,m\}$,}\\[3pt]
(g'(\sigma))_{ij}=(g'(\sigma))_{i'j'} & \mbox{if $\sigma_{i}=\sigma_{i'}$, $\sigma_{j}=\sigma_{j'}$, $i\neq j$, $i'\neq j'$ and $i,i',j,j'\in\{1,\ldots,m\}$,}\\[3pt]
(g'(\sigma))_{ij}=(g'(\sigma))_{ji}=0 & \mbox{if $\sigma_{i}=0$, $i\neq j$ and $i,j\in\{1,\ldots,m\}$.}
\end{array}
\right.
\]
\end{lemma}

Lemma \ref{lem:structure-g'-gen} is a simple extension of \cite[Lemma 2.1]{LSendov03} for  symmetric mappings. But one should  note that the  Jacobian matrix $g'(\sigma)$ of $g$ at the F-differentiable point $\sigma$ may not be symmetric since here $g$ is not assumed to be  the gradient mapping as in  \cite[Lemma 2.1]{LSendov03}. For example, let us consider the  absolutely symmetric mapping $g$ defined by \eqref{eq:example-Miao} in the introduction. Then $g$ is differentiable at $x=(2,1)$ by taking   $m=2$ and $\tau=\varepsilon=1$ in \eqref{eq:h-example-Miao}. However, it is easy to see that the Jacobian matrix $g'(x)$ is not symmetric.

  Let  $\eta(\sigma)\in\Re^m$  be the vector defined by
\begin{equation}\label{eq:def-eta}
(\eta(\sigma))_{i}:=\left\{ \begin{array}{ll} (g'({\sigma}))_{ii}-(g'({\sigma}))_{i(i+1)} & \mbox{if $\exists\, j\in\left\{1,\ldots,m\right\}$ and $j\ne i$  such that   ${\sigma}_i={\sigma}_j$},\\
(g'({\sigma}))_{ii} & \mbox{otherwise}\,,
\end{array} \right. \quad i\in\{1,\ldots,m\}\,.
\end{equation}
Define the corresponding {\it divided difference matrix}  ${\cal E}_{1}(\sigma)\in\Re^{m\times m}$, the {\it divided addition matrix} ${\cal E}_{2}(\sigma)\in\Re^{m\times m}$,  the {\it division matrix} ${\cal F}(\sigma)\in\Re^{m\times(n-m)}$, respectively, by
\begin{equation}\label{eq:def-matric-ED1}
({\cal E}_{1}(\sigma))_{ij}:=\left\{ \begin{array}{ll}  \displaystyle{\frac{g_{i}(\sigma)-g_{j}(\sigma)}{\sigma_{i}-\sigma_{j}}} & \mbox{if $\sigma_{i}\neq\sigma_{j}$}\,,\\
(\eta(\sigma))_{i} & \mbox{otherwise}\,,
\end{array} \right.  \quad i,j\in\{1,\ldots,m\}\,,
\end{equation}
\begin{equation}\label{eq:def-matric-ED2}
({\cal E}_{2}(\sigma))_{ij}:=\left\{ \begin{array}{ll}  \displaystyle{\frac{g_{i}(\sigma)+g_{j}(\sigma)}{\sigma_{i}+\sigma_{j}}} & \mbox{if $\sigma_{i}+\sigma_{j}\neq 0$}\,,\\
g'(\sigma))_{ii} & \mbox{otherwise}\,,
\end{array} \right.  \quad i,j\in\{1,\ldots,m\}\,,
\end{equation}
\begin{equation}\label{eq:def-matric-FD}
({\cal F}(\sigma))_{ij}:=\left\{\begin{array}{ll}
\displaystyle{\frac{g_{i}(\sigma)}{\sigma_{i}}} & \mbox{if $\sigma_{i}\neq 0$}\,,\\
(g'(\sigma))_{ii} & \mbox{otherwise,}
\end{array}\right. \quad  i\in\{1,\ldots,m\},\quad j\in\{1,\ldots, n-m\}\, .
\end{equation}
Define  the matrix ${\cal C}(\sigma)\in  \Re^{m\times m}$ to be the difference between $g'(\sigma)$ and ${\rm Diag}(\eta(\sigma))$, i.e.,
\begin{equation}\label{eq:def-matric-C}
{\cal C}(\sigma):=g'(\sigma)-{\rm Diag}(\eta(\sigma)) \,.
\end{equation}
Note that when the dependence of $\eta$, ${\cal E}_{1}$, ${\cal E}_{2}$, ${\cal F}$ and $\cal C$  on $\sigma$ is clear from the context, we often drop $\sigma$ from the corresponding notations.

Let $\overline{X}\in{\cal N}$ be given and denote $\overline \sigma = \sigma(\overline{X})$.
Denote $\overline{\eta}=\eta(\overline \sigma)\in\Re^m$ to be the vector defined by \eqref{eq:def-eta}. Let $\overline{\cal E}_{1}$, $\overline{\cal E}_{2}$, $\overline{\cal F}$ and $\overline{\cal C}$ be the real matrices defined in (\ref{eq:def-matric-ED1})--(\ref{eq:def-matric-C}) with respect to $\overline{\sigma}$.
Now, we are ready to state the result on the F-differentiability of spectral operators.

\begin{theorem}\label{thm:diff-spectral-op}
Suppose that the given matrix $\overline{X}\in{\cal N}$ has the SVD (\ref{eq:Y-eig-Z-SVD}). Then the spectral operator $G$ is  F-differentiable at $\overline{X}$ if and only if $g$ is F-differentiable at $\overline{\sigma}$. In that case, the derivative of $G$ at $\overline{X}$ is given by
\begin{equation}\label{eq:F-derivative-spectral-op}
G'(\overline{X})H= \overline{U}[\overline{\cal E}_{1}\circ S(A)+{\rm Diag}\left(\overline{\cal C}{\rm diag}(S( {A}))\right) +\overline{\cal E}_{2}\circ T( {A})\quad \overline{\cal F}\circ  {B} ]\overline{V}^\T  \quad
\forall\; H\in\V^{m\times n},
\end{equation}
where  $ {A}:=\overline{U}^\T H\overline{V}_1$ and  $ {B}:=\overline{U}^\T H\overline{V}_2$.
\end{theorem}
\noindent {\bf  Proof.} $``\Longleftarrow"$ For any $\V^{m\times n}\ni H=[H_1\quad H_2]\to{0}$ with $H_1\in\V^{m\times m}$ and $H_2\in\V^{m\times(n-m)}$, denote $X=\overline{X}+H$. Let $U\in{\mathbb O}^{m}$ and $V\in{\mathbb O}^{n}$ be such that
\begin{equation}\label{eq:def-Y-Z-F-diff}
X=U[\Sigma(X)\quad 0]V^\T \,.
\end{equation}
Denote $\sigma=\sigma(X)$. Let $G_{S}(X)$ and $G_{R}(X)$ be defined by (\ref{eq:def-GS}). Then, by Lemma \ref{lem:Gs-diff-formula}, we know that for any $H\to{0}$,
\begin{equation}\label{eq:smoothpart-F-diff-Spectral-op}
G_{S}(X)-G_{S}(\overline{X})=G'_{S}(\overline{X})H+O(\|H\|^{2})=G_S'(\overline{X})H+O(\|H\|^{2})\,,
\end{equation}
where $G_S'(\overline{X})H$ is given by \eqref{eq:Gs'-formula}. For $H\in\V^{m\times n}$ sufficiently small, we have ${\cal U}_{l}(X)={\sum_{i\in a_{l}}}u_{i}v_{i}^\T $, $l=1,\ldots,r$. Therefore,
\begin{equation}\label{eq:GR-def-F-diff-spectral-op}
G_{R}(X)=G(X)-G_{S}(X)=\sum_{l=1}^{r+1}\Delta_{l}(H)\,,
\end{equation}
where $\Delta_{l}(H)=\sum_{i\in a_{l}}(g_{i}(\sigma)-g_i(\overline{\sigma}))u_{i}v_{i}^\T$, $l=1,\ldots,r$ and $\Delta_{r+1}(H)=\sum_{i\in b}g_i(\sigma)u_{i}v_{i}^\T$.

Firstly, consider the case that $\overline{X}=[\Sigma(\overline{X})\quad 0 ]$. Then, from (\ref{eq:directi-diff}) and (\ref{eq:directi-diff-sigma}), for any  $H\in\V^{m\times n}$ sufficiently small, we have
\begin{equation}\label{eq:eig-singular-F-diff-1}
\sigma= \overline{\sigma}+h+O(\|H\|^{2})\,,
\end{equation} where $h:=\sigma'(\overline{X};H)\in\Re^{m}$ with
\begin{equation}\label{eq:eig-singular-F-diff-1-a-added}
(\sigma'(\overline{X};H))_{a_{l}}=\lambda(S(H_{a_{l}a_{l}})),\  l=1,\ldots,r \quad {\rm and} \quad (\sigma'(\overline{X};H))_{b}=\sigma([H_{bb}\quad H_{bc}])\,.
\end{equation}
Since $g$ is F-differentiable at $\overline{\sigma}$ and the singular value functions are globally Lipschitz continuous, we know from \eqref{eq:eig-singular-F-diff-1} that for any $H\in\V^{m\times n}$ sufficiently small,
\[
g(\sigma)-g(\overline{\sigma})=g(\overline{\sigma}+h+O(\|H\|^{2}))-g(\overline{\sigma})
=g'(\overline{\sigma})(h+O(\|H\|^{2}))+o(\|h\|)=g'(\overline{\sigma})h+o(\|H\|)\,.
\]
Since $u_{i}v_{i}^\T $, $i=1,\ldots,m$ are uniformly bounded, we have for $H$ sufficiently
small,
\[
\Delta_{l}(H)=
U_{a_{l}}{\rm Diag}((g'(\overline{\sigma})h)_{a_{l}})V_{a_{l}}^\T +o(\|H\|), \;
 l=1,\ldots,r,
 \quad
\Delta_{r+1}(H)=U_{b}{\rm Diag}((g'(\overline{\sigma})h)_{b})V_{b}^\T +o(\|H\|)\,.
\]
By (\ref{eq:diag-U-V}) in Proposition \ref{prop:diag-U-V}, we know that there exist $Q_{l}\in{\mathbb O}^{|a_{l}|}$,  $M\in{\mathbb O}^{|b|}$ and $N=[N_{1}\quad N_{2}]\in{\mathbb O}^{n-|a|}$ with $N_{1}\in\V^{(n-|a|)\times |b|}$ and $N_{2}\in\V^{(n-|a|)\times(n-m)}$ (depending on $H$) such that
\begin{eqnarray*}
U_{a_{l}}&=&\left[\begin{array}{c} O(\|H\|)  \\ Q_{l}+O(\|H\|) \\ O(\|H\|) \end{array}\right], \quad  V_{a_{l}}=\left[\begin{array}{c} O(\|H\|)  \\ Q_{l}+O(\|H\|) \\ O(\|H\|) \end{array}\right], \ l=1,\ldots,r\,,
\\
U_{b}&=&\left[\begin{array}{c} O(\|H\|) \\[3pt] M+O(\|H\|)\end{array}\right],\quad  [V_{b}\quad V_{c}]=\left[\begin{array}{c} O(\|H\|) \\[3pt]
N+O(\|H\|)\end{array}\right]  \,.
\end{eqnarray*}
Therefore, since $\|g'(\overline{\sigma})h\|=O(\|H\|)$,  we obtain that
\begin{eqnarray}
\Delta_{l}(H)&=& \left[\begin{array}{ccc} 0 & 0  & 0  \\[3pt] 0 & Q_{l}{\rm Diag}((g'(\overline{\sigma})h)_{a_{l}})Q_{l}^\T   & 0   \\[3pt] 0 & 0  & 0 \end{array}\right]
+o(\|H\|),\quad l=1,\ldots,r\,,
\label{eq:Delta-k-F-diff-1}
\\[3pt]
\Delta_{r+1}(H)&=& \left[\begin{array}{cc}0 & 0  \\0 & M{\rm Diag}((g'(\overline{\sigma})h)_{b})N_{1}^\T  \end{array}\right]+o(\|H\|)\,.
\label{eq:Delta-k-F-diff-2}
\end{eqnarray}
We know from \eqref{eq:def-eta} and Lemma \ref{lem:structure-g'-gen} that
$\overline{\eta}_{a_l}=\overline{\gamma}_le^{|a_l|}$ for some  $\overline{\gamma}_l\in \Re$, $l=1,\ldots,r$ and $\overline{\eta}_{r+1}=\overline{\gamma}_{r+1}e^{|b|}$ for some $\overline{\gamma}_{r+1}\in\Re$, where $ e^{p}$ is the vector  of all ones in $\Re^{p}$ and $\overline{\cal C}=g'(\overline{\sigma})-{\rm Diag}(\overline{\eta})\in\Re^{m\times m}$ has the following form
\begin{equation}\label{eq:structure-g'}
\overline{\cal C}=
\left[\begin{array}{cccc}
\overline{c}_{11}E^{|a_{1}||a_{1}|} & \cdots & \overline{c}_{1r}E^{|a_{1}||a_{r}|} & 0 \\
\vdots & \ddots & \vdots & \vdots \\
\overline{c}_{r1}E^{|a_{r}||a_{1}|} & \cdots & \overline{c}_{rr}E^{|a_{r}||a_{r}|} & 0 \\
0 & \cdots & 0 & 0
\end{array}\right]\,,
\end{equation}
where $E^{pq}\in\Re^{p\times q}$ is the $p$ by $q$ matrix of all ones and $\bar{c}\in\Re^{r\times r}$. Then we know from (\ref{eq:eig-singular-F-diff-1-a-added}) that
\[
 \left(g'(\overline{{\sigma}})h\right)_{a_l}=\left\{\begin{array}{ll}
\left(\overline{\cal C}h\right)_{a_l}+\overline{\gamma}_l\lambda(S(H_{a_l a_l})) & \mbox{if $l=1,\dots,r$,} \\[3pt]
\overline{\gamma}_{r+1}\sigma([H_{bb}\quad H_{bc}]) & \mbox{if $l=r+1$}\, ,
 \end{array} \right.
  \]
where for $l\in\{1,\ldots,r\}$, $\left(\overline{\cal C}h\right)_{a_l}={\sum_{l'=1}^{r}}\bar{c}_{ll'}{\rm tr}(S(H_{a_{l'}a_{l'}})) e^{|a_l|}=\left(\overline{\cal C}{\rm diag}(S(H_1))\right)_{a_l}$.
On the other hand, we know from  (\ref{eq:diag-H-decomp-ak}), (\ref{eq:diag-H-decomp-b}), (\ref{eq:eig-singular-F-diff-1}) and (\ref{eq:eig-singular-F-diff-1-a-added}) that for $H$ sufficiently close to $0$, and $l=1,\ldots,r$,
\begin{eqnarray*}
S(H_{a_{l}a_{l}}) &=& Q_{l}(\Sigma(X)_{a_{l}a_{l}}-\overline{\nu}_{l}I_{|a_{l}|})Q_{l}^\T +O(\|H\|^{2})
= Q_{l}\Lambda(S(H_{a_{l}a_{l}}))Q_{l}^\T +O(\|H\|^{2}),  \\[3pt]
\left[H_{bb}\quad H_{bc}\right] &=& M(\Sigma(X)_{bb}-\overline{\nu}_{r+1}I_{|b|})N_{1}^\T +O(\|H\|^{2})
=M\Sigma(\left[H_{bb}\quad H_{bc}\right]) N_1+O(\|H\|^2)\,.
\end{eqnarray*}
Therefore, from (\ref{eq:eig-singular-F-diff-1-a-added}),  (\ref{eq:Delta-k-F-diff-1}) and (\ref{eq:Delta-k-F-diff-2}), we obtain that
\begin{eqnarray*}
\Delta_l(H)&=&\left[\begin{array}{ccc} 0 & 0  & 0  \\[3pt] 0 & {\rm Diag}\left((\overline{\cal C}{\rm diag}(S(H_1)))_{a_l}\right)+\bar{\gamma}_l S(H_{a_la_l})  & 0   \\[3pt] 0 & 0  & 0 \end{array}\right]
+o(\|H\|),\quad l=1,\ldots,r\,,
\\[3pt]
\Delta_{r+1}(H)&=&\left[\begin{array}{ccc}0 & 0  & 0 \\0 & \overline{\gamma}_{r+1}H_{bb} & \overline{\gamma}_{r+1}H_{bc}  \end{array}\right]+o(\|H\|)\,.
\end{eqnarray*}
Thus, we know from (\ref{eq:GR-def-F-diff-spectral-op}) that for any $H$ sufficiently close to $0$,
\begin{equation}\label{eq:nonsmoothpart-F-diff-Spectral-op-diag}
G_{R}(X)=\left[ {\rm Diag}\left(\overline{\cal C}{\rm diag}(S(H_1))\right)\quad 0\right] +\left[\begin{array}{ccccc}\overline{\gamma}_{1} S(H_{a_1a_1})  & 0 & 0 & 0 & 0 \\0 & \ddots & 0 & 0 & 0 \\0 & 0 & \overline{\gamma}_{r} S(H_{a_ra_r})  & 0 & 0 \\0 & 0 & 0 & \overline{\gamma}_{r+1}H_{bb} & \overline{\gamma}_{r+1}H_{bc}\end{array}\right]+o(\|H\|)\,.
\end{equation}

Next, consider the general $\overline{X}\in\V^{m\times n}$. For any $H\in\V^{m\times n}$ sufficiently close to $0$, rewrite (\ref{eq:def-Y-Z-F-diff}) as
\[
[\Sigma(\overline{X})\quad 0 ]+\overline{U}^\T H\overline{V}=\overline{U}^\T U[\Sigma(X)\quad 0]V^\T \overline{V}\,.
\] Denote $\widetilde{U}:=\overline{U}^\T U$ and $\widetilde{V}:=\overline{V}^\T V$. Let $\widetilde{X}:=[\Sigma(\overline{X})\quad 0 ]+\overline{U}^\T H\overline{V}$. Then, since $ \overline{U}$ and $ \overline{V}$ are unitary  matrices, we know from (\ref{eq:nonsmoothpart-F-diff-Spectral-op-diag}) that
\begin{eqnarray}
G_R(X)&=&\overline{U} G_{R} (\widetilde{{X}}) \overline{V}^\T=\overline{U}\left[{\rm Diag}\left(\overline{\cal C}{\rm diag}(S( {A}))\right)\quad 0\right]\overline{V}^\T \nonumber\\
&&+\overline{U}\left[\begin{array}{ccccc}\overline{\gamma}_{1} S( {A}_{a_1a_1})  & 0 & 0 & 0 & 0 \\0 & \ddots & 0 & 0 & 0 \\0 & 0 & \overline{\gamma}_{r} S( {A}_{a_ra_r})  & 0 & 0 \\0 & 0 & 0 & \overline{\gamma}_{r+1} {A}_{bb} & \overline{\gamma}_{r+1} {B}_{bc}\end{array}\right]\overline{V}^\T+o(\|H\|).
\quad
\label{eq:nonsmoothpart-F-diff-Spectral-op}
\end{eqnarray}
Thus, by combining (\ref{eq:smoothpart-F-diff-Spectral-op}) and (\ref{eq:nonsmoothpart-F-diff-Spectral-op}) with (\ref{eq:Gs'-formula}) and noting that $G(\overline{X})=G_{S}(\overline{X})$, we obtain that for any  $H\in \V^{m\times n}$ sufficiently small,
\[
G(X)-G(\overline{X})=\overline{U}[\overline{\cal E}_{1}\circ S( {A})+{\rm Diag}\left(\overline{\cal C}{\rm diag}(S( {A}))\right)+\overline{\cal E}_{2}\circ T( {A}) \quad \overline{\cal F}\circ  {B} ]\overline{V}^\T +o(\|H\|)\,.
\] Therefore, we know that $G$ is F-differentiable at $\overline{X}$ and (\ref{eq:F-derivative-spectral-op}) holds.

$``\Longrightarrow"$ Suppose that  $G$ is F-differentiable at $\overline{X}$. Let $(\overline{U},\overline{V})\in{\mathbb O}^{m\times n}(\overline{X})$ be fixed.  For any $h\in\Re^{m}$, let $H=\overline{U}[{\rm Diag}(h)\quad 0]\overline{V}^\T \in\V^{m\times n}$. We know from Proposition \ref{prop:G-g} that for all $h$ sufficiently close to $0$,
$G(\overline{X}+H)=\overline{U}{\rm Diag}(g(\overline{\sigma}+h))\overline{V}_1^\T$. Therefore, we know from the assumption that for all $h$ sufficiently close to $0$,
 \[
{\rm Diag}(g(\overline{\sigma}+h)-g(\overline{\sigma}))= \overline{U}^\T \left(G(\overline{X}+H)-G(\overline{X})\right)\overline{V}_1=\overline{U}^\T G'(\overline{X})H \overline{V}_1+o(\|h\|)\,.
 \]
Hence, we know that $g$ is F-differentiable at $\overline{\sigma}$ and ${\rm Diag}( g'(\sigma)h)= \overline{U}^\T G'(\overline{X})H \overline{V}_1$. The proof  is competed.$\hfill \Box$

\vskip 10 true pt

Finally, we can present the continuous differentiability result of spectral operators  in the following theorem.

\begin{theorem}\label{thm:con-diff-spectral-op}
Suppose that $\overline{X}\in{\cal N}$ has the SVD (\ref{eq:Y-eig-Z-SVD}). Then, $G$ is continuously differentiable at $\overline{X}$ if and only if $g$ is continuously differentiable at $\overline{\sigma}=\sigma(\overline{X})$.
\end{theorem}
\noindent {\bf  Proof.} $``\Longleftarrow"$
By the assumption, we know from Theorem \ref{thm:diff-spectral-op} that there exists an open neighborhood ${\cal B}\subseteq{\cal N}$ of $\overline{X}$ such that the spectral operator $G$ is differentiable on ${\cal B}$, and for any $X\in{\cal B}$, the derivative $G'(X)$ is given by
\begin{equation}\label{eq:derv-Y-sep-op}
G'(X)H= U[{\cal E}_{1}\circ S(A)+ {\rm Diag}\left({\cal C}{\rm diag}(S(A))\right)+{\cal E}_{2}\circ T(A)
\quad {\cal F}\circ B ]V^\T  \quad \forall \, H\in\V^{m\times n}\,,
\end{equation} where $(U,V)\in{\mathbb O}^{m,n}(X)$, $A=U^\T HV_1$, $B=U^\T HV_2$ and $\eta$, ${\cal E}_{1}$, ${\cal E}_{2}$, ${\cal F}$ and  $\cal C$ are defined by \eqref{eq:def-eta}--\eqref{eq:def-matric-C} with respect to $\sigma=\sigma(X)$, respectively.   Next, we shall prove that
\begin{equation}\label{eq:conti-diff-sp-op}
\lim_{X\to \overline{X}}G'(X)H\to G'(\overline{X})H \quad \forall\, H\in\V^{m\times n}\,.
\end{equation}

Firstly, we will show that (\ref{eq:conti-diff-sp-op}) holds for the special case that $\overline{X}=[\Sigma(\overline{X})\quad 0]$ and $X=[\Sigma(X)\quad 0]\to\overline{X}$. Let $\{F^{(ij)}\}$ be the standard  basis of $\V^{m\times n}$, i.e., for each $i\in\{1,\ldots,m\}$ and $j\in\{1,\ldots,n\}$,
$F^{(ij)}\in\V^{m\times n}$ is a matrix whose entries are zeros, except the $(i,j)$-th entry is $1$ or $\sqrt{-1}$. Therefore, we only need to show (\ref{eq:conti-diff-sp-op}) holds for all $F^{(ij)}$. Note that since $\sigma(\cdot)$ is globally Lipschitz continuous, we know that for $X$ sufficiently close to $\overline{X}$, $\sigma_{i}\neq \sigma_{j}$ if $\overline{\sigma}_i\neq \overline{\sigma}_j$.

For each $i\in\{1,\ldots,m\}$ and $j\in\{1,\ldots,n\}$, write $F^{(ij)}$ in the following form
\[
F^{(ij)}=\big[F^{(ij)}_{1} \quad  F^{(ij)}_{2} \big]
\] with $F^{(ij)}_{1}\in\V^{m\times m}$ and $F^{(ij)}_{2}\in\V^{m\times(n-m)}$. Next, we consider the following cases.

{\bf  Case 1:} $i,j\in\{1,\ldots,m\}$ and $i=j$. In this case, since $g'$ is continuous at $\overline{\sigma}$, we know that if $F^{(ij)}$ is real, then
\[
\lim_{X\to \overline{X}}G'(X)F^{(ij)}=\lim_{X\to \overline{X}}\left[{\rm Diag}(g'(\sigma)e_{i})\quad 0  \right]=\left[{\rm Diag}(g'(\overline{\sigma})e_{i})\quad 0  \right]=G'(\overline{X})F^{(ij)}\,,
\]
where $e_{i}$ is the vector whose $i$-th entry is one, and zero otherwise; if $F^{(ij)}$ is complex, then
\[
\lim_{X\to \overline{X}}G'(X)F^{(ij)}=\lim_{X\to \overline{X}} \left[\displaystyle{\frac{g_{i}(\sigma)+g_{j}(\sigma)}{\sigma_{i}+\sigma_{j}}}T(F^{(ij)}_{1}) \quad 0\right] =\left[\displaystyle{\frac{g_{i}(\overline{\sigma})+g_{j}(\overline{\sigma})}{\overline{\sigma_{i}}+\overline{\sigma_{j}}}}T(F^{(ij)}_{1}) \quad 0\right] =G'(\overline{X})F^{(ij)}\,.
\]

{\bf  Case 2:} $i,j\in\{1,\ldots,m\}$, $i\neq j$, $\sigma_{i}=\sigma_{j}$ and $\overline{\sigma}_{i}=\overline{\sigma}_{j}>0$. Therefore, we know that there exists $l\in\{1,\ldots,r\}$ such that $i,j\in a_{l}$. Since $g'$ is continuous at $\overline{\sigma}$, we know from (\ref{eq:def-eta}) that
\begin{eqnarray*}
\lim_{X\to \overline{X}}G'(X)F^{(ij)}&=&\lim_{X\to \overline{X}} \left[\left((g'(\sigma))_{ii}-(g'(\sigma))_{ij}\right)S(F^{(ij)}_{1})+\displaystyle{\frac{g_{i}(\sigma)+g_{j}(\sigma)}{\sigma_{i}(X)+\sigma_{j}(X)}}T(F^{(ij)}_{1}) \quad 0 \right]\\
&=& \left[\left((g'(\overline{\sigma}))_{ii}-(g'(\overline{\sigma}))_{ij}\right)S(F^{(ij)}_{1})+\displaystyle{\frac{g_{i}(\overline{\sigma})+g_{j}(\overline{\sigma})}{\overline{\sigma}_{i}+\overline{\sigma}_{j}}}T(F^{(ij)}_{1}) \quad 0 \right]\\
&=&G'(\overline{X})F^{(ij)}\,.
\end{eqnarray*}

{\bf  Case 3:} $i,j\in\{1,\ldots,m\}$, $i\neq j$, $\sigma_{i}\neq\sigma_{j}$ and $\overline{\sigma}_{i}=\overline{\sigma}_{j}>0$. In this case, we know that
\[
G'(X)F^{(ij)}=\left[\displaystyle{\frac{g_{i}(\sigma)-g_{j}(\sigma)}{\sigma_{i}-\sigma_{j}}}S(F^{(ij)}_{1})+\displaystyle{\frac{g_{i}(\sigma)+g_{j}(\sigma)}{\sigma_{i}+\sigma_{j}}}T(F^{(ij)}_{1}) \quad 0 \right]\,.
\]
Let $s,t\in\Re^{m}$ be  two vectors defined by
\begin{equation}\label{eq:def-s-t}
s_{p}:=\left\{ \begin{array}{ll} \sigma_{p} & \mbox{if $p\neq i$,}\\
\sigma_{j} & \mbox{if $p=i$}
\end{array} \right. \quad {\rm and} \quad t_{p}:=\left\{ \begin{array}{ll}\sigma_{p}  & \mbox{if $p\neq i,j$,}\\
\sigma_{j}  & \mbox{if $p=i$},\\
\sigma_{i}   & \mbox{if $p=j$},
\end{array} \right. \quad p\in\{1,\ldots,m\}\,.
\end{equation}
It is clear that both $s$ and $t$ converge to $\overline{\sigma}$ as $X\to \overline{X}$. By noting that $g$ is absolutely symmetric on $\hat{\sigma}_{ {\cal N}}$,  we know from \eqref{eq:def-symmetric} that $g_{j}(\sigma)=g_{i}(t)$, since the vector $t$ is obtained from $\sigma$ by swapping the $i$-th and the $j$-th components. By the mean value theorem (cf. e.g., \cite[Page 68-69]{ORheinboldt70}), we have
\begin{eqnarray}
\frac{g_{i}(\sigma)-g_{j}(\sigma)}{\sigma_{i}-\sigma_{j}}&=&\frac{g_{i}(\sigma)-g_{i}(s)+g_{i}(s)-g_{j}(\sigma)}{\sigma_{i}-\sigma_{j}}=\frac{\displaystyle{\frac{\partial g_{i}(\xi)}{\partial\mu_{i}}}(\sigma_{i}-\sigma_{j})+g_{i}(s)-g_{j}(\sigma)}{\sigma_{i}-\sigma_{j}}\nonumber\\
&=&\frac{\partial g_{i}(\xi)}{\partial\mu_{i}}+\frac{g_{i}(s)-g_{i}(t)+g_{i}(t)-g_{j}(\sigma)}{\sigma_{i}-\sigma_{j}}\nonumber\\
&=&\frac{\partial g_{i}(\xi)}{\partial\mu_{i}}+\frac{\displaystyle{\frac{\partial g_{i}(\hat{\xi})}{\partial\mu_{j}}}(\sigma_{j}-\sigma_{i})+g_{i}(t)-g_{j}(\sigma)}{\sigma_{i}-\sigma_{j}}=\frac{\partial g_{i}(\xi)}{\partial\mu_{i}}-\frac{\partial g_{i}(\hat{\xi})}{\partial\mu_{j}}\,,\label{eq:mean-value-1}
\end{eqnarray}
where $\xi\in\Re^{m}$ lies between $\sigma$ and $s$ and $\hat{\xi}\in\Re^{m}$ is between $s$ and $t$. Consequently, we have $\xi\to\overline{\sigma}$ and $\widehat{\xi}\to\overline{\sigma}$ as $X\to \overline{X}$. By the continuity of $g'$, we have
\[
\lim_{X\to \overline{X}}\frac{g_{i}(\sigma)-g_{j}(\sigma)}{\sigma_{i}-\sigma_{j}}=(g'(\overline{\sigma}))_{ii}-(g'(\overline{\sigma}))_{ij}\quad {\rm and} \quad \lim_{X\to \overline{X}}\frac{g_{i}(\sigma)+g_{j}(\sigma)}{\sigma_{i}+\sigma_{j}}=\frac{g_{i}(\overline{\sigma})+g_{j}(\overline{\sigma})}{\overline{\sigma}_{i}+\overline{\sigma}_{j}}\,.
\] Therefore, we have
\[
\lim_{X\to \overline{X}}G'(X)F^{(ij)}=\left[\left((g'(\overline{\sigma}))_{ii}-(g'(\overline{\sigma}))_{ij}\right)S(F^{(ij)}_{1})+\displaystyle{\frac{g_{i}(\overline{\sigma})+g_{j}(\overline{\sigma})}{\overline{\sigma}_{i}+\overline{\sigma}_{j}}}T(F^{(ij)}_{1}) \quad 0 \right]=G'(\overline{X})F^{(ij)}\,.
\]

{\bf  Case 4:} $i,j\in\{1,\ldots,m\}$, $i\neq j$, $\sigma_{i}>0$ or $\sigma_{j}>0$ and $\overline{\sigma}_{i}\neq\overline{\sigma}_{j}$. Then, we have $\sigma_{i}>0$ or $\sigma_{j}>0$ and $\sigma_{i}\neq\sigma_{j}$. Since $g'$ is continuous at $\overline{\sigma}$, we know that
\begin{eqnarray*}
\lim_{X\to \overline{X}}G'(X)F^{(ij)}&=&\lim_{X\to \overline{X}} \left[\displaystyle{\frac{g_{i}(\sigma)-g_{j}(\sigma)}{\sigma_{i}-\sigma_{j}}}S(F^{(ij)}_{1})+\displaystyle{\frac{g_{i}(\sigma)+g_{j}(\sigma)}{\sigma_{i}+\sigma_{j}}}T(F^{(ij)}_{1}) \quad 0 \right]\\
&=&\left[\displaystyle{\frac{g_{i}(\overline{\sigma})-g_{j}(\overline{\sigma})}{\overline{\sigma}_{i}-\overline{\sigma}_{j}}}S(F^{(ij)}_{1})+\displaystyle{\frac{g_{i}(\overline{\sigma})+g_{j}(\overline{\sigma})}{\overline{\sigma}_{i}+\overline{\sigma}_{j}}}T(F^{(ij)}_{1}) \quad 0 \right]=G'(\overline{X})F^{(ij)}\,.
\end{eqnarray*}

{\bf  Case 5:} $j\in\{m+1,\ldots,n\}$ and $\overline{\sigma}_{i}>0$. Since $g'$ is continuous at $\overline{\sigma}$, we obtain that
\[
\lim_{X\to \overline{X}}G'(X)F^{(ij)}=\lim_{X\to \overline{X}}\left[0 \quad \displaystyle{\frac{g_{i}(\sigma)}{\sigma_{i}}}F^{(ij)}_{2} \right]=\left[0 \quad \displaystyle{\frac{g_{i}(\overline{\sigma})}{\overline{\sigma}_{i}}}F^{(ij)}_{2} \right]=G'(\overline{X})F^{(ij)}\,.
\]

{\bf  Case 6:} $i,j\in\{1,\ldots,m\}$, $i\neq j$, $\overline{\sigma}_{i}=\overline{\sigma}_{j}=0$ and $\sigma_{i}=\sigma_{j}>0$. Therefore, we know that
\[
G'(X)F^{(ij)}=\left[\left((g'(\sigma))_{ii}-(g'(\sigma))_{ij}\right)S(F^{(ij)}_{1})+\displaystyle{\frac{g_{i}(\sigma)+g_{j}(\sigma)}{\sigma_{i}+\sigma_{j}}}T(F^{(ij)}_{1}) \quad 0 \right]\,.
\]
We know from \eqref{eq:def-eta} and Lemma \ref{lem:structure-g'-gen} that
\begin{equation}\label{eq:limit-case6-1}
\lim_{X\to \overline{X}}(g'(\sigma))_{ii}=(g'(\overline{\sigma}))_{ii}=\overline{\eta}_{i}\quad {\rm and} \quad \lim_{X\to \overline{X}}(g'(\sigma))_{ij}=(g'(\overline{\sigma}))_{ij}= 0\,.
\end{equation}
Let $\hat{s},\hat{t}\in\Re^{m}$ be  two vectors defined by
\begin{equation}\label{eq:def-hat-s-t}
\hat{s}_{p}:=\left\{ \begin{array}{ll} \sigma_{p} & \mbox{if $p\neq i$,}\\
-\sigma_{j} & \mbox{if $p=i$}
\end{array} \right. \quad {\rm and} \quad \hat{t}_{p}:=\left\{ \begin{array}{ll}\sigma_{p}  & \mbox{if $p\neq i,j$,}\\
-\sigma_{j}  & \mbox{if $p=i$}\,,\\
-\sigma_{i}   & \mbox{if $p=j$}\,,
\end{array} \right. \quad p\in\{1,\ldots,m\}\,.
\end{equation}
Also, it clear that both $\hat{s}$ and $\hat{t}$ converge to $\overline{\sigma}$ as $X\to \overline{X}$. Again, by noting that $g$ is absolutely symmetric on $\hat{\sigma}_{ {\cal N}}$, we know from \eqref{eq:def-symmetric} that
\[
g_{i}(\sigma)=-g_{j}(\hat{t})\quad {\rm and} \quad g_{j}(\sigma)=-g_{i}(\hat{t})\,.
\]
By using similar arguments for deriving  (\ref{eq:mean-value-1}), we have
\begin{eqnarray}
\frac{g_{i}(\sigma)+g_{j}(\sigma)}{\sigma_{i}+\sigma_{j}}
=\frac{\partial g_{i}(\zeta)}{\partial\mu_{i}}+\frac{\partial g_{i}(\hat{\zeta})}{\partial\mu_{j}}\,,\label{eq:mean-value-2}
\end{eqnarray} where $\zeta\in\Re^{m}$ is between $\sigma$ and $\hat{s}$ and $\hat{\zeta}\in\Re^{m}$ is between $\hat{s}$ and $\hat{t}$. Consequently, we know that $\zeta,\hat{\zeta}\to\overline{\sigma}$ as $X\to \overline{X}$. By the continuity of $g'$, we know from \eqref{eq:def-eta} that
\begin{equation}\label{eq:limit-case6-2}
\lim_{X\to \overline{X}}\frac{g_{i}(\sigma)+g_{j}(\sigma)}{\sigma_{i}+\sigma_{j}}=(g'(\overline{\sigma}))_{ii}=\overline{\eta}_{i}\,.
\end{equation} Therefore, from (\ref{eq:limit-case6-1}) and (\ref{eq:limit-case6-2}), we have
\[
\lim_{X\to \overline{X}}G'(X)F^{(ij)}=\left[\overline{\eta}_{i}F^{(ij)}_{1} \quad 0 \right]=G'(\overline{\sigma})F^{(ij)}\,.
\]

{\bf  Case 7:} $i,j\in\{1,\ldots,m\}$, $i\neq j$, $\overline{\sigma}_{i}=\overline{\sigma}_{j}=0$, $\sigma_{i}\neq\sigma_{j}$ and $\sigma_{i}>0$ or $\sigma_{j}>0$.
Let $s$, $t$ and $\hat{s}$, $\hat{t}$ be defined by \eqref{eq:def-s-t} and (\ref{eq:def-hat-s-t}), respectively. By the continuity of $g'$, we know from (\ref{eq:mean-value-1}) and (\ref{eq:mean-value-2}) that
\begin{eqnarray*}
\lim_{X\to \overline{X}}G'(X)F^{(ij)}&=&\lim_{X\to \overline{X}} \left[\displaystyle{\frac{g_{i}(\sigma)-g_{j}(\sigma)}{\sigma_{i}-\sigma_{j}}}S(F^{(ij)}_{1})+\displaystyle{\frac{g_{i}(\sigma)+g_{j}(\sigma)}{\sigma_{i}+\sigma_{j}}}T(F^{(ij)}_{1}) \quad 0 \right]\\
&=&\left[ \overline{\eta}_{i}S(F^{(ij)}_{1})+\overline{\eta}_{i}T(F^{(ij)}_{1}) \quad 0 \right]=\left[\overline{\eta}_{i}F^{(ij)}_{1} \quad 0 \right]=G'(\overline{X})F^{(ij)}\,.
\end{eqnarray*}

{\bf  Case 8:} $i\neq j\in\{1,\ldots,m\}$, $\overline{\sigma}_{i}=\overline{\sigma}_{j}=0$ and $\sigma_{i}=\sigma_{j}=0$. By the continuity of $g'$, we obtain that
\[
\lim_{X\to \overline{X}}G'(X)F^{(ij)}=\lim_{X\to \overline{X}} \left[(g'(\sigma))_{ii}F^{(ij)}_{1} \quad 0 \right]=\left[(g'(\overline{\sigma}))_{ii}F^{(ij)}_{1} \quad 0 \right]=\left[\overline{\eta}_{i}F^{(ij)}_{1} \quad 0 \right]=G'(\overline{X})F^{(ij)}\,.
\]

{\bf  Case 9:} $j\in\{m+1,\ldots,n\}$, $\overline{\sigma}_{i}=0$ and $\sigma_{i}>0$. We know that
\[
G'(X)F^{(ij)}=\left[0 \quad \displaystyle{\frac{g_{i}(\sigma)}{\sigma_{i}}}F^{(ij)}_{2} \right]\,.
\]
Let $\tilde{s}\in\Re^{m}$ be a vector given by
\[
\tilde{s}_{p}:=\left\{ \begin{array}{ll} \sigma_{p} & \mbox{if $p\neq i$,}\\
0 & \mbox{if $p=i$,}
\end{array} \right. \quad p\in\{1,\ldots,m\}\,.
\]
 Therefore, we have $\tilde{s}$ converges to $\overline{\sigma}$ as $X\to\overline{X}$. Since $g$ is absolutely symmetric on $\hat{\sigma}_{ {\cal N}}$, we know that $g_{i}(\tilde{s})=0$. Also, by the mean value theorem, we have
\[
\frac{g_{i}(\sigma)}{\sigma_{i}}=\frac{g_{i}(\sigma)-g_{i}(\tilde{s})}{\sigma_{i}}=\frac{\partial g_{i}(\rho)}{\partial\mu_{i}}\,,
\] where $\rho\in\Re^{m}$ is between $\sigma$ and $\tilde{s}$. Consequently, we have $\rho$ converges to $\overline{\sigma}$ as $X\to \overline{X}$. By the continuity of $g'$, we know from \eqref{eq:def-eta} that
\[
\lim_{X\to \overline{X}}\frac{g_{i}(\sigma)}{\sigma_{i}}=(g'(\overline{\sigma}))_{ii}=\overline{\eta}_{i}\,.
\] Thus,
\[
\lim_{X\to \overline{X}}G'(X)F^{(ij)}=\lim_{X\to \overline{X}} \left[0 \quad \displaystyle{\frac{g_{i}(\sigma)}{\sigma_{i}}}F^{(ij)}_{2} \right]=\left[0 \quad \overline{\eta}_{i}F^{(ij)}_{2} \right]=G'(\overline{X})F^{(ij)}\,.
\]

{\bf  Case 10:} $j\in\{m+1,\ldots,n\}$, $\overline{\sigma}_{i}=0$ and $\sigma_{i}=0$. By the continuity of $g'$, we know that
\[
\lim_{X\to \overline{X}}G'(X)F^{(ij)}=\lim_{X\to \overline{X}} \left[0 \quad (g'(\sigma))_{ii}F^{(ij)}_{2} \right]=\left[0 \quad (g'(\overline{\sigma}))_{ii}F^{(ij)}_{2} \right]=G'(\overline{X})F^{(ij)}\,.
\]

Finally, we consider the general case that
\[
X=U\left[\Sigma(X)\quad 0 \right]V^\T \quad {\rm and}\quad \overline{X}= \overline{U}\left[\Sigma(\overline{X})\quad 0\right]\overline{V}^\T \,.
\]
 By noting from Theorem \ref{thm:diff-spectral-op} that $G$ is F-differential at $X$ if and only if $G$ is F-differential at $\left[\Sigma (X) \quad  0\right]$ and for any $H\in\V^{m\times n}$,
 \[
    G^\prime(X) H   =U \left ( G^\prime(\left[\Sigma (X) \quad  0\right]) (U^THV)\right) V^T\, ,
 \]
we know from the above analysis that $G$ is continuously differentiable at $\overline{X}$.

$``\Longrightarrow"$
Suppose that $G$ is continuously differentiable at $\overline{X}$. Let $(\overline{U},\overline{V})\in{\mathbb O}^{m\times n}(\overline{X})$ be fixed. For any $\sigma\in\Re^{m}$, define $X:=\overline{U}[{\rm Diag}(\sigma)\quad 0]\overline{V}^\T$. For any $h\in\Re^{m}$, let $H:=\overline{U}[{\rm Diag}(h)\quad 0]\overline{V}^\T$. From the proof of the second part of Theorem \ref{thm:diff-spectral-op}, we know from the assumption that for all $\sigma$ sufficiently close to $\overline{\sigma}$,
\[
{\rm Diag}( g'(\sigma)h) =\overline{U}^{\T}(G'(X)H)\overline{V}_{1},\quad h\in\Re^{m}\,.
\]
Consequently,  $g$ is also continuously differentiable at $\overline{\sigma}$.
$\hfill \Box$

\medskip
\begin{remark}
{
In order to compute \eqref{eq:F-derivative-spectral-op}, it appears that one needs to
compute and store $\overline{V}_2\in \V^{n\times (n-m)}$ explicitly, which would  incur huge memory cost if $n \gg m$. Fortunately, due to the special form of $\overline{\cal F}$, the 
explicit computation of $\overline{V}_2$ can be avoided as we shall show next. 
Let $\bar{f}=(\bar{f}_1,\dots,\bar{f}_m)^T$ be defined by
$$
\bar{f}_i = \left\{\begin{array}{ll}
\displaystyle{{g_{i}(\bar{\sigma})}/{\bar{\sigma}_{i}}} & \mbox{if $\bar{\sigma}_{i}\neq 0$}\,,\\[3pt]
(g'(\bar{\sigma}))_{ii} & \mbox{otherwise.}
\end{array}\right.
$$
Observe
that the term in \eqref{eq:F-derivative-spectral-op} involving $\overline{V}_2$ is given by
\begin{eqnarray*}
\overline{U} (\overline{\cal F}\circ (\overline{U}^\T H\overline{V}_{2}) )\overline{V}_2^\T
= \overline{U} {\rm Diag}(\bar{f})  \overline{U}^\T H \overline{V}_2\overline{V}_2^\T
=  \overline{U} {\rm Diag}(\bar{f})  \overline{U}^\T H(I_n-\overline{V}_1\overline{V}_1^\T)
=  \overline{U} {\rm Diag}(\bar{f})  \overline{U}^\T (H - (H\overline{V}_1)\overline{V}_1^\T).
\end{eqnarray*}
Thus in numerical implementation, the large matrix $\overline{V}_2$ 
is not needed. 
}
\end{remark}

\section{Lipschitz continuity, Bouligand differentiability, G-semismoothness, and Clarke's
generalized Jacobian}
\label{section:semismoothness}

\subsection{Lipschitz continuity}

In this subsection, we analyze the local Lipschitz continuity of  the spectral operator $G$ defined on a nonempty set ${\cal N}$. Let  $\overline{X}\in{\cal N}$ be given. Assume that $g$ is locally Lipschitz continuous near $\overline{\sigma}=\sigma(\overline{X})$ with module $L>0$. Therefore, there exists a positive constant $\delta_{0}>0$ such that
\[
\|g(\sigma)-g(\sigma')\|\leq L\|\sigma-\sigma'\| \quad \forall\, \sigma,\sigma'\in B(\overline{\sigma},\delta_{0}):=\left\{y\in\hat{\sigma}_{ {\cal N}}\mid \|y-\overline{\sigma}\|\leq \delta_0\right\}\,.
\]
By using the absolutely symmetric property of $g$ on $\hat{\sigma}_{ {\cal N}}$, we obtain the following simple proposition.

\begin{proposition}\label{prop:gigj}
There exist a positive constant $L'>0$ and a positive constant $\delta>0$ such that for any $\sigma\in B(\overline{\sigma},\delta)$,
\begin{eqnarray}
|g_{i}(\sigma)-g_{j}(\sigma)|&\leq& L'|\sigma_{i}-\sigma_{j}|\quad \forall\, i,j\in\{1,\ldots, m\},\ i\neq j,\;\; \sigma_{i}\neq \sigma_{j}\,,\label{eq:gigj-1}\\
|g_{i}(\sigma)+g_{j}(\sigma)|&\leq& L'|\sigma_{i}+\sigma_{j}|\quad \forall\, i,j\in\{1,\ldots, m\},\;\;  \sigma_{i}+\sigma_{j}>0\,,\label{eq:gigj-2}\\
|g_{i}(\sigma)|&\leq& L'|\sigma_{i}|\quad \forall\, i\in\{1,\ldots,m\},\;\; \sigma_{i}>0\label{eq:gigj-3}\,.
\end{eqnarray}
\end{proposition}
\noindent {\bf  Proof.}
It is easy to check that there exists a positive constant $\delta_{1}>0$ such that for any $\sigma\in B(\overline{\sigma},\delta_{1})$,
\begin{eqnarray}
|\sigma_{i}-\sigma_{j}| &\geq &\delta_{1}>0\quad \forall\, i,j\in\{1,\ldots,m\},\ i\neq j,
\;\; \overline{\sigma}_{i}\neq\overline{\sigma}_{j}\,,
\label{eq:bounded-cond-1}
\\[3pt]
|\sigma_{i}+\sigma_{j}|&\geq& \delta_{1}>0\quad \forall\, i,j\in\{1,\ldots,m\},\;\; \overline{\sigma}_{i}+\overline{\sigma}_{j}>0\,,
\label{eq:bounded-cond-2}
\\[3pt]
|\sigma_{i}|&\geq& \delta_{1}>0\quad \forall\, i\in\{1,\ldots,m\},\;\;
 \overline{\sigma}_{i}>0\,.
\label{eq:bounded-cond-3}
\end{eqnarray}
Let $\delta:=\min\{\delta_{0},\delta_{1}\}>0$. Denote $\tau:=\displaystyle{\max_{i,j}}\{|g_{i}(\overline{\sigma})-g_{j}(\overline{\sigma})|, |g_{i}(\overline{\sigma})+g_{j}(\overline{\sigma})|, |g_{i}(\overline{\sigma})|\}\ge 0$, $L_{1}:=(2L\delta+\tau)/\delta$ and $L':=\max\{L_{1},\sqrt{2}L\}$. Let $\sigma$ be any fixed vector in $B(\overline{\sigma},\delta)$.

Firstly, we consider the case that $i,j \in\{1,\ldots,m\}$, $i\neq j$ and $\sigma_{i}\neq\sigma_{j}$. If $\overline{\sigma}_{i}\neq \overline{\sigma}_{j}$, then from (\ref{eq:bounded-cond-1}), we know that
\begin{eqnarray}
|g_{i}(\sigma)-g_{j}(\sigma)|&=&|g_{i}(\sigma)-g_{i}(\overline{\sigma})+g_{i}(\overline{\sigma})-g_{j}(\overline{\sigma})+g_{j}(\overline{\sigma})-g_{j}(\sigma)|\nonumber\\
&\leq&2\|g(\sigma)-g(\overline{\sigma})\|+\tau \leq \frac{2L\delta+\tau}{\delta}|\sigma_{i}-\sigma_{j}|= L_{1}|\sigma_{i}-\sigma_{j}| \,.\label{eq:bouded-1-1}
\end{eqnarray}
If $\overline{\sigma}_{i}= \overline{\sigma}_{j}$, define $t\in\Re^{m}$ by
\[
t_{p}:=\left\{ \begin{array}{ll}\sigma_{p}  & \mbox{if $p\neq i,j$,}\\
\sigma_{j}  & \mbox{if $p=i$},\\
\sigma_{i}   & \mbox{if $p=j$},
\end{array} \right. \quad p=1,\ldots,m\,.
\] Then, we have $\|t-\overline{\sigma}\|=\|\sigma-\overline{\sigma}\|\leq \delta$. Moreover, since $g$ is absolutely symmetric on $\hat{\sigma}_{ {\cal N}}$, we have $g_{i}(t)=g_{j}(\sigma)$. Therefore
\begin{equation}\label{eq:bouded-1-2}
|g_{i}(\sigma)-g_{j}(\sigma)|=|g_{i}(\sigma)-g_{i}(t)| \leq \|g(\sigma)-g(t)\|\leq L\|\sigma-t\|=\sqrt{2}L|\sigma_{i}-\sigma_{j}|\,.
\end{equation}
Thus, the inequality (\ref{eq:gigj-1}) follows from (\ref{eq:bouded-1-1}) and (\ref{eq:bouded-1-2}) immediately.

Secondly, consider the case $i,j\in\{1,\ldots, m \}$ and $\sigma_{i}+\sigma_{j}>0$. If $\overline{\sigma}_{i}+\overline{\sigma}_{j}>0$, it follows from (\ref{eq:bounded-cond-2}) that
\begin{eqnarray}
|g_{i}(\sigma)+g_{j}(\sigma)|&=&|g_{i}(\sigma)-g_{i}(\overline{\sigma})+g_{i}(\overline{\sigma})+g_{j}(\overline{\sigma})-g_{j}(\overline{\sigma})+g_{j}(\sigma)|\nonumber\\
&\leq& 2\|g(\sigma)-g(\overline{\sigma})\|+\tau \leq\frac{2L\delta+\tau}{\delta}|\sigma_{i}+\sigma_{j}|=L_{1}|\sigma_{i}+\sigma_{j}| \,.\label{eq:bouded-2-1}
\end{eqnarray} If $\overline{\sigma}_{i}+\overline{\sigma}_{j}=0$, i.e., $\overline{\sigma}_{i}=\overline{\sigma}_{j}=0$, define the vector $\hat{t}\in\Re^{m}$ by
\[
\hat{t}_{p}:=\left\{ \begin{array}{ll}\sigma_{p}  & \mbox{if $p\neq i,j$,}\\
-\sigma_{j}  & \mbox{if $p=i$},\\
-\sigma_{i}   & \mbox{if $p=j$},
\end{array} \right. \quad p=1,\ldots,m\,.
\] By noting that $\overline{\sigma}_{i}=\overline{\sigma}_{j}=0$, we obtain that $\|\hat{t}-\overline{\sigma}\|=\|\sigma-\overline{\sigma}\|\leq \delta$.
Again, since $g$ is absolutely symmetric on $\hat{\sigma}_{ {\cal N}}$, we have $g_{i}(t)=g_{j}(\sigma)$, we have $g_{i}(\hat{t})=-g_{j}(\sigma)$. Therefore,
\begin{equation}\label{eq:bouded-2-2}
|g_{i}(\sigma)+g_{j}(\sigma)|=|g_{i}(\sigma)-g_{i}(\hat{t})|\leq \|g(\sigma)-g(\hat{t})\|\leq L\| \sigma-\hat{t}\|=\sqrt{2}L|\sigma_{i}+\sigma_{j}|\,.
\end{equation}
Thus the inequality (\ref{eq:gigj-2}) follows from  (\ref{eq:bouded-2-1}) and (\ref{eq:bouded-2-2}).

Finally, we consider the case that $i\in\{1,\ldots,m \}$ and $\sigma_{i}>0$\,. If $\overline{\sigma}_{i}>0$, then we know from (\ref{eq:bounded-cond-3}) that
\begin{eqnarray}
|g_{i}(\sigma)|&=&|g_{i}(\sigma)-g_{i}(\overline{\sigma})+g_{i}(\overline{\sigma})|\leq |g_{i}(\sigma)-g_{i}(\overline{\sigma})|+|g_{i}(\overline{\sigma})|\nonumber\\
&\leq&\|g(\sigma)-g(\overline{\sigma})\|+\tau\leq \frac{2L\delta+\tau}{\delta}|\sigma_{i}|\leq L_{1}|\sigma_{i}|\,.\label{eq:bouded-3-1}
\end{eqnarray} If $\overline{\sigma}_{i}=0$, define $s\in\Re^{m}$  by
\[
s_{p}:=\left\{ \begin{array}{ll}\sigma_{p}  & \mbox{if $p\neq i$,}\\
0  & \mbox{if $p=i$,}
\end{array} \right. \quad p=1,\ldots,m\,.
\]
Then, since $\sigma_{i}>0$, we know that $\|s-\overline{\sigma}\|<\|\sigma-\overline{\sigma}\|\leq \delta$. Moreover, since $g$ is absolutely symmetric on $\hat{\sigma}_{ {\cal N}}$, we have $g_{i}(t)=g_{j}(\sigma)$, we know that $g_{i}(s)=0$. Therefore, we have
\begin{equation}\label{eq:bouded-3-2}
|g_{i}(\sigma)|=|g_{i}(\sigma)-g_{i}(s)|\leq \|g(\sigma)-g(s)\|\leq L\|\sigma-s\|\leq L|\sigma_{i}|\,.
\end{equation} Thus,  the inequality (\ref{eq:gigj-1}) follows from (\ref{eq:bouded-3-1}) and (\ref{eq:bouded-3-2}) immediately. This completes the proof.
$\hfill \Box$
\vskip 10 true pt

For any fixed $0<\omega\leq\delta_{0}/\sqrt{m}$ and $y\in B(\overline{\sigma},\delta_{0}/(2\sqrt{m})):=\{\|y-\overline{\sigma}\|_{\infty}\leq\delta_{0}/(2\sqrt{m})\}$, the function $g$ is integrable on $V_{\omega}(y):=\{z\in\Re^{m}\,|\,\|y-z\|_{\infty}\leq \omega/2\}$ (in the sense of Lebesgue). Therefore, we know that the function
\begin{equation}\label{eq:def-Steklov averaged function}
g(\omega,y):=\frac{1}{\omega ^{m}}\int_{V_{\omega}(y)}g(z)dz
\end{equation} is well-defined on $(0,\delta_{0}/\sqrt{m}\,]\times B(\overline{\sigma},\delta_{0}/(2\sqrt{m}))$ and is said to be Steklov averaged function \cite{Steklov1907} of $g$. For the sake of  convenience, we always define $g(0,y)=g(y)$. Since $g$ is absolutely symmetric on $\hat{\sigma}_{ {\cal N}}$, we have $g_{i}(t)=g_{j}(\sigma)$, it is easy to check that for each fixed $0<\omega\leq\delta_{0}/\sqrt{m}$, the function $g(\omega,\cdot)$ is also absolutely symmetric on $B(\overline{\sigma},\delta_{0}/(2\sqrt{m}))$. By the definition, we know that $g(\cdot,\cdot)$ is locally Lipschitz continuous on $(0,\delta_{0}/\sqrt{m}\,]\times B(\overline{\sigma},\delta_{0}/(2\sqrt{m}))$ with the module $L$.  Meanwhile, by elementary calculations, we know that $g(\cdot,\cdot)$ is continuously differentiable on $(0,\delta_{0}/\sqrt{m}\,]\times B(\overline{\sigma},\delta_{0}/(2\sqrt{m}))$ and for any  fixed $\omega\in(0,\delta_{0}/\sqrt{m}\,]$ and $y\in B(\overline{\sigma},\delta_{0}/(2\sqrt{m}))$,
\begin{equation*}\label{eq:g'-uni-bounded}
\|g'_{y}(\omega,y)\|\leq L\,.
\end{equation*} Moreover, we know that $g(\omega,\cdot)$ converges to $g$ uniformly on the compact set $B(\overline{\sigma},\delta_{0}/(2\sqrt{m}))$ as $\omega\downarrow 0$. By using the formula (\ref{eq:F-derivative-spectral-op}), the following results can be obtained from Theorem \ref{thm:con-diff-spectral-op} and Proposition \ref{prop:gigj} directly.

\begin{proposition}\label{prop:G'-uni-bounded}
Suppose that $g$ is locally Lipschitz continuous near $\overline{\sigma}$, Let $g(\cdot,\cdot)$ be the corresponding Steklov averaged function defined in  (\ref{eq:def-Steklov averaged function}). Then, for any given $\omega\in(0,\delta_{0}/\sqrt{m}\,]$, the spectral operator $G(\omega,\cdot)$ with respect to $g(\omega,\cdot)$ is continuously differentiable on $B(\overline{X},\delta_{0}/(2\sqrt{m})):=\{X\in{\cal X}\,|\, \|\sigma(X)-\overline{\sigma}\|_{\infty}\leq\delta_{0}/(2\sqrt{m})\}$, and there exist two positive constants $\delta_{1}>0$ and $\overline{L}>0$ such that
\begin{equation}\label{eq:G'-uni-bounded}
\|G'(\omega,X)\|\leq \overline{L}\quad \forall\, 0<\omega\leq \min\{\delta_{0}/\sqrt{m},\delta_{1}\}\ {\rm and}\ X\in B(\overline{X},\delta_{0}/(2\sqrt{m}))\,.
\end{equation} Moreover, $G(\omega,\cdot)$ converges to $G$ uniformly in the compact set $B(\overline{X},\delta_{0}/(2\sqrt{m}))$ as $\omega\downarrow 0$.
 \end{proposition}

Proposition \ref{prop:G'-uni-bounded} allows us to derive the following result on the local Lipschitz continuity of  spectral operators.

\begin{theorem}\label{thm:Lip-spectral-op}
Suppose that $\overline{X}$ has the SVD (\ref{eq:Y-eig-Z-SVD}).  The spectral operator $G$ is locally Lipschitz continuous near $\overline{X}$ if and only if $g$ is locally Lipschitz continuous near $\overline{\sigma}=\sigma(\overline{X})$.
\end{theorem}
\noindent {\bf  Proof.} $``\Longleftarrow"$ Suppose that $g$ is locally Lipschitz continuous near $\overline{\sigma}=\sigma(\overline{X})$ with module $L>0$, i.e., there exists a positive constant $\delta_{0}>0$ such that
\[
\|g(\sigma)-g(\sigma')\|\leq L\|\sigma-\sigma'\| \quad \forall\, \sigma,\sigma'\in B(\overline{\sigma},\delta_{0})\,.
\] By Proposition \ref{prop:G'-uni-bounded}, for any  $\omega\in(0,\delta_{0}/\sqrt{m}\,]$, the spectral operator $G(\omega,\cdot)$ defined with respect to the Steklov averaged function $g(\omega, \cdot)$ is continuously differentiable. Since $G(\omega,\cdot)$ converges to $G$ uniformly in the compact set $B(\overline{X},\delta_{0}/(2\sqrt{m}))$ as $\omega\downarrow 0$, we know that for any $\varepsilon>0$, there exists a constant $\delta_{2}>0$ such that for any $0<\omega\leq\delta_{2}$,
\[
\|G(\omega,X)-G(X)\|\leq \varepsilon\quad \forall\,X\in B(\overline{X},\delta_{0}/(2\sqrt{m}))\,.
\] Fix any $X,X'\in B(\overline{X},\delta_{0}/(2\sqrt{m}))$ with $X\neq X'$. By Proposition \ref{prop:G'-uni-bounded}, we know that there exists $\delta_{1}>0$ such that (\ref{eq:G'-uni-bounded}) holds. Let $\bar{\delta}:=\min\{\delta_{1}, \delta_{2},\delta_{0}/\sqrt{m}\}$. Then, by the mean value theorem, we know that
\begin{eqnarray*}
\|G(X)-G(X')\|&=&\|G(X)-G(\omega,X)+G(\omega,X)-G(\omega,X')+G(\omega,X')-G(X')\|\\
&\leq&2\varepsilon+\|\int_{0}^{1}G'(\omega,X+t(X-X'))dt\|\\
&\leq&\overline{L}\|X-X'\|+2\varepsilon\quad \forall\, 0<\omega<\bar{\delta}\,.
\end{eqnarray*} Since $X,X'\in B(\overline{X},\delta_{0}/(2\sqrt{m}))$ and $\varepsilon>0$ are arbitrary, by letting $\varepsilon\downarrow 0$, we obtain that
\[
\|G(X)-G(X')\|\leq\overline{L}\|X-X'\|\quad \forall\,X,X'\in B(\overline{X},\delta_{0}/(2\sqrt{m}))\,.
\] Thus $G$ is locally Lipschitz continuous near $\overline{X}$.

$``\Longrightarrow"$ Suppose that $G$ is locally Lipschitz continuous near $\overline{X}$ with module $L>0$, i.e., there exists an open neighborhood ${\cal B}$ of $\overline{X}$ in ${\cal N}$ such that for any $X,X'\in{\cal B}$,
\[
\|G(X)-G(X')\|\leq L\|X-X'\|\,.
\]
Let $(\overline{U},\overline{V})\in{\mathbb O}^{m\times n}(\overline{X})$ be fixed. For any $y\in\hat{\sigma}_{\cal N}$, we   define $Y:=\overline{U}\left[{\rm Diag}(y)\quad 0 \right]\overline{V}^\T$. Then, we know from Proposition \ref{prop:G-g} that $G(Y)=\overline{U}\left[{\rm Diag}(g(y))\quad 0 \right]\overline{V}^\T$. Therefore, we obtain that there exists  an open neighborhood ${\cal B}_{\overline{\sigma}}$ of $\overline{\sigma}$ in $\hat{\sigma}_{\cal N}$ such that
\[
\|g(y)-g(y')\|=\|G(Y)-G(Y')\|\leq L\|Y-Y'\|=L\|y-y'\|\quad \forall\, y,y'\in{\cal B}_{\overline{\sigma}}\,.
\]This completes the proof.
$\hfill \Box$

\subsection{Bouligand-differentiability}

In this section, we study the $\rho$-order Bouligand-differentiability of $G$ with $0<\rho\leq 1$, which is stronger than the directional differentiability. Let ${\cal Z}$ be a finite dimensional real Euclidean space equipped with an inner product $\langle \cdot,\cdot \rangle$ and its induced norm $\|\cdot\|$. Let ${\cal O}$ be an open set in ${\cal Z}$ and ${\cal Z}'$ be another finite dimensional real Euclidean space. The function $F: {\cal O}\subseteq{\cal Z}  \to  {\cal Z}^\prime$ is said to be  {\it B(ouligand)-differentiable} at $z\in{\cal O}$ if for any $h\in{\cal Z}$ with $h\to 0$,
\begin{equation*}\label{eq:def-B-diff}
F(z+h)-F(z)-F'(z;h)=o(\|h\|)\,.
\end{equation*}
A stronger notion than B-differentiability is $\rho$-order B-differentiability with $\rho>0$. The function $F: {\cal O}\subseteq{\cal Z}  \to  {\cal Z}^\prime$ is said to be  {\it $\rho$-order B-differentiable} at $z\in{\cal O}$ if for any $h\in{\cal Z}$ with $h\to 0$,
\begin{equation*}\label{eq:def-rho-order-B-diff}
F(z+h)-F(z)-F'(z;h)=O(\|h\|^{1+\rho})\,.
\end{equation*}

Let $\overline{X}\in\V^{m\times n}$ be given. We have the following results on the $\rho$-order B-differentiability of spectral operators.
\begin{theorem}\label{thm:rho-order-B-diff-spectral-op}
Suppose that $\overline{X}\in{\cal N}$ has the SVD (\ref{eq:Y-eig-Z-SVD}).  Let $0<\rho\leq 1$ be given. Then,
\begin{itemize}
\item[(i)] if $g$ is locally Lipschitz continuous near $\sigma(\overline{X})$ and $\rho$-order B-differentiable at $\sigma(\overline{X})$, then $G$ is $\rho$-order B-differentiable at $\overline{X}$;
\item[(ii)] if $G$ is $\rho$-order B-differentiable at $\overline{X}$, then $g$ is $\rho$-order B-differentiable at $\sigma(\overline{X})$.
\end{itemize}
\end{theorem}
\noindent {\bf  Proof.} Without loss of generality, we only consider the case that $\rho =1$.

(i) For any $H\in\V^{m\times n}$, denote $X=\overline{X}+H$. Let $U\in{\mathbb O}^{m}$ and $V\in{\mathbb O}^{n}$ be such that
\begin{equation}\label{eq:def-Y-Z-B-diff}
X=U[\Sigma(X)\quad 0]V^\T \,.
\end{equation}
Denote $\sigma=\sigma(X)$. Let $G_{S}(X)$ and $G_{R}(X)$ be defined by (\ref{eq:def-GS}). Therefore, by Lemma \ref{lem:Gs-diff-formula}, we know that for any $H\to{0}$,
\begin{equation}\label{eq:smoothpart-B-diff-Spectral-op}
G_{S}(X)-G_{S}(\overline{X})=G'_{S}(\overline{X})H+O(\|H\|^{2})=G_S'(\overline{X})H+O(\|H\|^{2})\,,
\end{equation}
where $G_S'(\overline{X})H$ is given by \eqref{eq:Gs'-formula}. For $H\in\V^{m\times n}$ sufficiently small, we have ${\cal U}_{l}(X)={\sum_{i\in a_{l}}}u_{i}v_{i}^\T $, $l=1,\ldots,r$. Therefore,  we know that
\begin{equation}\label{eq:GR-def-B-diff-spectral-op}
G_{R}(X)=G(X)-G_{S}(X)=\sum_{l=1}^{r+1}\Delta_{l}(H)\,,
\end{equation}
where
\[
\Delta_{l}(H)=\sum_{i\in a_{l}}(g_{i}(\sigma)-g_i(\overline{\sigma}))u_{i}v_{i}^\T  \quad l=1,\ldots,r \quad {\rm and} \quad \Delta_{r+1}(H)=\sum_{i\in b}g_i(\sigma)u_{i}v_{i}^\T \,.
\]

We first consider the case that $\overline{X}=[\Sigma(\overline{X})\quad 0 ]$. Then,  we know from (\ref{eq:directi-diff}) and (\ref{eq:directi-diff-sigma}) that for any  $H$ sufficiently small,
\begin{equation}\label{eq:eig-singular-B-diff-1}
\sigma=\overline{\sigma}+\sigma'(\overline{X};H)+O(\|H\|^{2})\,,
\end{equation}
where $\sigma'(\overline{X};H)=\left(\lambda(S(H_{a_{1}a_{1}})),\ldots,\lambda(S(B_{a_{r}a_{r}})),\sigma([H_{bb}\quad H_{bc}])\right)\in\Re^{m}$. Denote $h:=\sigma'(\overline{X};H)$. Since $g$ is locally Lipschitz continuous near $\overline{\sigma}$ and $1$-order B-differentiable at $\overline{\sigma}$, we know that for any $H$ sufficiently small,
\begin{equation*}
g(\sigma)-g(\overline{\sigma})=g(\sigma+h+O(\|H\|^{2}))-g(\overline{\sigma})= g(\sigma+h)-g(\overline{\sigma})+O(\|H\|^{2})=g'(\overline{\sigma};h)+O(\|H\|^{2})\,.
\end{equation*}
Let $\phi=g'(\overline{\sigma};\cdot)$. Since $u_{i}v_{i}^\T $, $i=1,\ldots,m$ are uniformly bounded, we obtain that for $H$ sufficiently small,
\begin{eqnarray*}
\Delta_{l}(H) &=& U_{a_{l}}{\rm Diag}(\phi_{l}(h))V_{a_{l}}^\T +O(\|H\|^{2}),\quad  l=1,\ldots,r\,,
\\
\Delta_{r+1}(H)&=& U_{b}{\rm Diag}(\phi_{r+1}(h))V_{b}^\T +O(\|H\|^{2})\,.
\end{eqnarray*}
Again, we know from (\ref{eq:diag-U-V}) that there exist $Q_{l}\in{\mathbb O}^{|a_{l}|}$, $M\in{\mathbb O}^{|b|}$ and $N=[N_{1}\quad N_{2}]\in{\mathbb O}^{n-|a|}$ with $N_{1}\in\V^{(n-|a|)\times |b|}$ and $N_{2}\in\V^{(n-|a|)\times(n-m)}$ (depending on $H$) such that
\begin{eqnarray*}
U_{a_{l}} &=&\left[\begin{array}{c} O(\|H\|)  \\ Q_{l}+O(\|H\|) \\ O(\|H\|) \end{array}\right], \quad  V_{a_{l}}=\left[\begin{array}{c} O(\|H\|)  \\ Q_{l}+O(\|H\|) \\ O(\|H\|) \end{array}\right], \ l=1,\ldots,r\,,
\\[3pt]
U_{b}&=&\left[\begin{array}{c} O(\|H\|) \\[3pt] M+O(\|H\|)\end{array}\right],\quad  [V_{b}\quad V_{c}]=\left[\begin{array}{c} O(\|H\|) \\[3pt]
N+O(\|H\|)\end{array}\right]  \,.
\end{eqnarray*}
Since $g$ is locally Lipschitz continuous near $\overline{\sigma}$ and directionally differentiable at $\overline{\sigma}$, we know from \cite[Theorem A.2]{Robinson87} or \cite[Lemma 2.2]{QSun93} that the directional derivative $\phi$ is globally Lipschitz continuous on $\Re^{m}$. Thus, for $H$ sufficiently small, we have $\|\phi(h)\|=O(\|H\|)$. Therefore,  we obtain that
\begin{eqnarray}
\Delta_{l}(H) &=& \left[\begin{array}{ccc} 0 & 0  & 0  \\[3pt] 0 & Q_{l}{\rm Diag}(\phi_{l}(h))Q_{l}^\T   & 0   \\[3pt] 0 & 0  & 0 \end{array}\right]
+O(\|H\|^{2}),\quad l=1,\ldots,r\,,
\label{eq:Delta-k-B-diff-1}
\\
\Delta_{r+1}(H)&=&\left[\begin{array}{cc}0 & 0  \\0 & M{\rm Diag}(\phi_{r+1}(h))N_{1}^\T  \end{array}\right]+O(\|H\|^{2})\,.
\label{eq:Delta-k-B-diff-2}
\end{eqnarray}
Again, we know from \eqref{eq:diag-H-decomp-ak} and \eqref{eq:diag-H-decomp-b} that
\begin{eqnarray}
S(H_{a_{l}a_{l}}) &=& Q_{l}(\Sigma(X)_{a_{l}a_{l}}-\overline{\nu}_{l}I_{|a_{l}|})Q_{l}^\T +O(\|H\|^{2}),\quad l=1,\ldots,r\,,
\label{eq:H-H'-2}
\\[3pt]
[H_{bb}\quad H_{bc}] &=& M(\Sigma(X)_{bb}-\overline{\nu}_{r+1}I_{|b|})N_{1}^\T +O(\|H\|^{2})\,.
\label{eq:H-H'-3}
\end{eqnarray}
Since $g$ is locally Lipschitz continuous near $\overline{\sigma}=\sigma(\overline{X})$,  we know from Theorem \ref{thm:Lip-spectral-op} that the spectral operator $G$ is locally Lipschitz continuous near $\overline{X}$. Therefore, we know from Theorem \ref{prop:H-dir-diff-spectr-op} and Remark \ref{remark:Lip-H-diff-spectral-op} that $G$ is directional differentiable at $\overline{X}$. Thus, from  \cite[Theorem A.2]{Robinson87} or \cite[Lemma 2.2]{QSun93}, we know that $G'(\overline{X},\cdot)$ is globally Lipschitz continuous on $\V^{m\times n}$. Thus, the corresponding spectral operator $\Phi $ defined by (\ref{eq:dir-derivative-spectral-op}) is globally Lipschitz continuous on ${\cal W}$. Hence, we know from (\ref{eq:GR-def-B-diff-spectral-op}) that for $H$ sufficiently small,
\begin{equation}\label{eq:nonsmoothpart-B-diff-Spectral-op-diag}
G_{R}(X)= {\widehat
\Phi}(D(H))+O(\|H\|^{2}) \, ,
\end{equation} where
$D(H)=\left(S(H_{a_{1}a_{1}}),\ldots,S(H_{a_{r}a_{r}}), H_{b\bar{a}} \right)$ and ${\widehat
\Phi}(\cdot)$ is defined by (\ref{eq:def-Diag-Phi}).

Next, consider the general case that $\overline{X}\in\V^{m\times n}$. For any $H\in\V^{m\times n}$, rewrite (\ref{eq:def-Y-Z-B-diff}) as
\[
[\Sigma(\overline{X})\quad 0 ]+\overline{U}^\T H\overline{V}=\overline{U}^\T U[\Sigma(Z)\quad 0]V^\T \overline{V}\,.
\] Denote $\widetilde{U}:=\overline{U}^\T U$ and $\widetilde{V}:=\overline{V}^\T V$. Let $\widetilde{X}:=[\Sigma(\overline{X})\quad 0 ]+\overline{U}^\T H\overline{V}$. Then, since $ \overline{U}$ and $ \overline{V}$ are unitary matrices, we know from (\ref{eq:nonsmoothpart-B-diff-Spectral-op-diag}) that
\begin{equation}\label{eq:nonsmoothpart-B-diff-Spectral-op}
G_{R}(X)=\overline{U} {\widehat
\Phi}(D(H)) \overline{V}^\T  +O(\|H\|^{2})\,,
\end{equation}
where $D(H)=\left(S(\widetilde{H}_{a_{1}a_{1}}),\ldots,S(\widetilde{H}_{a_{r}a_{r}}), \widetilde{H}_{b\bar{a}} \right)$ and $\widetilde{H}=\overline{U}^\T H\overline{V}$. Thus, by combining (\ref{eq:smoothpart-B-diff-Spectral-op}) and (\ref{eq:nonsmoothpart-B-diff-Spectral-op}) and noting that $G(X)=G_{S}(\overline{X})$, we obtain that for any  $H\in\V^{m\times n}$ sufficiently close to $0$,
\[
G(X)-G(\overline{X})-G'(\overline{X};H)=O(\|H\|^{2})\,,
\] where $G'(\overline{X};H)$ is given by (\ref{eq:dir-diff-spectr-op}). This implies that $G$ is $1$-order B-differentiable at $\overline{X}$.

(ii) Suppose that $G$ is $1$-order B-differentiable at $\overline{X}$. Let $(\overline{U},\overline{V})\in{\mathbb O}^{m\times n}(\overline{X})$ be fixed.  For any $h\in\Re^{m}$, let $H=\overline{U}[{\rm Diag}(h)\quad 0]\overline{V}^\T \in\V^{m\times n}$. We know from Proposition \ref{prop:G-g} that for all $h$ sufficiently close to $0$, $G(\overline{X}+H)=\overline{U}{\rm Diag}(g(\overline{\sigma}+h))\overline{V}_1^\T$. Therefore, we know from the assumption that
 \[
 {\rm Diag}(g(\overline{\sigma}+h)-g(\overline{\sigma}))= \overline{U}^\T \left(G(\overline{X}+H)-G(\overline{X})\right)\overline{V}_1=\overline{U}^\T G'(\overline{X};H)\overline{V}_1+O(\|H\|^{2})\,.
  \]
This shows that $g$ is $1$-order B-differentiable at $\overline{\sigma}$. The proof  is competed.
$\hfill \Box$
\vskip 5 true pt

\subsection{G-semismoothness}

Let ${\cal Z}$ and ${\cal Z}'$ be two finite dimensional real Euclidean spaces and ${\cal O}$ be an open set in ${\cal Z}$. Suppose that $F: {\cal O}  \subseteq {\cal Z}  \to  {\cal Z}^\prime$ is a locally Lipschitz continuous function on ${\cal O}$. Then, according to Rademacher's theorem, $F$ is almost everywhere differentiable (in the sense of Fr\'{e}chet)  in ${\cal O}$.  Let ${\cal D}_F$ be the set of points in ${\cal O}$ where $F$ is differentiable. Let $F^\prime(z)$ be the derivative of $F$ at $z\in {\cal D}_F$. Then the  {\it B(ouligand)-subdifferential} of $F$ at $z\in {\cal O}$ is denoted by \cite{Qi93}:
\[
\partial_{B}F(z):=\left\{ \lim_{{\cal D}_{F}\ni z^{k}\to z} F'(z^{k}) \right\}\,
\] and   {\it Clarke's generalized Jacobian} of $F$ at $z\in {\cal O}$ \cite{Clarke83}  takes the form:
\[
\partial F(z)={\rm conv}\{\partial_{B} F(z)\}\,,
\] where ``conv'' stands for the convex hull in the usual sense of convex analysis \cite{Rockafellar70}. The function $F$ is said to be G-semismooth  at a point $z \in {\cal O}$  if for any $y\to z $ and $V \in \partial F(y)$,
\[
F(y) -F(z) - V(y-z) = o(\|y- z\| )\, .
\]
A stronger notion than G-semismoothness is $\rho$-order G-semismoothness with $\rho>0$. The function $F$ is said to be  $\rho$-order G-semismooth at $z$  if
 for any $y\to z $ and $V \in \partial F(y)$,
\[
F(y) - F(z) - V(y- z) = O(\|y-z\|^{1+\rho})\, .
\]
In particular, the function $F$ is said to be strongly G-semismooth at $z$ if $F$ is $1$-order G-semismooth at $z$. Furthermore, the  function $F$ is said to be  ($\rho$-order, strongly) semismooth at $z\in {\cal O}$  if (i)
the directional derivative of $F$ at $z$ along any direction $d\in {\cal Z}$, denoted by $F^\prime(z;d)$,  exists; and (ii) $F$ is   ($\rho$-order, strongly) G-semismooth.

The following result taken from \cite[Theorem 3.7]{SSun02} provides a convenient tool for proving the G-semismoothness of Lipschitz functions.
\begin{lemma}\label{lem:semismoothness-equiv}
Let $F: {\cal O}  \subseteq {\cal Z} \to {\cal Z}^\prime $ be a locally Lipschitz continuous function on the open set ${\cal O}$. Let  $\rho >0$ be a constant. $F$ is $\rho $-order G-semismooth (G-semismooth) at $z$ if and only if  for any ${\cal D}_F\ni y\to  z $,
\begin{equation}\label{eq:G-semismoothness-equiv}
F(y) - F(z) - F^\prime(y)(y- z) = O(\|y- z\|^{1+\rho}) \quad \big(=o(\|y- z\|)\big) \, .
 \end{equation}
\end{lemma}

Let  $\overline{X}\in{\cal N}$ be given. Assume that $g$ is locally Lipschitz continuous near $\overline{\sigma}=\sigma(\overline{X})$. Thus, from Theorem \ref{thm:Lip-spectral-op} we know that the corresponding spectral operator $G$ is locally Lipschitz continuous near $\overline{X}$. The following theorem is on the G-semismoothness of the spectral operator $G$.

\begin{theorem}\label{thm:rho-order-G-semismooth-spectral-op}
Suppose that $\overline{X}\in{\cal N}$ has the decomposition (\ref{eq:Y-eig-Z-SVD}).  Let $0<\rho\leq 1$ be given. $G$ is $\rho$-order G-semismooth at $\overline{X}$ if and only if $g$ is $\rho$-order G-semismooth at $\overline{\sigma}$.
\end{theorem}
\noindent {\bf  Proof.} Without loss of generality, we only consider the case that $\rho=1$.

$``\Longleftarrow"$ For any $H\in\V^{m\times n}$, denote $X=\overline{X}+H$. Let $U\in{\mathbb O}^{m}$ and $V\in{\mathbb O}^{n}$ be such that
\begin{equation}\label{eq:def-Y-Z-semismooth}
 X=U[\Sigma(X)\quad 0]V^\T \,.
\end{equation}
Denote $\sigma=\sigma(X)$. $G_{S}$ and $G_{R}$ are two mappings defined in (\ref{eq:def-GS}). We know from Lemma \ref{lem:Gs-diff-formula} that there exists an open neighborhood ${\cal B}\subseteq{\cal N}$ of $\overline{X}$ such that $G_{S}$ twice continuously differentiable on ${\cal B}$ and
\begin{eqnarray}
&& G_{S}(X)-G_{S}(\overline{X})=\sum_{l=1}^{r}\bar{g}_{l}\,{\cal U}'_{l}(X)H+O(\|H\|^{2}) \nonumber\\
&=&\sum_{l=1}^{r} \bar{g}_{l}\left\{U[\Gamma_{l}(X)\circ S(U^\T HV_{1} )+\Xi_{l}(X)\circ
T(U^\T HV_{1})]V_{1}^\T +U(\Upsilon_{l}(X)\circ U^\T HV_{2} )V_{2}^\T \right\} +O(\|H\|^{2})\,,\nonumber\\
\label{eq:smoothpart-semismooth-Spectral-op}
\end{eqnarray}
where for each $l\in\{1,\ldots,r\}$, $\Gamma_{l}(X)$, $\Xi_{l}(X)$  and $\Upsilon_{l}(X)$ are given by (\ref{eq:Gamma_k})-(\ref{eq:tildeXi_k}) with $X$, respectively. By taking a smaller ${\cal B}$ if necessary, we assume that for any $X\in{\cal B}$ and $l,l'\in\{1,\ldots,r\}$,
\begin{equation}\label{eq:sigma_X-diff}
\sigma_i(X)>0,\quad \sigma_i(X)\neq \sigma_j(X)\quad \forall\, i\in a_l,\ j\in a_{l'}\ {\rm and}\ l\neq l'\,.
\end{equation} Since $g$ is locally Lipschitz continuous near $\overline{\sigma}$, we know that for any $H$ sufficiently small,
\[
\bar{g}_{l}=g_{i}(\sigma)+O(\|H\|)\quad \forall\, i\in a_l,\quad l=1,\ldots,r\,.
\] Let ${\cal E}_{1}$, ${\cal E}_{2}$ and ${\cal F}$ (depending on $X$) be the matrices defined by (\ref{eq:def-matric-E1})-(\ref{eq:def-matric-F}). By noting that $U\in{\mathbb O}^{m}$ and $V\in{\mathbb O}^{n}$ are uniformly bounded, since $g$ is locally Lipschitz continuous near $\overline{\sigma}$, we know that for any $X\in{\cal B}$ (shrinking ${\cal B}$ if necessary),
\begin{equation}\label{eq:GX-GbarX-semismooth}
G_{S}(X)-G_{S}(\overline{X})=U\left[{\cal E}^0_{1}\circ S(U^\T HV_{1} ) + {\cal E}^0_{2}\circ
T(U^\T HV_{1} )\quad {\cal F}^0\circ U^\T HV_{2} \right]V^\T +O(\|H\|^{2})\,,
\end{equation}
where ${\cal E}^0_{1}$, ${\cal E}^0_{2}$ and ${\cal F}^0$ are the corresponding real matrices defined in \eqref{eq:def-matric-E1}-\eqref{eq:def-matric-F}, respectively.

Let $X\in{\cal D}_{G}\cap{\cal B}$, where ${\cal D}_{G}$ is the set of points in $\V^{m\times n}$ for which $G$ is (F-)differentiable. Define the corresponding index sets in $\{1,\ldots,m\}$ for $X$ by $a':=\{i\mid \sigma_i(X)>0\}$ and $b':=\{i\mid \sigma_i(X)=0\}$. By \eqref{eq:sigma_X-diff}, we have
\begin{equation}
a'\supseteq a \quad {\rm and} \quad  b'\subseteq b\,.
\end{equation}
Let ${\cal E}_{1}$, ${\cal E}_{2}$, ${\cal F}$ and ${\cal C}$ be the corresponding real matrices defined in (\ref{eq:def-matric-ED1})--(\ref{eq:def-matric-C}), respectively. We know from  Theorem \ref{thm:diff-spectral-op} that
\begin{equation}\label{eq:diff-Y-spectral-op-semismooth}
G'(X)H= U[{\cal E}_{1}\circ S(U^\T HV_1)+{\cal E}_{2}\circ T(U^\T HV_1)+ {\rm Diag}\left({\cal C}{\rm diag}(S(U^\T HV_1))\right)
\quad {\cal F}\circ U^\T HV_2 ]V^\T  \,,
\end{equation}
where $\eta$, ${\cal E}_{1}$, ${\cal E}_{2}$, ${\cal F}$ and  $\cal C$ are defined by \eqref{eq:def-eta}--\eqref{eq:def-matric-C} with respect to $\sigma$, respectively.  Denote \[
\Delta(H):=G'(X)H-(G_{S}(X)-G_{S}(\overline{X})).
\]
Moreover, since there exists an integer $j\in\{0,\ldots,|b|\}$ such that $|a'|=|a|+j$, we can define two index sets $b_1:=\{|a|+1,\ldots,|a|+j\}$ and $b_2:=\{|a|+j+1,\ldots,|a|+|b|\}$ such that  $a'=a\cup b_1$ and $b'=b_2$.
From (\ref{eq:GX-GbarX-semismooth}) and (\ref{eq:diff-Y-spectral-op-semismooth}), we obtain that
\begin{equation}\label{eq:Delta2-semismooth-spectral-op}
\Delta(H)=U
{\widehat R} (H)
V^\T +O(\|H\|^{2})\,,
\end{equation}
 where  ${\widehat R} (H)\in \V^{m\times n}$ is defined by
\[
{\widehat R}(H):= \left[\begin{array}{cc} {\rm Diag}\left( R_{1}(H), \dots, R_{r}(H) \right)  & 0
\\
0 &  R_{r+1}(H)
\end{array}\right],
\]
 \begin{eqnarray}\label{eq:R-k-semismooth-spectral-op-2}
 R_{l}(H) &=& ({\cal E}_{1})_{a_{l}a_{l}}\circ S(U_{a_{l}}^\T HV_{a_{l}})+{\rm Diag}\left(({\cal C}{\rm diag}(S(U^\T HV_1)))_{a_la_l}\right), \  l=1,\ldots,r,
\\
 \label{eq:R-k-semismooth-spectral-op-3}
 R_{r+1}(H)&=&\left[\begin{array}{ccc}
({\cal E}_{1})_{b_1b_1}\circ S(U_{b_1}^\T HV_{b_1})+{\rm Diag}\left(({\cal C}{\rm diag}(S(U^\T HV_1)))_{b_1b_1}\right)  & 0 & 0 \\
0 & \gamma U_{b_2}^\T HV_{b_2}  & \gamma U_{b_2}^\T HV_{2}
 \end{array}\right]
\quad
 \end{eqnarray}
and $\gamma:=(g'(\sigma))_{ii}$ for any $i\in b_2$.
By (\ref{eq:Y-eig-Z-SVD}),  we obtain from (\ref{eq:def-Y-Z-semismooth}) that
\[
\left[ \Sigma(\overline{X}) \quad 0\right]+\overline{U}^\T H\overline{V}=\overline{U}^\T U\left[ \Sigma(X)\quad 0 \right]V^\T \overline{V}\,.
\] Let $\widehat{H}:=\overline{U}^\T H\overline{V}$, $\widehat{U}:=\overline{U}^\T U$ and $\widehat{V}:=\overline{V}^\T V$. Then, $U^\T HV=\widehat{U}^\T \overline{U}^\T H\overline{V}\widehat{V}=\widehat{U}^\T \widehat{H}\widehat{V}$. We know from (\ref{eq:diag-U-V}) that there exist $Q_{l}\in{\mathbb O}^{|a_{l}|}$, $l=1,\ldots,r$ and $M\in{\mathbb O}^{|b|}$, $N\in{\mathbb O}^{n-|a|}$  such that
\begin{eqnarray*}
&U_{a_{l}}^\T HV_{a_{l}}=\widehat{U}_{a_{l}}^\T  \widehat{H} \widehat{V}_{a_{l}}=Q_{l}^\T  \widehat{H}_{a_{l}a_{l}}Q_{l}+O(\|H\|^{2}),\quad l=1,\ldots,r\,,&
\\
&\left[ U_{b}^\T HV_{b}\quad U_{b}^\T HV_{2} \right]
 = \left[ \widehat{U}_{b}^\T  \widehat{H} \widehat{V}_{b}\quad \widehat{U}_{b}^\T  \widehat{H} \widehat{V}_{2} \right]=M^\T \left[\widehat{H}_{bb}\quad \widehat{H}_{bc}\right]N+O(\|H\|^{2})\,.&
\end{eqnarray*}
Moreover, from (\ref{eq:diag-H-decomp-ak}) and (\ref{eq:diag-H-decomp-b}), we obtain that
\begin{eqnarray*}
&S(U_{a_{l}}^\T HV_{a_{l}})=Q_{l}^\T  S(\widehat{H}_{a_{l}a_{l}})Q_{l}+O(\|H\|^{2})=\Sigma(X)_{a_{l}a_{l}}-\Sigma(\overline{X})_{a_{l}a_{l}}+O(\|H\|^{2}),\quad l=1,\ldots,r\,,&
\\[3pt]
&\left[ U_{b}^\T HV_{b}\quad U_{b}^\T HV_{2} \right]=M^\T \left[\widehat{H}_{bb}\quad \widehat{H}_{bc}\right]N=\left[\Sigma(X)_{bb}-\Sigma(\overline{X})_{bb} \quad 0\right]+O(\|H\|^{2})\,.&
\end{eqnarray*}
Denote $h=\sigma'(X;H)\in\Re^{m}$. Since the single value functions are strongly semismooth \cite{SSun03}, we know that
\begin{eqnarray*}
& S(U_{a_{l}}^\T HV_{a_{l}})={\rm Diag}(h_{a_{l}})+O(\|H\|^{2}), \quad l=1,\ldots,r\,,&
\\[3pt]
&S(U_{b_1}^\T HV_{b_1})={\rm Diag}(h_{b_1})+O(\|H\|^{2}),\quad \left[ U_{b_2}^\T HV_{b_2}\quad U_{b_2}^\T HV_{2} \right]=\left[ {\rm Diag}(h_{b_2})\quad 0 \right]+O(\|H\|^{2}).
\nonumber
&
\end{eqnarray*}
Therefore, since ${\cal C}=g'(\sigma)-{\rm Diag}(\eta)$, by (\ref{eq:R-k-semismooth-spectral-op-2}) and (\ref{eq:R-k-semismooth-spectral-op-3}), we obtain from (\ref{eq:Delta2-semismooth-spectral-op}) that
\begin{equation}\label{eq:Delta-semismooth-spectral-op}
\Delta(H)=U\left[{\rm Diag}\left(g'(\sigma)h\right)  \quad 0\right]V^\T=U\left[{\rm Diag}\left(g'(\sigma)h\right)  \quad 0\right]V^\T +O(\|H\|^{2})\,.
\end{equation}

On the other hand, for $X$ sufficiently close to $\overline{X}$, we have ${\cal U}_{l}(X)={\sum_{i\in a_{l}}}u_{i}v_{i}^\T $, $l=1,\ldots,r$. Therefore,
\begin{equation}\label{eq:GR-semismooth-spectral-op-1}
G_{R}(X)=G(X)-G_{S}(X)=\sum_{l=1}^{r}\sum_{i\in a_{l}}[g_i(\sigma)-g_i(\overline{\sigma})]u_{i}v_{i}^\T +\sum_{i\in b}g_i(\sigma)u_iv_i^\T \,.
\end{equation} We know from Theorem \ref{thm:diff-spectral-op} that $G$ is differentiable at $X$ if and only if $g$ is differentiable at $\sigma$.  Since $g$ is $1$-order G-semismooth at $\overline{\sigma}$ and $\sigma(\cdot)$ is strongly semismooth, we obtain that for any $X\in{\cal D}_{G}\cap{\cal B}$ (shrinking ${\cal B}$ if necessary),
\begin{equation*}
g(\sigma)-g(\overline{\sigma})=g'(\sigma)(\sigma-\overline{\sigma})+O(\|H\|^{2})=g'(\sigma)(h+O(\|H\|^{2}))+O(\|H\|^{2})=g'(\sigma)h+O(\|H\|^{2})\,.
\end{equation*} Then, since $U\in{\mathbb O}^{m}$ and $U\in{\mathbb O}^{n}$ are uniformly bounded, we obtain from (\ref{eq:GR-semismooth-spectral-op-1}) that
\[
G_{R}(X)=U\left[{\rm Diag}\left(g'(\sigma)h\right)  \quad 0\right]V^\T +O(\|H\|^{2})\,.
\] Thus, from (\ref{eq:Delta-semismooth-spectral-op}), we obtain that $\Delta(H)=G_{R}(X)+O(\|H\|^{2})$. That is, for any $X\in{\cal D}_{G}$ converging to $\overline{X}$,
\[
G(X)-G(\overline{X})-G'(X)H=-\Delta(H)+G_{R}(X)=O(\|H\|^{2})\,.
\]

$``\Longrightarrow"$ Suppose that $G$ is $1$-order G-semismooth at $\overline{X}$.  Let $(\overline{U},\overline{V})\in{\mathbb O}^{m\times n}(\overline{X})$ be fixed.  Assume that $\sigma=\overline{\sigma}+h\in {\cal D}_{g}$ and $h\in\Re^{m}$ is sufficiently small. Let $X=\overline{U}\left[{\rm Diag}(\sigma)\quad 0 \right] \overline{V}^\T $ and $H=\overline{U}\left[{\rm Diag}(h)\quad 0 \right] \overline{V}^\T $. Then, $X\in{\cal D}_{G}$ and converges to $\overline{X}$ if $h$ goes to zero. We know from Proposition \ref{prop:G-g} that for all $h$ sufficiently close to $0$, $G(X)=\overline{U}{\rm Diag}(g({\sigma}))\overline{V}_1^\T$. Therefore,   for any $h$ sufficiently close to $0$,
  \[
 {\rm Diag}(g(\overline{\sigma}+h)-g(\overline{\sigma}))= \overline{U}^\T  \left(G(X)-G(\overline{X})\right)\overline{V}_1=\overline{U}^\T G'(X)H\overline{V}_1+O(\|H\|^{2})\,.
  \]
Hence, since obviously ${\rm Diag}(g'(\sigma)h)=\overline{U}^\T G'(X)H\overline{V}_1$, we know that for $h$ sufficiently small,
 $
 g(\overline{\sigma}+h)-g(\overline{\sigma})=g'(\overline{\sigma})h+O(\|h\|^{2})
 $.  Then, $g$ is $1$-order G-semismooth at $\overline{\sigma}$.
$\hfill \Box$

\subsection{Characterization of Clarke's generalized Jacobian}

Let  $\overline{X}\in{\cal N}$ be given. In this section, we   assume that $g$ is locally Lipschitz continuous near $\overline{\sigma}=\sigma(\overline{X})$ and directionally differentiable at $\overline{\sigma}$.
Therefore, from Theorem \ref{thm:Lip-spectral-op}, Theorem \ref{prop:H-dir-diff-spectr-op} and Remark \ref{remark:Lip-H-diff-spectral-op}, we know that the corresponding spectral operator $G$ is locally Lipschitz continuous near $\overline{X}$ and  directionally differentiable at $\overline{X}$. Furthermore, we define the function $d:\Re^{m}\to\Re^{m}$ by
\begin{equation}\label{eq:def_d}
d(h):=g(\overline{\sigma}+h)-g(\overline{\sigma})-g'(\overline{\sigma};h),\quad h\in\Re^{m}\,.
\end{equation} Thus, since $g$ is absolutely symmetric on the nonempty open set $\hat{\sigma}_{ {\cal N}}$, we know from \eqref{eq:dir-diff-symmetric} and \eqref{eq:block-structure-P_sigma} that $d$ is also a mixed symmetric mapping, with respect to
$ \mathbb{P}^{|a_1|}\times\ldots\times\mathbb{P}^{|a_r|}\times\pm\mathbb{P}^{|b|}$,
 over $\Re^{|a_{1}|}\times \ldots \times  \Re^{|a_{r}|}\times\Re^{|b|}$.
Moreover, since $g$ is locally Lipschitz continuous near $\overline{\sigma}$ and directional differentiable at $\overline{\sigma}$, we know that $g$ is B-differentiable at $\overline{\sigma}$ (cf. \cite{Shapiro90}). Thus,  $d$ is differentiable at zero with the derivative $d'(0)=0$. Furthermore, if we assume that the function $d$ is also strictly differentiable at zero, then we have
\begin{equation}\label{eq:theta-strict-diff}
\lim_{w,w'\to 0\atop w\neq w'}\frac{d(w)-d(w')}{\|w-w'\|}=0\,.
\end{equation} Thus, by using the mixed symmetric property of $d$, one can easily obtain the following results. We omit the details of the proof here.

\begin{lemma}\label{lem:limit-deri-theta}
Let $d:\Re^m\to\Re^m$ be the function given by \eqref{eq:def_d}. Suppose that $d$ is strictly differentiable at zero. Let $\{w^k\}$ be a given sequence in $\Re^m$ converging to zero. Then, if there exist $i,j\in a_{l}$ for some $l\in\{1,\ldots,r\}$ or $i, j\in b$ such that $w^k_i\neq w^k_j$ for all $k$ sufficiently large, then
\begin{equation}\label{eq:limit-deri-theta-1}
\lim_{k\to \infty}\frac{d_i(w^k)-d_j(w^k)}{w^k_i-w^k_j}=0\,;
\end{equation} if there exist $i,j \in b$ such that $w^k_i+w^k_j\neq 0$ for all $k$ sufficiently large, then
\begin{equation}\label{eq:limit-deri-theta-2}
\lim_{k\to \infty}\frac{d_i(w^k)+d_j(w^k)}{w^k_i+w^k_j}=0\,;
\end{equation}
and if there exists $i\in b$ such that $w^k_i\neq 0$ for all $k$ sufficiently large, then
\begin{equation}\label{eq:limit-deri-theta-3}
\lim_{k\to \infty}\frac{d_i(w^k)}{w^k_i}=0\,.
\end{equation}
\end{lemma}

Let $\Psi(\cdot):=G'(\overline{X};\cdot):\V^{m\times n}\to\V^{m\times n}$ be the directional derivative of $G$ at $\overline{X}$. We know from (\ref{eq:dir-diff-spectr-op}) that for any $Z\in\V^{m\times n}$,
\begin{equation}\label{eq:def-Psi}
\Psi(Z) =G'(\overline{X};Z)= \overline{U}\left[\overline{\cal E}^0_1\circ S(\overline{U}^\T Z\overline{V}_{1})+\overline{\cal E}^0_2\circ T(\overline{U}^\T Z\overline{V}_{1})\quad \overline{\cal F}^0\circ \overline{U}^\T Z\overline{V}_{2}  \right]\overline{V}^\T  +\overline{U} {\widehat \Phi}(D(Z))\overline{V}^\T ,
\end{equation}
where $D(Z)=\left(S(\widetilde{Z}_{a_1a_1}),\ldots,S(\widetilde{Z}_{a_ra_r}),\widetilde{Z}_{b\bar{a}}\right)\in{\cal W}$, $\widetilde{Z}=\overline{U}^\T Z\overline{V}$
and ${\widehat \Phi}(\cdot):{\cal W}\to \V^{m\times n}$ is given by (\ref{eq:def-Diag-Phi}) with ${\Phi}(\cdot):{\cal W}\to{\cal W}$ being the spectral operator defined by \eqref{eq:dir-derivative-spectral-op} with respect to the mixed symmetric mapping $\phi(\cdot):=g'(\overline{\sigma};\cdot)$.
 Since the spectral operator $G$ is locally Lipschitz continuous near $\overline{X}$, we know that $\Psi(\cdot)=G'(\overline{X};\cdot)$ is globally Lipschitz continuous (cf. \cite[Theorem A.2]{Robinson87} or \cite[Lemma 2.2]{QSun93}). Therefore, $\partial_B \Psi(0)$ and $\partial \Psi(0)$ are well-defined. Furthermore, we have the following useful results on the characterization of   the B-subdifferential   and Clarke's subdifferential
    of the spectral operator $G$ at $\overline{X}$.

\begin{theorem}\label{thm:B-subdiff-spectral-op}
Suppose that the given $\overline{X}\in{\cal N}$ has the decomposition (\ref{eq:Y-eig-Z-SVD}). Suppose that there exists an open neighborhood ${\cal B}\subseteq\Re^{m}$ of $\overline{\sigma}$ in $\hat{\sigma}_{\cal N}$ such that $g(\cdot)$ is differentiable at $\sigma\in{\cal B}$ if and only if $g'(\overline{\sigma};\cdot)$ is differentiable at $\sigma-\overline{\sigma}$. Assume further that the function $d:\Re^{m}\to\Re^{m}$ defined by \eqref{eq:def_d} is strictly differentiable at zero. Then, we have
\[
\partial_{B}G(\overline{X})=\partial_{B}\Psi(0) \quad {\rm and} \quad \partial G(\overline{X})=\partial\Psi(0)\,.
\]
\end{theorem}
\noindent {\bf  Proof.} We only need to prove the result for the B-subdifferentials. Let ${\cal V}$ be any element of $\partial_{B}G(\overline{X})$. Then, there exists a sequence $\{X^k\}$ in ${\cal D}_{G}$ converging to $\overline{X}$ such that ${\cal V}=\displaystyle{\lim_{k\to\infty}} G'(X^k)$. For each $X^k$, let $U^k\in{\mathbb O}^{m}$ and $V^k\in{\mathbb O}^{n}$ be the  matrices such that
\[
X^k=U^k[\Sigma(X^k)\quad 0](V^k)^\T \,.
\] For each $X^k$, denote $\sigma^k=\sigma(X^k)$. Then, we know from Theorem \ref{thm:diff-spectral-op} that for each $k$, $\sigma_{k}\in{\cal D}_{g}$. For $k$ sufficiently large,  we know from Lemma \ref{lem:Gs-diff-formula} that for each $k$, $G_{S} $ is twice continuously differentiable at $\overline{X}$. Thus,  $\displaystyle{\lim_{k\to\infty}}G'_{S}(X^k)=G'_{S}(\overline{X})$. Hence, we have for any $H\in\V^{m\times n}$,
\begin{equation}\label{eq:limit-GS-spectral-op}
\lim_{k\to \infty}G'_{S}(X^k)H=G'_{S}(\overline{X})H=\overline{U}\left[\overline{\cal E}^0_1\circ S(\overline{U}^\T H\overline{V}_{1})+\overline{\cal E}^0_2\circ T(\overline{U}^\T H\overline{V}_{1})\quad \overline{\cal F}^0\circ \overline{U}^\T H\overline{V}_{2}  \right]\overline{V}^\T \,.
\end{equation}
Moreover, we know that the mapping $G_{R}=G-G_{S}$ is also differentiable at each $X^k$ for $k$ sufficiently large. Therefore, we have
\begin{equation}\label{eq:U-Gs-Gr}
{\cal V}=\lim_{k\to \infty}G'(X^k)=G'_{S}(\overline{X})+\lim_{k\to \infty}G_{R}'(X^k)\,.
\end{equation}
From the continuity of the singular value function $\sigma(\cdot)$, by taking a subsequence if necessary, we assume that for each $X^k$ and $l,l'\in\{1,\ldots,r\}$, $\sigma_i(X^k)>0$, $\sigma_{i}(X^k)\neq \sigma_{j}(X^k)$ for any $i\in a_{l}$, $j\in a_{l'}$ and $l\neq l'$. Since $\{U^{k}\}$ and $\{V^{k}\}$ are uniformly bounded, by taking subsequences if necessary, we may also assume that  $\{U^{k}\}$ and $\{V^{k}\}$ converge and denote the limits by $U^{\infty}\in{\mathbb O}^{m}$ and $V^{\infty}\in{\mathbb O}^{n}$, respectively. It is clear that $(U^{\infty},V^{\infty})\in{\mathbb O}^{m,n}(\overline{X})$. Therefore, we know from Proposition \ref{prop:PSigma=SigmaW} that there exist $Q_{l}\in{\mathbb O}^{|a_{l}|}$, $l=1,\ldots,r$,  $Q'\in{\mathbb O}^{|b|}$ and $Q''\in{\mathbb O}^{n-|a|}$ such that $U^{\infty}=\overline{U}M$ and $V^{\infty}=\overline{V}N$, where $M={\rm Diag}(Q_{1},\ldots,Q_{r},Q')\in{\mathbb O}^{m}$ and $N={\rm Diag}(Q_{1},\ldots,Q_{r},Q'')\in{\mathbb O}^{n}$. Let $H\in\V^{m\times n}$ be arbitrarily given. For each $k$, denote $\widetilde{H}^k:=(U^k)^\T HV^k$. Since $\{(U^k,V^k)\}\in{\mathbb O}^{m,n}(X^k)$ converges to $(U^{\infty},V^{\infty})\in{\mathbb O}^{m,n}(\overline{X})$, we know that $\displaystyle{\lim_{k\to\infty}}\widetilde{H}^k=(U^{\infty})^\T HV^{\infty}$. For the notational simplicity, we denote $\widetilde{H}:=\overline{U}^\T H\overline{V}$ and $\widehat{H}:=(U^{\infty})^\T HV^{\infty}$.

For $k$ sufficiently large, we know from \eqref{eq:U'_k} and \eqref{eq:F-derivative-spectral-op} that for any $H\in\V^{m\times n}$, $G_R'(X^k)H=U^k\Delta^k(V^k)^\T $ with
\[
\Delta^k:=\left[\begin{array}{cc}
{\rm Diag}\left(\Delta_1^k ,\dots,\Delta_{r}^k  \right) & 0\\
0  & \Delta_{r+1}^k
\end{array}\right]\in\V^{m\times n}\,,
\]
where for each $k$, $\Delta_l^k=({\cal E}_{l}(\sigma^k))_{a_{l}a_{l}}\circ S(\widetilde{H}^k_{a_{l}a_{l}})+{\rm Diag}(({\cal C}(\sigma){\rm diag}(S(\widetilde{H}^k)))_{a_{l}})$, $l=1,\ldots,r$,
\[
\Delta_{r+1}^k=\left[({\cal E}_{1}(\sigma^k))_{bb}\circ S(\widetilde{H}^k_{bb})+{\rm Diag}(({\cal C}(\sigma){\rm diag}(S(\widetilde{H}^k)))_{b})+({\cal E}_{2}(\sigma^k))_{bb}\circ T(\widetilde{H}^k_{bb})\quad ({\cal F}_{2}(\sigma^k))_{bc}\circ \widetilde{H}^k_{bc}\right]\,
\]
and ${\cal E}_{1}(\sigma^k)$, ${\cal E}_{2}(\sigma^k)$, ${\cal F}(\sigma^k)$ and ${\cal C}(\sigma^k)$ are defined for $\sigma^k$ by (\ref{eq:def-matric-ED1})-(\ref{eq:def-matric-FD}), respectively. Again, since $\{U^{k}\}$ and $\{V^{k}\}$ are uniformly bounded, we know that
\begin{equation}\label{eq:Gr'H-limit}
\lim_{k\to \infty}G_{R}'(X^k)H=U^{\infty}(\lim_{k\to\infty}\Delta^k)(V^{\infty})^\T =\overline{U}M(\lim_{k\to\infty}\Delta^k)N^\T \overline{V}^\T \,.
\end{equation}

Next, we shall show that ${\cal V}\in\partial_B \Psi(0)$. For each $k$, denote $w^k:=\sigma^k-\overline{\sigma}\in\Re^{m}$. Moreover, for each $k$, we can define $W^k_l:=Q_l{\rm Diag}(w^k_{a_l})Q_l^\T \in\S^{|a_l|}$, $l=1,\ldots,r$ and $W^k_{r+1}:=Q'[{\rm Diag}(w^k_b)\quad 0]Q''^\T \in\V^{|b|\times(n-|a|)}$. Therefore, it is clear that for each $k$, $W^k:=(W^k_1,\ldots,W^k_l,W^k_{r+1})\in{\cal W}$ and $\kappa(W^k)=w^k$. Moreover, since ${\lim_{k\to \infty}}\sigma^k=\overline{\sigma}$, we know that ${\lim_{k\to \infty}}W^k=0$ in ${\cal W}$.
 From the assumption, we know that $\phi(\cdot)=g'(\overline{\sigma};\cdot)$ and $d(\cdot)$ are differentiable at each $w^k$ and $\phi'(w^{k})=g'(\sigma^{k})-d'(w^{k})$ for all $w^{k}$. Since $d$ is strictly differentiable at zero, it can be checked easily that $\lim_{k\to\infty}d'(w^{k})=d'(0)=0$. By taking a subsequence if necessary, we may assume that $\lim_{k\to\infty}g'(\sigma^k)$ exists. Therefore, we have
\begin{equation}\label{eq:lim-g'-phi'}
\lim_{k\to\infty}\phi'(w^k)=\lim_{k\to\infty}g'(\sigma^k)\,.
\end{equation}
Since $\Phi$ is the spectral operator with respect to the mixed symmetric mapping $\phi$, from Theorem \ref{thm:diff-spectral-op-block} in Section \ref{section:extension}\footnote{We could present the results in this subsection after introducing Theorem \ref{thm:diff-spectral-op-block} in Section \ref{section:extension}. We include it here for the sake of readability and notational convenience.} we know that $\Phi$ is differentiable at $W\in{\cal W}$  if and only if $\phi$ is differentiable at $\kappa(W)$.
 Recall that  ${\widehat \Phi}: {\cal W}\to\V^{m\times n}$ is defined by (\ref{eq:def-Diag-Phi}).
Then, for $k$ sufficiently large, $ {\widehat \Phi}$ is differentiable at $W^k$. Moreover, for each $k$, we  define the matrix $C^k\in\V^{m\times n}$ by
\[
C^k=\overline{U} \left[\begin{array}{cc} {\rm Diag}\left( W^k_{1},\dots,W^k_{r}\right) & 0 \\[2mm]
0 & W^k_{r+1}
\end{array}\right]
\overline{V}^\T \,.
\]  Then, we know that for $k$ sufficiently large, $\Psi $ is differentiable at $C^k$ and ${\lim_{k\to \infty}}C^k=0$ in $\V^{m\times n}$.
Thus, we know from \eqref{eq:def-Psi} that for each $k$,
\[
\Psi'(C^k)H=G'_S(\overline{X})H+\overline{U}\left[{\widehat \Phi}'(W^k)D(H)\right]\overline{V}^\T
\quad\forall \; H\in\V^{m\times n}\, ,
\]
where $D(H)=\left(S(\widehat{H}_{a_1a_1}),\ldots,S(\widehat{H}_{a_ra_r}),\widehat{H}_{b\bar{a}}\right)\in{\cal W}$ is defined by \eqref{eq:def-D-mapping} and ${\widehat \Phi^\prime}(W^k)D(H)$ can be derived from  \eqref{eq:F-derivative-spectral-op-Phi}. By comparing with \eqref{eq:U-Gs-Gr} and \eqref{eq:Gr'H-limit}, we know that the conclusion then follows if we show that
\begin{equation}\label{eq:K=limRt-1}
\lim_{k\to\infty}\Delta^k=\lim_{k\to\infty}M^\T {\widehat \Phi}'(W^k)D(H)N\,.
\end{equation}
 For any $(i,j)\in\{1,\ldots,m\}\times\{1,\ldots,n\}$, consider the following cases.

{\bf  Case 1:} $i=j$. It is easy to check that for each $k$,
\[
(\Delta^k)_{ii}=(g'(\sigma^k)h^k)_i \quad {\rm and}\quad  \left(M^\T {\widehat \Phi}'(W^k)D(H)N\right)_{ii}=(\phi'(w^k)\widehat{h})_i\,,
\]
where $h^k=\left({\rm diag}(S(\widetilde{H}^{k}_{aa})),{\rm diag}(\widetilde{H}^{k}_{bb})\right)$ and $\widehat{h}=\left({\rm diag}(S(\widehat{H}_{aa})),{\rm diag}(\widehat{H}_{bb})\right)$. Therefore, we know from \eqref{eq:lim-g'-phi'} that
\[
\lim_{k\to\infty}(\Delta^k)_{ii}=\lim_{k\to\infty}(g'(\sigma^k)h^k)_i=\lim_{k\to\infty}(\phi'(w^k)\widehat{h})_i=\lim_{k\to\infty}\left(M^\T {\widehat \Phi}'(W^k)D(H)N\right)_{ii}\,.
\]

{\bf  Case 2:} $i,j\in a_{l}$ for some $l\in\{1,\ldots,r\}$, $i\neq j$ and $\sigma^k_{i}\neq\sigma^k_{j}$ for $k$ sufficiently large. We obtain that for $k$ sufficiently large,
\begin{eqnarray*}
&(\Delta^k)_{ij}=\frac{g_i(\sigma^k)-g_j(\sigma^k)}{\sigma^k_{i}-\sigma^k_{j}}(S(\widetilde{H}^{k}_{a_la_l}))_{ij},&
\\
&
\left(M^\T {\widehat \Phi}'(W^k)D(H)N\right)_{ij}=\frac{\phi_i(w^k)-\phi_j(w^k)}{w^k_{i}-w^k_{j}}(S(\widehat{H}_{a_la_l}))_{ij} \,.
&
\end{eqnarray*}
Since $\overline{\sigma}_{i}=\overline{\sigma}_{j}$ and $g_{i}(\overline{\sigma})=g_{j}(\overline{\sigma})$, we know that for $k$ sufficiently large,
\begin{eqnarray}
\frac{g_i(\sigma^k)-g_j(\sigma^k)}{\sigma^k_{i}-\sigma^k_{j}}&=&\frac{g_i(\overline{\sigma}+w^k)-g_j(\overline{\sigma}+w^k)}{w^k_i-w^k_j}=\frac{g_i(\overline{\sigma}+w^k)-g_i(\overline{\sigma})+g_j(\overline{\sigma})-g_j(\overline{\sigma}+w^k)}{w^k_i-w^k_j}\nonumber\\
&=&\frac{d_i(w^k)-d_j(w^k)}{w^k_i-w^k_j}+\frac{\phi_i(w^k)-\phi_j(w^k)}{w^k_i-w^k_j}\,.\label{eq:difference-g-theta-phi-1}
\end{eqnarray} Therefore, we know from \eqref{eq:limit-deri-theta-1} that
\[
\lim_{k\to\infty}\frac{g_i(\sigma^k)-g_j(\sigma^k)}{\sigma^k_{i}-\sigma^k_{j}}(S(\widetilde{H}^{k}_{a_la_l}))_{ij}=\lim_{k\to\infty}\frac{\phi_i(w^k)-\phi_j(w^k)}{w^k_{i}-w^k_{j}}(S(\widehat{H}_{a_la_l}))_{ij}\,,
\] which implies $\displaystyle{\lim_{k\to\infty}}(\Delta^k)_{ij}=\displaystyle{\lim_{k\to\infty}}\left(M^\T {\widehat \Phi}'(W^k)D(H)N\right)_{ij}$.

{\bf  Case 3:} $i,j\in a_{l}$ for some $l\in\{1,\ldots,r\}$, $i\neq j$ and $\sigma^k_{i}= \sigma^k_{j}$ for $k$ sufficiently large. We have for $k$ sufficiently large,
\begin{eqnarray*}
&(\Delta^k)_{ij}=\left((g'(\sigma^k))_{ii}-(g'(\sigma^k))_{ij}\right)(S(\widetilde{H}^{k}_{a_la_l}))_{ij},&
\\[3pt]
&\left(M^\T {\widehat \Phi}'(W^k)D(H)N\right)_{ij}=\left((\phi'(w^k))_{ii}-(\phi'(w^k))_{ij}\right)(S(\widehat{H}_{a_la_l}))_{ij} \,.&
\end{eqnarray*}
Therefore, we obtain from \eqref{eq:lim-g'-phi'} that
\[
\lim_{k\to\infty}\left((g'(\sigma^k))_{ii}-(g'(\sigma^k))_{ij}\right)(S(\widetilde{H}^{k}_{a_la_l}))_{ij}=\lim_{k\to\infty}\left((\phi'(w^k))_{ii}-(\phi'(w^k))_{ij}\right)(S(\widehat{H}_{a_la_l}))_{ij}\,.
\]
Thus, we have $\displaystyle{\lim_{k\to\infty}}(\Delta^k)_{ij}=\displaystyle{\lim_{k\to\infty}}\left(M^\T {\widehat \Phi}'(W^k)D(H)N\right)_{ij}$.

{\bf  Case 4:} $i,j\in b$, $i\neq j$ and $\sigma^k_{i}=\sigma^k_{j}>0$ for $k$ sufficiently large.
We have for $k$ large,
\begin{eqnarray*}
 &(\Delta^k)_{ij}=\left((g'(\sigma^k))_{ii}-(g'(\sigma^k))_{ij}\right)(S(\widetilde{H}^{k}_{bb}))_{ij}
 +\frac{g_i(\sigma^k)+g_j(\sigma^k)}{\sigma^k_{i}+\sigma^k_{j}}(T(\widetilde{H}^{k}_{bb}))_{ij},&
\\[3pt]
&\left(M^\T {\widehat \Phi}'(W^k)D(H)N\right)_{ij}=\left((\phi'(w^k))_{ii}-(\phi'(w^k))_{ij}\right)(S(\widehat{H}_{bb}))_{ij}+\frac{\phi_i(w^k)+\phi_j(w^k)}{w^k_{i}+w^k_{j}}(T(\widehat{H}_{bb}))_{ij}\,.&
\end{eqnarray*}
Since $\overline{\sigma}_{i}=\overline{\sigma}_{j}=0$ and $g_i(\overline{\sigma})=g_j(\overline{\sigma})=0$, we get
\begin{eqnarray}
\frac{g_i(\sigma^k)+g_j(\sigma^k)}{\sigma^k_{i}+\sigma^k_{j}}
&=&\frac{d_i(w^k)+d_j(w^k)}{w^k_i+w^k_j}+\frac{\phi_i(w^k)+\phi_j(w^k)}{w^k_i+w^k_j}\,.\label{eq:difference-g-theta-phi-2}
\end{eqnarray}
Therefore, we know from (\ref{eq:limit-deri-theta-2}) and (\ref{eq:lim-g'-phi'}) that $\displaystyle{\lim_{k\to\infty}}(\Delta^k)_{ij}=\displaystyle{\lim_{k\to\infty}}\left(M^\T {\widehat \Phi}'(W^k)D(H)N\right)_{ij}$.

{\bf  Case 5:} $i,j\in b$, $i\neq j$ and $\sigma^k_{i}\neq\sigma^k_{j}$ for $k$ sufficiently large. For large $k$, we have
\begin{eqnarray*}
&(\Delta^k)_{ij}=\frac{g_i(\sigma^k)-g_j(\sigma^k)}{\sigma^k_{i}-\sigma^k_{j}}(S(\widetilde{H}^{k}_{bb}))_{ij}+
\frac{g_i(\sigma^k)+g_j(\sigma^k)}{\sigma^k_{i}+\sigma^k_{j}}(T(\widetilde{H}^{k}_{bb}))_{ij},&
\\[3pt]
&\left(M^\T {\widehat \Phi}'(W^k)D(H)N\right)_{ij}=\frac{\phi_i(w^k)-\phi_j(w^k)}{w^k_{i}-w^k_{j}}(S(\widehat{H}_{bb}))_{ij}+\frac{\phi_i(w^k)+\phi_j(w^k)}{w^k_{i}+w^k_{j}}(T(\widehat{H}_{bb}))_{ij}\,.&
\end{eqnarray*}
Thus, by (\ref{eq:difference-g-theta-phi-1}) and (\ref{eq:difference-g-theta-phi-2}), we know from (\ref{eq:limit-deri-theta-1}) and (\ref{eq:limit-deri-theta-2}) that $\displaystyle{\lim_{k\to\infty}}(\Delta^k)_{ij}=\displaystyle{\lim_{k\to\infty}}\left(M^\T {\widehat \Phi}'(W^k)D(H)N\right)_{ij}$.

{\bf  Case 6:} $i,j\in b$, $i\neq j$ and $\sigma^k_{i}=\sigma^k_{j}=0$ for $k$ sufficiently large. We know for $k$ large,
\begin{eqnarray*}
&(\Delta^k)_{ij}=\left((g'(\sigma^k))_{ii}-(g'(\sigma^k))_{ij}\right)(S(\widetilde{H}^{k}_{bb}))_{ij}+
(g'(\sigma^k))_{ii}(T(\widetilde{H}^{k}_{bb}))_{ij},
&
\\[3pt]
&\left(M^\T {\widehat \Phi}'(W^k)D(H)N\right)_{ij}=\left((\phi'(w^k))_{ii}-(\phi'(w^k))_{ij}\right)(S(\widehat{H}_{bb}))_{ij}+(\phi'(w^k))_{ii}(T(\widehat{H}_{bb}))_{ij}\,.&
\end{eqnarray*}
Again, we obtain from \eqref{eq:lim-g'-phi'} that $\displaystyle{\lim_{k\to\infty}}(\Delta^k)_{ij}=\displaystyle{\lim_{k\to\infty}}\left(M^\T {\widehat \Phi}'(W^k)D(H)N\right)_{ij}$.

{\bf  Case 7:} $i\in b$, $j\in c$ and $\sigma^k_{i}>0$ for $k$ sufficiently large. We have for $k$ sufficiently large,
\begin{equation*}
(\Delta^k)_{ij}=\frac{g_i(\sigma^k)}{\sigma^k_{i}}(\widetilde{H}^{k}_{bc})_{ij},
\quad \left(M^\T {\widehat \Phi}'(W^k)D(H)N\right)_{ij}=\frac{\phi_i(w^k)}{w^k_{i}}(\widehat{H}_{bc})_{ij}\,.
\end{equation*}
Since $\overline{\sigma}_{i}=0$ and $g_{i}(\overline{\sigma})=0$, we get
\[
\frac{g_i(\sigma^k)}{\sigma^k_{i}}=\frac{g_i(\overline{\sigma}+w^k)-g_i(\overline{\sigma})}{w^k_{i}}=\frac{d_i(w^k)}{w^k_{i}}+\frac{\phi_i(w^k)}{w^k_{i}}\,.
\] Therefore,  by (\ref{eq:limit-deri-theta-3}), we obtain that $\displaystyle{\lim_{k\to\infty}}(\Delta^k)_{ij}=\displaystyle{\lim_{k\to\infty}}\left(M^\T {\widehat \Phi}'(W^k)D(H)N\right)_{ij}$.

{\bf  Case 8:} $i\in b$, $j\in c$ and $\sigma^k_{i}=0$ for $k$ sufficiently large. We have
for $k$ sufficiently large,
\begin{eqnarray*}
&(\Delta^k)_{ij}=(g'(\sigma^k))_{ii}(\widetilde{H}^{k}_{bc})_{ij},\quad
\left(M^\T {\widehat \Phi}'(W^k)D(H)N\right)_{ij}=(\phi'(w^k))_{ii}(\widehat{H}_{bc})_{ij}\,.&
\end{eqnarray*}
Therefore, by (\ref{eq:lim-g'-phi'}), we obtain that $\displaystyle{\lim_{k\to\infty}}(\Delta^k)_{ij}=\displaystyle{\lim_{k\to\infty}}\left(M^\T {\widehat \Phi}'(W^k)D(H)N\right)_{ij}$.

Thus, we know that \eqref{eq:K=limRt-1} holds. Therefore, by \eqref{eq:U-Gs-Gr} and \eqref{eq:Gr'H-limit}, we obtain that ${\cal V}\in\partial_{B}\Psi(0)$.

Conversely, suppose that ${\cal V}\in\partial_{B}\Psi(0)$ is arbitrarily chosen. Then, from the definition of $\partial_{B}\Psi(0)$, we know that there exists a sequence $\{C^k\}\subseteq\V^{m\times n}$ converging to zero such that $\Psi$ is differentiable at each $C^k$ and ${\cal V}={\lim_{k\to\infty}}\Psi'(C^k)$. For each $k$, we know from \eqref{eq:def-Psi} that $\Psi$ is differentiable at $C^k$ if and only if the spectral operator $\Phi:{\cal W}\to {\cal W}$ is differentiable at $W^k:=D(C^k)=\left(S(\widetilde{C}^k_{a_{1}a_{1}}),\ldots,S(\widetilde{C}^k_{a_{r}a_{r}}),\widetilde{C}^k_{b\bar{a}}\right)\in{\cal W}$, where for each $k$, $\widetilde{C}^k=\overline{U}^\T C^k\overline{V}$.
Moreover, for each $k$, we have the following decompositions
\begin{eqnarray*}
S(\widetilde{C}^k_{a_{l}a_{l}})=Q^k_{l}\Lambda(S(\widetilde{C}^k_{a_{l}a_{l}}))(Q^k_{l})^\T ,\; l=1,\ldots,r,
\quad
\widetilde{C}^k_{b\bar{a}}={Q'}^k\left[ \Sigma(\widetilde{C}^k_{b\bar{a}})\quad 0 \right]({Q''}^k)^\T \,,
\end{eqnarray*}
where $Q^k_{l}\in{\mathbb O}^{|a_{l}|}$,   ${Q'}^k\in{\mathbb O}^{|b|}$ and ${Q''}^k\in{\mathbb O}^{n-|a|}$. For each $k$,
 let
\begin{eqnarray*}
&w^k:=\left(\lambda(S(\widetilde{C}^k_{a_{1}a_{1}})),\ldots,\lambda(S(\widetilde{C}^k_{a_{r}a_{r}})),
\sigma(\widetilde{C}^k_{b\bar{a}})\right)\in\Re^m,&
\\
&M^k:={\rm Diag}\Big(Q_1^k,\dots,Q_{r}^k,  {Q^\prime}^k\Big) \in{\mathbb O}^m,
\quad
N^k:= {\rm Diag}\Big(Q_1^k,\dots,Q_{r}^k,  {Q^{\prime\prime}}^k\Big)
\in{\mathbb O}^n.&
\end{eqnarray*}
Since $\{M^k\}$ and $\{N^k\}$ are uniformly bounded, by taking subsequences if necessary, we know that there exist $Q_l\in{\mathbb O}^{|a_l|}$, $Q'\in{\mathbb O}^{|b|}$ and $Q''\in{\mathbb O}^{n-|b|}$ such that
\[
\lim_{k\to\infty}M^k=M:=
{\rm Diag}\Big(Q_1,\dots,Q_{r},  {Q^\prime}\Big)
 \quad
\lim_{k\to\infty}N^k=N:=
{\rm Diag}\Big(Q_1,\dots,Q_{r},  {Q^{\prime\prime}}\Big).
\]
For each $k$, by  \eqref{eq:F-derivative-spectral-op-Phi} (in Section \ref{section:extension}), we know that for any $H\in\V^{m\times n}$,
\begin{equation}\label{eq:Psi'CkH}
\Psi'(C^k)H = \overline{U}\left[\overline{\cal E}^0_1\circ S(\overline{U}^\T H\overline{V}_{1})+\overline{\cal E}^0_2\circ T(\overline{U}^\T H\overline{V}_{1})\quad \overline{\cal F}^0\circ \overline{U}^\T H\overline{V}_{2}  \right]\overline{V}^\T + \overline{U}\left[{\widehat \Phi}'(W^k)D(H)\right]\overline{V}^\T,
\end{equation}
where $D(H)\in{\cal W}$ is defined by \eqref{eq:def-D-mapping}. Let $R^k: =\Phi'_k(W^k)D(H)$,  $k=1, \ldots,  r+1$.

For each $k$, define $\sigma^k:=\overline{\sigma}+w^k\in\Re^m$.
Since ${\lim_{k\to\infty}}w^k=0$ and for each $k$, $w^k_i\geq 0$ for all $i\in b$,
we have $\sigma^k\geq 0$ for $k$ sufficiently large.
Therefore, for $k$ sufficiently large, we are able to define
\[
X^k:=\overline{U}M[{\rm Diag}(\sigma^k)\quad 0]N^\T \overline{V}^\T \in\V^{m\times n}\,.
\]
For simplicity, denote $U=\overline{U}M\in{\mathbb O}^{m}$ and $ V=\overline{V}N\in{\mathbb O}^{n}$. It is clear that the sequence $\{X^k\}$ converges to $\overline{X}$. From the assumption, we know that $g$ is differentiable at each $\sigma^k$ and $d$ is differentiable at each $w^{k}$ with $g'(\sigma^{k})=\phi'(w^{k})+d'(w^{k})$ for all $\sigma^{k}$. Therefore, by Theorem \ref{thm:diff-spectral-op}, we know that $G$ is differentiable at each $X^k$. By taking subsequences if necessary, we may assume that $\lim_{k\to\infty}\phi'(w^{k})$ exists. Thus, since $d$ is strictly differentiable at zero,  we know that \eqref{eq:lim-g'-phi'} holds.
Since the derivative formula \eqref{eq:Gs'-formula} is independent of $(U,V)\in{\mathbb O}^{m,n}(\overline{X})$, we know from \eqref{eq:F-derivative-spectral-op} that for any $H\in\V^{m\times n}$,
\begin{eqnarray}
G'(X^k)H&=&\overline{U}\left[\overline{\cal E}^0_1\circ S(\overline{U}^\T H\overline{V}_{1})+\overline{\cal E}^0_2\circ T(\overline{U}^\T H\overline{V}_{1})\quad \overline{\cal F}^0\circ \overline{U}^\T H\overline{V}_{2}  \right]\overline{V}^\T \nonumber\\
&&+\;\overline{U}
\left[\begin{array}{cc} {\rm Diag}\left(Q_1\Omega_1^kQ_1^\T,\dots,
Q_r\Omega_r^kQ_r^\T\right)  & 0 \\[2mm]
0 & Q'\Omega_{r+1}^kQ''^\T  \end{array}\right]
\overline{V}^\T \, ,
\label{eq:G'm-form-conver}
\end{eqnarray}
where for each $k$, $\Omega_l^k=({\cal E}_{l}(\sigma^k))_{a_{l}a_{l}}\circ S(\widehat{H}_{a_{l}a_{l}})+{\rm Diag}(({\cal C}(\sigma^k){\rm diag}(S(\widehat{H})))_{a_{l}})$, $l=1,\ldots,r$ and
\[
\Omega_{r+1}^k=\left[({\cal E}_{1}(\sigma^k))_{bb}\circ S(\widehat{H}_{bb})+{\rm Diag}(({\cal C}(\sigma^k){\rm diag}(S(\widehat{H})))_{b})+({\cal E}_{2}(\sigma^k))_{bb}\circ T(\widehat{H}_{bb})\quad ({\cal F}_{2}(\sigma^k))_{bc}\circ \widehat{H}_{bc}\right]\,,
\]
  ${\cal E}_{1}(\sigma^k)$, ${\cal E}_{2}(\sigma^k)$ and ${\cal F}(\sigma^k)$ are defined   by (\ref{eq:def-matric-ED1})--(\ref{eq:def-matric-FD}), respectively and $\widehat{H}:=M^\T \overline{U}^\T H\overline{V}N=M^\T \widetilde{H}N$. Therefore, by comparing \eqref{eq:Psi'CkH} and \eqref{eq:G'm-form-conver}, we know that the inclusion ${\cal V}\in \partial_B G(\overline{X})$ follows if we show that
\begin{equation}\label{eq:lim-R-Omega}
\lim_{k\to \infty} \left(R_1^k,\dots, R_r^k, R_{r+1}^k \right)
=\lim_{k\to \infty}
\left(Q_1\Omega_1^kQ_1^\T,\dots,Q_r\Omega_r^kQ_r^\T,
 Q'\Omega_{r+1}^kQ''^\T \right)\,.
\end{equation}
Similarly to the proofs to {\bf  Cases 1-8} in the first part, by using   (\ref{eq:lim-g'-phi'}) and \eqref{eq:limit-deri-theta-1}--\eqref{eq:limit-deri-theta-3} in Lemma \ref{lem:limit-deri-theta}, we can show that \eqref{eq:lim-R-Omega} holds. For simplicity, we omit the details here. Therefore, we obtain that $\partial_{B}G(\overline{X})=\partial_{B}\Psi(0)$. This completes the proof. $\hfill \Box$

\section{Extensions}\label{section:extension}

In this section, we consider the spectral operators defined on the Cartesian product
of several  real or complex matrices. The corresponding properties, including   continuity,   directional differentiability,   (continuous) differentiability,   locally Lipschitzian continuity,   $\rho$-order B-differentiability,   $\rho$-order G-semismoothness and the characterization of Clarke's generalized Jacobian, can be studied in the same fashion as those in Section \ref{section:spectraloperator} and Section \ref{section:semismoothness}. Instead of presenting the    proofs here, we refer the readers to the PhD thesis of Ding \cite{Ding12} to work out details.

Without loss of generality, from now on, we assume that ${\cal X}=\S^{m_1}\times\V^{m_2\times n_2}$ with $m=m_1+m_2$. For any $X=\left(X_1,X_2\right)\in\S^{m_1}\times\V^{m_2\times n_2}$, denote $\kappa(X)=\left(\lambda(X_1),\sigma(X_2)\right)$. Let ${\cal N}$ be a given nonempty open set in ${\cal X}$. Suppose that ${\bf g}:\Re^m\to\Re^m$ is mixed symmetric, with respect to ${\cal P}\equiv{\mathbb P}^{m_1}\times\pm{\mathbb P}^{m_2}$, on an open set $\hat{\kappa}_{\cal N}$ in $\Re^m$ containing $\kappa_{\cal N}=\left\{\kappa(X)\mid X\in{\cal N}\right\}$. Let $G:{\cal X}\to{\cal X}$ be the corresponding spectral operator defined in Definition \ref{def:def-spectral-op}.

Let $\overline{X}=(\overline{X}_1,\overline{X}_2)\in{\cal N}$ be given. Suppose the given $\overline{X}_1\in\S^{m_1}$ and $\overline{X}_2\in\V^{m_2\times n_2}$ have the following decompositions
\begin{equation}\label{eq:decomp-X1-X2}
\overline{X}_1=\overline{P}{\rm Diag}(\lambda(\overline{X}_1))\overline{P}^\T \quad {\rm and} \quad \overline{X}_2=\overline{U}[{\rm Diag}(\sigma(\overline{X}_2))\quad 0]\overline{V}^\T \,,
\end{equation}
where $\overline{P}\in{\mathbb O}^{m_1}$, $\overline{U}\in{\mathbb O}^{m_2}$ and $\overline{V}=\left[ \overline{V}_{1} \quad  \overline{V}_{2} \right] \in{\mathbb O}^{n_2}$ with $\overline{V}_{1}\in\V^{n_2\times m_2}$ and $\overline{V}_{2}\in\V^{n_2\times (n_2-m_2)}$. Denote $\overline{\lambda}:=\lambda(\overline{X}_1)$, $\overline{\sigma}:=\sigma(\overline{X}_2)$ and $\overline{\kappa}:=\left(\overline{\lambda},\overline{\sigma}\right)$. We use $\overline{\nu}_{1}>\ldots>\overline{\nu}_{r_1}$ to denote the distinct eigenvalues of $\overline{X}_1$ and $\overline{\nu}_{r_1+1}>\ldots>\overline{\nu}_{r_1+r_2}>0$ to denote the distinct nonzero singular values of $\overline{X}_2$. Define the index sets
\begin{equation*}\label{eq:ak-W1-W2}
\left\{\begin{array}{ll}
a_{l}:=\{i\,|\,\overline{\lambda}_{i}=\overline{\nu}_{l}, \ 1\le i\le m_1\} & l=1,\ldots,r_1\,,\\[3pt]
a_{l}:=\{i\,|\,\overline{\sigma}_{i}=\overline{\nu}_{l}, \ 1\le i\le m_2\} & l=r_1+1,\ldots,r_1+r_2\,.
\end{array}\right.
\end{equation*}
Define $b:=\{i\,|\,\overline{\sigma}_{i}=0, \ 1\le i\le m_2\}$. We have the following result on the continuity of spectral operators.

\begin{theorem}
Let  $\overline{X}=(\overline{X}_1,\overline{X}_2)\in{\cal N}$ be given. Suppose that $\overline{X}_1$ and $\overline{X}_2$ have the decompositions \eqref{eq:decomp-X1-X2}. The spectral operator $G$ is continuous at $\overline{X}$ if and only if  ${\bf g}$ is continuous at $\kappa(\overline{X})$.
\end{theorem}

In order to present the results on the directional differentiability of spectral operators of matrices, we introduce some notations. For the given mixed symmetric mapping ${\bf g} = ({\bf g}_1, {\bf g}_2):\Re^m\to\Re^{m_1}\times \Re^{m_2}$, define the matrices ${\cal A}^0(\overline{\kappa})\in\S^{m_1}$, ${\cal E}^0_1(\overline{\kappa})\in\S^{m_2}$, ${\cal E}^0_2(\overline{\kappa})\in\V^{m_2\times m_2}$ and ${\cal F}^0(\overline{\kappa})\in\V^{m_2\times(n_2-m_2)}$ with respect to $\overline{\kappa}=(\overline{\lambda},\overline{\sigma})$ by
\[
({\cal A}^0(\overline{\kappa}))_{ij}:=\left\{ \begin{array}{ll}  \displaystyle{\frac{({\bf g}_{1}(\overline{\kappa}))_{i}-({\bf g}_{1}(\overline{\kappa}))_{j}}{\overline{\lambda}_{i}-\overline{\lambda}_{j}}} & \mbox{if $\overline{\lambda}_{i}\neq\overline{\lambda}_{j}$,}\\
0 & \mbox{otherwise,}
\end{array} \right.  \quad  i,j\in\{1,\ldots,m_{1}\}\,,
\]
\[
({\cal E}^0_{1}(\overline{\kappa}))_{ij}:=\left\{ \begin{array}{ll}  \displaystyle{\frac{({\bf g}_{2}(\overline{\kappa}))_{i}-({\bf g}_{2}(\overline{\kappa}))_{j}}{\overline{\sigma}_i-\overline{\sigma}_{j}}} & \mbox{if $\overline{\sigma}_i\neq\overline{\sigma}_j$}\,,\\
0 & \mbox{otherwise}\,,
\end{array} \right.  \quad i,j\in\{1,\ldots,m_2\}\,,
\]
\[
({\cal E}^0_{2}(\overline{\kappa}))_{ij}:=\left\{ \begin{array}{ll}  \displaystyle{\frac{({\bf g}_{2}(\overline{\kappa}))_{i}+({\bf g}_{2}(\overline{\kappa}))_{j}}{\overline{\sigma}_i+\overline{\sigma}_{j}}} & \mbox{if $\overline{\sigma}_i+\overline{\sigma}_j\neq 0$}\,,\\
0 & \mbox{otherwise}\,,
\end{array} \right.  \quad i,j\in\{1,\ldots,m_2\}
\] and
\[
({\cal F}^0(\overline{\kappa}))_{ij}:=\left\{\begin{array}{ll}
\displaystyle{\frac{({\bf g}_{2}(\overline{\kappa}))_{i}}{\overline{\sigma}_i}} & \mbox{if $\overline{\sigma}_i\neq 0$}\,,\\
0 & \mbox{otherwise.}
\end{array}\right. \quad  i\in\{1,\ldots,m_2\},\quad j\in\{1,\ldots, n_2-m_2\}\,.
\]

Suppose that ${\bf g}$ is directionally differentiable at $\overline{\kappa}$. Then, we know that the directional derivative $ {\bf g}'(\overline{\kappa};\cdot)=\left({\bf g}'_1(\overline{\kappa};\cdot),{\bf g}'_2(\overline{\kappa};\cdot)\right):\Re^{m_1+m_2}\to\Re^{m_1+m_2}$ satisfies that for any $(Q_1,Q_2)\in{\cal P}_{\overline{\kappa}}$ and any $({\bf h}_1,{\bf h}_2)\in\Re^{m_1}\times\Re^{m_2}$,
\begin{equation}\label{eq:phi_symmtric_2_block}
\Big({\bf g}'_1(\overline{\kappa};(Q_1{\bf h}_1,Q_2{\bf h}_2)),{\bf g}'_2(\overline{\kappa};(Q_1{\bf h}_1,Q_2{\bf h}_2))\Big)=\Big(Q_1{\bf g}'_1(\overline{\kappa};({\bf h}_1,{\bf h}_2)),Q_2{\bf g}'_2(\overline{\kappa};({\bf h}_1,{\bf h}_2))\Big)\,,
\end{equation}
where ${\cal P}_{\overline{\kappa}}$ is the subset of ${\cal P}\equiv{\mathbb P}^{m_1}\times\pm{\mathbb P}^{m_2}$ defined with respect to $\overline{\kappa}$ by
\[
{\cal P}_{\overline{\kappa}}:=\left\{\left(Q_1,Q_2\right)\in{\mathbb P}^{m_1}\times\pm{\mathbb P}^{m_2}\mid(\overline{\lambda},\overline{\sigma})=(Q_1\overline{\lambda},Q_2\overline{\sigma}) \right\}\,.
\] It is easy to check that $\left(Q_1,Q_2\right)\in{\cal P}_{\overline{\kappa}}$ if and only if there exist $Q_1^{l}\in{\mathbb P}^{|a_l|}$, $l=1,\ldots,r_1$, $Q_2^{l}\in{\mathbb P}^{|a_l|}$, $l=r_{1}+1,\ldots,r_1+r_{2}$ and $Q_{2}^{r_1+r_2+1}\in\pm{\mathbb P}^{|b|}$ such that
\begin{equation}\label{eq:Q-form-2-block}
Q_1={\rm Diag}\left(Q^{1}_1,\dots,Q^{r_1}_{1}\right)\in{\mathbb P}^{m_1} \quad {\rm and} \quad Q_2={\rm Diag}\left(Q_{2}^{r_{1}+1},\dots,Q_{2}^{r_{1}+r_{2}},Q^{r_1+r_2+1}_{2}\right)\in\pm{\mathbb P}^{m_2}\,.
\end{equation}
Denote $\phi(\cdot):={\bf g}'(\overline{\kappa};\cdot)$. For any $h\in\Re^{m}$, rewrite $\phi(h)\in\Re^m$  as
$\phi(h)=\left(\phi_{1}(h),\ldots,\phi_{r_{1}+r_2}(h),\phi_{r_{1}+r_2+1}(h)\right)$ with $\phi_{l}(h)\in \Re^{|a_l|}$ for $l=1,\ldots,r_1+r_2$ and
$\phi_{r_1+r_2+1}(h)\in \Re^{|b|}$. Therefore, we know from  \eqref{eq:phi_symmtric_2_block} and \eqref{eq:Q-form-2-block} that the directional derivative $\phi$ is mixed symmetric mapping, with respect to ${\mathbb P}^{|a_1|}\times\ldots\times{\mathbb P}^{|a_{r_1+r_2}|}\times\pm{\mathbb P}^{|b|}$, over $\Re^{|a_1|}\times\ldots\times\Re^{|a_{r_1+r_2}|}\times\Re^{|b|}$.
Denote \[{\cal W}:=\S^{|a_{1}|}\times\ldots\times\S^{|a_{r_1+r_2}|}\times\V^{|b|\times(|b|+n_2-m_2)}.\]
 Let $\Phi:{\cal W}\to{\cal W}$ be the corresponding spectral operator defined in Definition \ref{def:def-spectral-op} with respect to the mixed symmetric mapping $\phi$, i.e., for any $W=\left(W_{1}, \ldots,W_{r_{1}+r_2}, W_{r_{1}+r_2+1}\right)\in{\cal W}$,  \[\Phi(W) =\big(\Phi_{1}(W),\ldots, \Phi_{r_{1}+r_2}(W), \Phi_{r_{1}+r_2+1}(W) \big)\]
  with
\[
\Phi_{l}(W)=\left\{ \begin{array}{ll}
\widetilde{R}_{l}{\rm Diag}(\phi_{l}(\kappa(W)))\widetilde{R}_{l}^\T  & \mbox{if $l=1,\ldots,r_{1}+r_2$,}
\\[3pt]
\widetilde{M} {\rm Diag}(\phi_{r_{1}+r_2+1}(\kappa(W))) \widetilde{N}_{1}^\T  & \mbox{if $l=r_{1}+r_2+1$,}
\end{array}\right.
\]
where $\kappa(W)=\left(\lambda(W_{1}),\ldots,\lambda(W_{r_{1}+r_2}),\sigma(W_{r_{1}+r_2+1})\right)\in\Re^{m}$,
$\widetilde{R}_{l}\in{\mathbb O}^{|a_{l}|}(W_{l})$, and $(\widetilde{M},\widetilde{N})\in{\mathbb O}^{|b|,|b|+n_2-m_2}(W_{r_{1}+r_2+1})$, $\widetilde{N}=\big[\widetilde{N}_{1}\quad \widetilde{N}_{2}\big]$ with $\widetilde{N}_{1}\in\V^{(|b|+n_2-m_2)\times|b|}$, $\widetilde{N}_{2}\in\V^{(|b|+n_2-m_2)\times(n_2-m_2)}$. Then, the first divided directional difference ${\bf g}^{[1]}(\overline{X};H)\in{\cal X}$ of ${\bf g}$ at $\overline{X}$ along the direction $H=(H_1, H_2)\in{\cal X}$ is defined by
\[
{\bf g}^{[1]}(\overline{X};H):=\left({\bf g}_{1}^{[1]}(\overline{X};H), \; {\bf g}_{2}^{[1]}(\overline{X};H)\right)
\] with
\begin{eqnarray*}
{\bf g}_{1}^{[1]}(\overline{X};H) &= & {\cal A}^0(\overline{\kappa})\circ\overline{P}^\T H_1\overline{P}+
{\rm Diag}\Big(\Phi_{1}(D(H)) ,\dots, \Phi_{r_{1}}(D(H)) \Big)
\in\S^{m_{1}},
\\[3pt]
{\bf g}_{2}^{[1]}(\overline{X};H)&=&\left[{\cal E}^0_1(\overline{\kappa})\circ S(\overline{U}^\T H_2\overline{V}_{1})+{\cal E}^0_2(\overline{\kappa})\circ T(\overline{U}^\T H_2\overline{V}_{1})\quad {\cal F}^0(\overline{\kappa})\circ \overline{U}^\T H_2\overline{V}_{2}  \right]\\
&&+ \; \left[\begin{array}{cc} {\rm Diag}\left(\Phi_{r_{1}+1}(D(H)),\dots, \Phi_{r_{1}+r_2}(D(H))\right) & 0 \\[2mm]
0 &    \Phi_{r_{1}+r_2+1}(D(H))
\end{array} \right]
\in\V^{m_2\times n_2}\,,
\end{eqnarray*}
where
\begin{equation*}
D(H)=\left(\overline{P}_{a_1}^\T H_1\overline{P}_{a_1},\ldots,\overline{P}_{a_{r_1}}^\T H_1\overline{P}_{a_{r_1}}, S(\overline{U}_{a_{r_1+1}}^\T H_2\overline{V}_{a_{r_1+1}}),\ldots,S(\overline{U}_{a_{r_1+r_2}}^\T H_2\overline{V}_{a_{r_1+r_2}}),\overline{U}_{b}^\T H_2[\overline{V}_{b}\quad \overline{V}_{2}] \right)\in{\cal W}\,.
\end{equation*}

Now, we are ready to state the results on the directional differentiability of the spectral operator $G$.

\begin{theorem}\label{prop:H-dir-diff-spectr-op-block}
Let $\overline{X}=(\overline{X}_1,\overline{X}_2)\in{\cal N}$ be given. Suppose that $\overline{X}_1$ and $\overline{X}_2$ have the decompositions \eqref{eq:decomp-X1-X2}. The spectral operator $G$ is Hadamard directionally differentiable at $\overline{X}$ if and only if ${\bf g}$ is Hadamard directionally differentiable at $\kappa(\overline{X})$. In that case, $G$ is directionally differentiable at $\overline{X}$ and the directional derivative at $\overline{X}$ along any direction $H\in{\cal X}$ is given by
\begin{equation*}\label{eq:dir-diff-spectr-op-block}
G'(\overline{X};H)=\left({\overline{P}}{\bf g}_{1}^{[1]}(\overline{X};H)\overline{P}^\T , \; \overline{U}{\bf g}_{2}^{[1]}(\overline{X};H) \overline{V}^\T \right)\,.
\end{equation*}
\end{theorem}

In order to present the derivative formulas of spectral operators, we introduce the following notations. For the given $\overline{X}=(\overline{X}_1,\overline{X}_2)\in{\cal N}$, suppose that  ${\bf g}$ is F-differentiable at  $\overline{\kappa}=\kappa(\overline{X})$. Denote by ${\bf g}'(\overline{\kappa})\in\Re^{m\times m}$ the Jacobian matrix of ${\bf g}$ at $\overline{\kappa}$. Let $\eta_1(\overline{\kappa})\in\Re^{m_1}$ and $\eta_2(\overline{\kappa})\in\Re^{m_2}$ be the vectors defined by
\[
\left(\eta_1(\overline{\kappa})\right)_i:=\left\{ \begin{array}{ll} ({\bf g}_1'(\overline{\kappa}))_{ii}-({\bf g}_1'(\overline{\kappa}))_{i(i+1)} & \mbox{if $\exists\, j\in\left\{1,\ldots,m_1\right\}$ and $j\neq i$ such that $\overline{\lambda}_i=\overline{\lambda}_j$},\\
({\bf g}_1'(\overline{\kappa}))_{ii} & \mbox{otherwise}\,,
\end{array} \right. \quad i\in\{1,\ldots,m_1\}
\]
and
\[
\left(\eta_2(\overline{\kappa})\right)_i:=\left\{ \begin{array}{ll} ({\bf g}_2'(\overline{\kappa}))_{ii}-({\bf g}_2'(\overline{\kappa}))_{i(i+1)} & \mbox{if $\exists\, j\in\left\{1,\ldots,m_2\right\}$ and $j\neq i$ such that $\overline{\sigma}_i=\overline{\sigma}_j$},\\
({\bf g}_2'(\overline{\kappa}))_{ii} & \mbox{otherwise}\,,
\end{array} \right. \quad i\in\{1,\ldots,m_2\} \,.
\]
Define the corresponding {\it divided difference matrices} ${\cal A}(\overline{\kappa})\in\Re^{m_1\times m_1}$ and  ${\cal E}_1(\overline{\kappa})\in\Re^{m_2\times m_2}$, the {\it divided addition matrix} ${\cal E}_2(\overline{\kappa})\in\Re^{m_2\times m_2}$, the {\it division  matrix} ${\cal F}(\overline{\kappa})\in\Re^{m_2\times(n_2-m_2)}$, respectively, by
\[
({\cal A}(\overline{\kappa}))_{ij}:=\left\{ \begin{array}{ll}  \displaystyle{\frac{({\bf g}_{1}(\overline{\kappa}))_{i}-({\bf g}_{1}(\overline{\kappa}))_{j}}{\overline{\lambda}_{i}-\overline{\lambda}_{j}}} & \mbox{if $\overline{\lambda}_{i}\neq\overline{\lambda}_{j}$,}\\
\left(\eta_1(\overline{\kappa})\right)_i & \mbox{otherwise,}
\end{array} \right.  \quad  i,j\in\{1,\ldots,m_{1}\}\,,
\]
\[
({\cal E}_1(\overline{\kappa}))_{ij}:=\left\{ \begin{array}{ll}  \displaystyle{\frac{({\bf g}_{2}(\overline{\kappa}))_{i}-({\bf g}_{2}(\overline{\kappa}))_{j}}{\overline{\sigma}_{i}-\overline{\sigma}_{j}}} & \mbox{if $\overline{\sigma}_{i}\neq\overline{\sigma}_{j}$,}\\
\left(\eta_2(\overline{\kappa})\right)_i & \mbox{otherwise,}
\end{array} \right.  \quad  i,j\in\{1,\ldots,m_{2}\}\,,
\]
\[
({\cal E}_2(\overline{\kappa}))_{ij}:=\left\{ \begin{array}{ll}  \displaystyle{\frac{({\bf g}_{2}(\overline{\kappa}))_{i}+({\bf g}_{2}(\overline{\kappa}))_{j}}{\overline{\sigma}_{i}+\overline{\sigma}_{j}}} & \mbox{if $\overline{\sigma}_{i}+\overline{\sigma}_{j}\neq 0$,}\\
({\bf g}'_{2}(\overline{\kappa}))_{ii} & \mbox{otherwise,}
\end{array} \right.  \quad  i,j\in\{1,\ldots,m_{2}\}\,,
\]
\[
({\cal F}(\overline{\kappa}))_{ij}:=\left\{ \begin{array}{ll}  \displaystyle{\frac{({\bf g}_{2}(\overline{\kappa}))_{i}}{\overline{\sigma}_{i}}} & \mbox{if $\overline{\sigma}_{i}\neq 0$,}\\
({\bf g}'_{2}(\overline{\kappa}))_{ii} & \mbox{otherwise,}
\end{array} \right.  \quad  i\in\{1,\ldots,m_2\},\quad j\in\{1,\ldots, n_2-m_2\}.
\]
Define the  matrices ${\cal C}_1(\overline{\kappa})\in\Re^{m_1\times m}$ and ${\cal C}_2(\overline{\kappa})\in\Re^{m_2\times m}$ by
\[
{\cal C}_1(\overline{\kappa})={\bf g}_1'(\overline{\kappa})-\left[{\rm Diag}\left(\eta_1(\overline{\kappa} )\right) \quad 0 \right]\quad {\rm and} \quad {\cal C}_2(\overline{\kappa})={\bf g}_2'(\overline{\kappa})-\left[0 \quad {\rm Diag}\left(\eta_2(\overline{\kappa} )\right) \right]\,.
\]
Then, we have the following results on the F-differentiability of spectral operators.
\begin{theorem}\label{thm:diff-spectral-op-block}
Let $\overline{X}=(\overline{X}_1,\overline{X}_2)\in{\cal N}$ be given. Suppose that $\overline{X}_1$ and $\overline{X}_2$ have the decompositions \eqref{eq:decomp-X1-X2}. The spectral operator $G$ is (continuously) differentiable at $\overline{X}$ if and only if ${\bf g}$ is (continuously) differentiable at $\overline{\kappa}=\kappa(\overline{X})$. In that case, the derivative of $G$ at $\overline{X}$ is given by for any $H=(H_1,H_2)\in{\cal X}$,
\begin{equation}\label{eq:F-derivative-spectral-op-Phi}
\begin{array}{l}
G'(\overline{X})(H)
 =\left(\overline{P} [{\cal A}(\overline{\kappa})\circ \overline{P}^\T H_1\overline{P} +{\rm Diag}\left({\cal C}_1(\overline{\kappa})h\right)]\overline{P}^\T ,\right.\\
\left. \overline{U}\left[{\cal E}_{1}(\overline{\kappa})\circ S(\overline{U}^\T H_2\overline{V}_1)+{\rm Diag}\left({\cal C}_2(\overline{\kappa})h\right)+{\cal E}_{2}(\overline{\kappa})\circ T(\overline{U}^\T H_2\overline{V}_1) \quad {\cal F}(\overline{\kappa})\circ \overline{U}^\T H_2\overline{V}_2 \right]\overline{V}^\T \right)\,,
\end{array}
\end{equation}
where $h:=\big({\rm diag}\big(\overline{P}^\T H_1\overline{P}\big),{\rm diag}\big(S(\overline{U}^\T H_2\overline{V}_1)\big)\big)\in\Re^{m}$.
\end{theorem}

The following theorem is on the locally Lipschitzian continuity of spectral operators of matrices.
\begin{theorem}\label{thm:Lip-spectral-op-block}
Let $\overline{X}=(\overline{X}_1,\overline{X}_2)\in{\cal N}$ be given. Suppose that $\overline{X}_1$ and $\overline{X}_2$ have the decompositions \eqref{eq:decomp-X1-X2}. Then, the spectral operator $G$ is locally Lipschitz continuous near $\overline{X}$ if and only if ${\bf g}$ is locally Lipschitz continuous near $\kappa(\overline{X})$.
\end{theorem}

For the $\rho$-order B(ouligand)-differentiability of spectral operators, we have the following theorem.
\begin{theorem}\label{thm:rho-order-B-diff-spectral-op-block}
Let $\overline{X}=(\overline{X}_1,\overline{X}_2)\in{\cal N}$ be given. Suppose that $\overline{X}_1$ and $\overline{X}_2$ have the decompositions \eqref{eq:decomp-X1-X2}. Let $0<\rho\leq 1$ be given. Then,
\begin{itemize}
\item[(i)] if $\bf g$ is locally Lipschitz continuous near $\kappa(\overline{X})$ and $\rho$-order B-differentiable at $\kappa(\overline{X})$, then $G$ is $\rho$-order B-differentiable at $\overline{X}$;
\item[(ii)] if $G$ is $\rho$-order B-differentiable at $\overline{X}$, then $\bf g$ is $\rho$-order B-differentiable at $\kappa(\overline{X})$.
\end{itemize}
\end{theorem}

Suppose that ${\bf g}$ is locally Lipschitz continuous near $\kappa(\overline{X})$. Thus, we know from Theorem \ref{thm:Lip-spectral-op-block} that the corresponding spectral operator $G$ is also locally Lipschitz continuous near $\overline{X}$. We have the following theorem on the G-semismoothness of spectral operators.

 \begin{theorem}\label{thm:rho-order-G-semismooth-spectral-op-block}
Let $\overline{X}=(\overline{X}_1,\overline{X}_2)\in{\cal N}$ be given. Suppose that $\overline{X}_1$ and $\overline{X}_2$ have the decompositions \eqref{eq:decomp-X1-X2}. Let $0<\rho\leq 1$ be given. Then, $G$ is $\rho$-order G-semismooth at $\overline{X}$ if and only if  ${\bf g}$ is $\rho$-order G-semismooth at $\kappa(\overline{X})$.
\end{theorem}

Finally, we assume that ${\bf g}$ is locally Lipschitz continuous near $\overline{\kappa}=\kappa(\overline{X})$ and directionally differentiable at $\overline{\kappa}$. From Theorems
\ref{prop:H-dir-diff-spectr-op-block} and Theorem \ref{thm:Lip-spectral-op-block}, the spectral operator $G$ is also locally Lipschitz continuous near $\overline{X}$ and directionally differentiable at $\overline{X}$. Then, we have the following results on the characterization of the B-subdifferential $\partial_B G(\overline{X})$ and Clarke's subdifferential $\partial G(\overline{X})$.

\begin{theorem}\label{thm:B-subdiff-spectral-op-block}
Let $\overline{X}=(\overline{X}_1,\overline{X}_2)\in{\cal N}$ be given. Suppose that $\overline{X}_1$ and $\overline{X}_2$ have the decompositions \eqref{eq:decomp-X1-X2}. Suppose that there exists an open neighborhood ${\cal B}\subseteq\Re^{m}$ of $\overline{\kappa}$ in $\hat{\kappa}_{\cal N}$ such that ${\bf g}(\cdot)$ is differentiable at $\kappa\in{\cal B}$ if and only if $\phi={\bf g}'(\overline{\kappa};\cdot)$ is differentiable at $\kappa-\overline{\kappa}$. Assume that the function $d:\Re^{m}\to\Re^{m}$ defined by
\[
d(h)={\bf g}(\overline{\kappa}+h)-{\bf g}(\overline{\kappa})-{\bf g}'(\overline{\kappa};h),\quad h\in\Re^{m}
\] is strictly differentiable at zero. Then, we have
\[
\partial_{B}G(\overline{X})=\partial_{B}\Psi(0)\quad {\rm and} \quad \partial G(\overline{X})=\partial\Psi(0)\,,
\]where $\Psi(\cdot):=G'(\overline{X};\cdot):{\cal X}\to{\cal X}$ is the directional derivative of $G$ at $\overline{X}$.
\end{theorem}

\section{Conclusions}
In this paper, we introduced and studied a class of matrix-valued functions, termed spectral operators of matrices, which frequently arise and play a crucial role in various applications including matrix optimization problems, matrix completion, multi-dimensional shape analysis and others. Several fundamental properties of spectral operators, including   well-definedness, continuity,   directional differentiability,  Fr\'{e}chet-differentiability,  locally Lipschitz continuity,   $\rho$-order B(ouligand)-differentiability ($0<\rho\leq 1$),   $\rho$-order G-semismooth ($0<\rho\leq 1$) and the characterization of Clarke's generalized Jacobian, are studied systematically. These results   provide the necessary theoretical foundations 
{for} many applications. Consequently, one is able to use these results to design some efficient numerical methods for solving large-scale matrix optimization problems arising from many important applications. For instance, Chen et al. \cite{CLSToh12} proposed an efficient and robust semismooth Newton-CG dual proximal point algorithm for solving large scale matrix spectral norm approximation problems. In \cite{CLSToh12}, the properties of the spectral operator, such as the semismoothness and the characterization of Clarke's generalized Jacobian, {played} an important role in the convergence analysis of the proposed algorithm. The work done in this paper on spectral operators of matrices is by no means complete. Due to the
{rapid advances in the applications of matrix optimization in different fields,}
spectral operators of matrices will become even more important and many other properties of spectral operators are waiting to be explored.

\end{document}